\documentclass[12pt,a4paper,twoside,chapter=TITLE,section=TITLE]{abntex2}%{book}
\usepackage[tablename=Table]{caption}
\usepackage{relsize}
\usepackage[utf8]{inputenc}
\usepackage{amsmath,amssymb}
\usepackage{amsthm}
\usepackage[all]{xy}
\newtheorem{remk}{Remark}
\newtheorem{lema}{Lemma}
\newtheorem{prop}{Proposition}
\newtheorem{teo}{Theorem}
\newtheorem{defi}{Definition}
\newtheorem{exem}{Example}
\usepackage{mathtools}
\usepackage[active]{srcltx}
\usepackage{hyperref,ifpdf}
\usepackage{graphicx}
\date{May, 2018}
\DeclareMathOperator{\sym}{Sym}
\settocdepth{subsection}
\hypersetup{
	colorlinks=true,
	linkcolor=red,
	filecolor=magenta,      
	urlcolor=cyan,
}
\def\rk{\ensuremath{\operatorname{rank}}}
\def\W{\ensuremath{\mathcal W}}

\def\cl#1{\ensuremath{\mathcal{#1}}}
\def\p#1{\ensuremath{\mathbb P^{#1}}}
\def\pd#1{\ensuremath{\check{\vphantom{|}\mathbb P}^{#1}}}
\def\bb#1{\ensuremath{\mathbb{#1}}}
\def\rar{\dasharrow}
\def\cao{\c c\~ao }
\def\co{\c co }
\def\id#1{\ensuremath{\langle #1\rangle}}
\def\red#1{{\color{red}#1}}
\def\blu#1{{\color{blue}#1}}
\def\ba#1{\begin{array}{#1}}\def\ea{\end{array}}
\def\na#1{\noalign{\vskip#1pt}}\def\ea{\end{array}}
\def\be{\begin{equation}}\def\ee{\end{equation}}

\newcommand{\hooklongrightarrow}{\lhook\joinrel\longrightarrow}
\def\twoheadlar{\twoheadlongrightarrow}
\newcommand{\twoheadlongrightarrow}{\relbar\joinrel{\!\!\twoheadrightarrow}}
\def\hooklar{\lhook\joinrel\longrightarrow}
\def\X{\ensuremath{\bb X}}
\def\vez{\ensuremath{\times}}

\def\ov#1{\ensuremath{\overline{#1}}}
\def\un#1{\ensuremath{\underline{#1}}}
\def\mb#1{\ensuremath{\mathbf{#1}}}
\def\ve{\ensuremath{^\vee}}
\def\ffi{\ensuremath{\varphi_0}}
\def\sym{\ensuremath{\operatorname{Sym}}}

\newcommand{\hookdownarrow}{
\ifpdf
\mathrel{\rotatebox[origin=c]{-90}{$\hookrightarrow$}}\else
\ba c \ve\\\na{-16}
\downarrow\ea\fi
}

\begin{document}
\title{Degree of the Exceptional Component of 
the Space of Holomorphic
Foliations of 
\\
Degree Two and Codimension One in 
$\mathbb{P}^3$}
\author{Artur Afonso Guedes Rossini}
\orientador{Israel Vainsencher}
\preambulo{Tese apresentada ao Departamento de 
Matem\'atica
	da Universidade Federal de Minas Gerais, como
        requisito para a obten\c c\~ao
	do grau de doutor em Matem\'atica.}
\local{Belo Horizonte - MG}
\instituicao{UFMG}

\imprimircapa
\imprimirfolhaderosto

\chapter*{Agradecimentos}

\textcolor{white}{}

Ao professor e orientador Israel Vainsencher por todo incentivo, apoio, amizade e dedicação ao longo desse trabalho.

À minha esposa Aline por toda paciência e compreensão nesse período, pelo carinho de sempre e por todo amor que inclui em minha vida.

Aos meus pais Sergio e Laura por tudo, e em particular por sempre me incentivarem a chegar até aqui.

À professora Viviana Ferrer Cuadrado (UFF) pelas conversas e solicitações sempre atendidas.

Aos amigos que fizeram essa jornada ser ainda mais agradável, em especial aos grandes companheiros Pedro, Bruno, Gilberto e Divane.

Aos professores do núcleo de Matemática do IF Sudeste MG -- Campus Juiz de Fora pela amizade e o apoio para minha qualificação.

A todos que, de alguma forma, contribuíram para esse trabalho, meus agradecimentos.

\chapter*{Resumo}

\nocite{Abhyankar}
\nocite{AluffiFaber}
\nocite{jinvariante}
\nocite{Coray}
\nocite{CukAceaMassri}
\nocite{CukiermanVI}
\nocite{3264}
\nocite{Fulton}
\nocite{hartshorne2013algebraic}
\nocite{Ionescu}
\nocite{daniel}
\nocite{meirelles}
\nocite{alcides}

\noindent O prop\'osito deste trabalho \'e obter o grau
da componente excepcional, $E(3)$, do espa\c co de \linebreak
folhea\c
c\~oes holomorfas de grau 2 e codimens\~ao um em
$\mathbb{P}^3$.
Trata-se de uma componente de dimens\~ao treze,
descrita no c\'elebre trabalho de Alcides Lins-Neto e
Dominique Cerveau, \cite{alcides}. $E(3)$
\'e o fecho da \'orbita, sob a a\cao natural de
Aut\p3,
da folhea\cao
definida pela forma diferencial 
$$\omega=(3fdg-2gdf)/x_0, \text{\ 
\ onde
}\ \ f=x_0^2x_3-x_0x_1x_2+\dfrac{x_1^3}{3},\ \ 
g=x_0x_2-\dfrac{x_1^2}{2}.
$$

Nossa tarefa inicial é descrever uma caracteriza\cao
geométrica para o par $g,f$. Isto nos levará à
constru\cao de um espa\co de parâmetros como um fibrado
explícito sobre a variedade de bandeiras completas de \p3.
Fazendo uso de ferramentas
 de teoria da interse\c c\~ao
equivariante, podemos calcular o n\'umero desejado
como uma integral sobre o nosso espa\co de parâmetros.

\bigskip

\noindent \textbf{Palavras-chave:} Folhea\c c\~oes
Holomorfas.
 Componente Excepcional. Grau.

\chapter*{Abstract}

\noindent The purpose of this work is to obtain the
degree of the exceptional component, $E(3)$, of the space
of holomorphic foliations of degree two and codimension
one in $\mathbb{P}^3$. As shown in the celebrated work by
Dominique Cerveau and Alcides Lins Neto \cite{alcides},
$E(3)$ is a 13-dimensional component. It is the
closure of the orbit under the natural action of
 Aut\p3
of the foliation defined by the differential form
$$\omega=(3fdg-2gdf)\big/x_0,
\ \mbox{ where } \ f=x_0^2x_3-x_0x_1x_2+\dfrac{x_1^3}{3},\ 
g=x_0x_2-\dfrac{x_1^2}{2}
.$$
Our first task is to unravel a geometric characterization
of the pair $g,f$. This leads us to 
 the construction of a parameter space as 
an explicit fiber bundle over the 
variety of complete flags.
Using tools from equivariant intersection theory,
especially Bott's formula,
the  degree is expressed as an integral over our
parameter
space.

\bigskip

\noindent \textbf{Keywords:} Holomorphic Foliations. Exceptional Component. Degree.

\renewcommand{\listtablename}{List of Tables}
\renewcommand{\chaptername}{}
\renewcommand\tablename{Table}
\renewcommand{\proofname}{Proof}
\renewcommand{\contentsname}{Contents}
\renewcommand{\bibname}{Bibliography}
\renewcommand{\figurename}{Figure}

\cleardoublepage
\listoftables*

\cleardoublepage
\tableofcontents*
\textual
\chapter{Introduction}

A holomorphic foliation of codimension one and degree $d$
in the complex projective space $\mathbb{P}^n$ is given
by a differential 1-form 
$$\omega = A_0dx_0 + \cdots%A_1dx_1 + A_2dx_2 
+ A_ndx_n 
$$
where $A_0,\dots,A_n$ are homogeneous polynomials of degree $d+1$, satisfying the conditions
\begin{enumerate}
	\item (projectivity) $A_0x_0 + \cdots +A_nx_n = 0$
	\item (integrability) $\omega \wedge d\omega = 0.$
\end{enumerate}

These conditions define a closed scheme
$\overline{\mathcal{F}(d,n)}$ in the projective
space of global sections of the twisted cotangent sheaf
$\Omega_{\mathbb{P}^n}^1(d+2)$. 
It is further required that the singular locus of
$\omega$ (defined by the common zeros of the polynomials
$A_i$, $i=0,\dots,n$)  be of codimension greater than
or equal to two. This  
 corresponds to the inexistence of a common factor of
 positive degree for the polynomials $A_i$ thus defining
 a  Zariski open subset $\mathcal{F}(d,n)$. 
Naturally, the closure 
$\overline{\mathcal{F}(d,n)}$ possesses a
decomposition into irreducible components.

In \cite{alcides}, a full description of the
components of $\overline{\mathcal{F}(2,n)}
,n\geq3
$ is
given. We learn that there are six components of $\overline{\mathcal{F}(2,3)}$:

\begin{enumerate}
	\item Linear pull-backs $\cl S(2,3):$  For 
a degree 2 foliation $\mathcal{F}$ in $\mathbb{P}^2$ and a linear rational map $\alpha: \mathbb{P}^3 \dashrightarrow \mathbb{P}^2$, take the pull-back $\mathcal{F}^* = \alpha ^*(\mathcal{F})$.
	\item \label{R22}  Rational $\cl R(2,2):$ 
Take homogeneous polynomials $f,g$ both of degree 2 and
make the 1-form \ $\omega = fdg - gdf$.
	\item \label{R23} Rational $\cl R(1,3):$ Take
          homogeneous polynomials $f$ of degree 1 and $g$
          of degree 3, and 
form \ $\omega = fdg - 3gdf$.
	\item Logarithmic $\mathcal{L}(1,1,1,1):$ Take four linear polynomials $f_1 \cdots f_4$ and complex numbers $\lambda _1, \dots, \lambda _4$ with $\lambda_1+ \cdots +  \lambda_4 = 0$. Then, write
	$$\omega =
        \lambda_1f_2f_3f_4df_1+\lambda_2f_1f_3f_4df_2+\lambda_3f_1f_2f_4df_3+\lambda_4f_1f_2f_3df_4
.$$
	\item Logarithmic $\mathcal{L}(1,1,2):$ Take
          $f_1,f_2$ linear polynomials, $f_3$ quadratic,
          and complex numbers $\lambda _1, \dots, \lambda _3$ with $\lambda_1+ \lambda_2 + 2\lambda_3 = 0$. Then, write
$$\omega =
\lambda_1f_2f_3df_1+\lambda_2f_1f_3df_2+\lambda_3f_1f_2df_3.
$$

\item The exceptional component $E(3):$ This
 is the orbit closure of the
foliation defined by the differential form
\begin{equation}\label{fgi}
\omega=\frac{3f_0dg_0-2g_0df_0}{x_0},
\ \mbox{ where } \ f_0=x_0^2x_3-x_0x_1x_2+\dfrac{x_1^3}{3},\ 
g_0=x_0x_2-\dfrac{x_1^2}{2}
.\end{equation}

We will study in detail the component $E(3)$ along this work.
The main point is to describe the geometry of the family
of pairs $g,f$ in the orbit closure.
\end{enumerate}

Two important invariants of a projective scheme
are its dimension and degree.  

For some components of spaces of foliations the degree
has been found in recent works.
For instance, in \cite{CukiermanVI} the authors managed
to compute 
the degrees of the rational components of type $\cl
R(n,d,d)$ in $\mathbb{P}^n, \ n \geq 3, d>1$: 
$$\textrm{deg} \ \cl R(n,d,d) = \dfrac{1}{N_d-1}
\displaystyle{2N_d-2 \choose N_d},
$$
where $N_d = \displaystyle{n+d \choose d} - 1$. As a
particular case, the degree of the rational component
$\cl R(2,2)$ is 1430. Also the degrees of rational
components of type $\cl R(n,d_0,d_1)$ where $d_0$ divides
$d_1$ were found: 

$$\textrm{deg} \ \cl R(n,d_0,d_1) = \displaystyle{N_{d_1}+N_{d_0}-1 \choose N_{d_0}} - \dfrac{d_1}{d_0} \displaystyle{N_{d_1}+N_{d_0}-1 \choose N_{d_0}-1},$$
and so the degree of $\cl R(1,3)$ is 700.

The degrees of the rational components $\cl R(n,2,3)$ were also computed in \cite{CukiermanVI} for $n \leq 5$ and in \cite{daniel} this is extended to the degree of rational components $\cl R(n,2,2r+1)$ for general $n,r$.

The case of linear pullbacks was done by Viviane Ferrer and Israel Vainsencher (to appear). For linear pullbacks of foliations of degree $d$ in $\mathbb{P}^2$, the degree of the component $\cl S(d,3)$ is
$$\textrm{deg} \ \cl S(d,3) = \dfrac{1}{162}\dfrac{(d+4)!}{(d+1)!}(d^2+6d+11)(d^2+2d+3),$$
so we have $\textrm{deg} \ \cl S(2,3) = 1320.$

To the best of my knowledge, for the  logarithmic
components the degrees have not been determined. In  \cite{CukAceaMassri} the authors showed that these components are generically reduced. 

In this thesis we handle the construction of an adequate
parameter space for the exceptional component, enabling
us to achieve  the calculation of the degree of 
$$E(3)\subset\p{}(H^0(\Omega^1_{\p3}(4)))=\p{44}.
 $$

As customary in enumerative questions, the bulk of the
work lies in the description of an appropriate parameter
space. The description of
the exceptional component $E(3)$ could in principle suggest
using Aluffi \& Faber's approach to the calculation of
degrees of orbits closures, cf.\,\cite{AluffiFaber}. We took 
another path, essentially trying and deciphering an
algebro-geometric characterization of that pair
(cubic,quadric) from the
explicit expressions (\ref{fgi}) given just above.  Eventually, this lead to our construction
of an explicit smooth, projective variety $\bb Y_4$ 
(cf.\, Proposition \ref{Y4}, page \pageref{Y4}) endowed with a
generically injective morphism onto $E(3)$. 

Our strategy  for the exceptional component starts with a
complete flag in $\mathbb{P}^3$ 
$$p \textrm{ (point)} \in \ell \textrm{ (line) } 
\subset v \textrm { (plane)},
$$
over which we describe suitable cubic forms $f$ and 
quadratic forms $g$ that fulfill some special
conditions. The polynomials $g,f$ will 
give us the differential 1-form
$$\omega = \dfrac{3fdg-2gdf}{h},
$$
where $h$ is an equation of the plane $v$ of the
given flag. 

It turns out that $\omega$ can be zero for certain values
of otherwise acceptable $g,f$, for instance, $g=x_0^2, \  f=x_0^3$.
The main technical difficulty is how
 to handle the indeterminacy locus of the rational map
$$(g,f) \mapsto \omega
$$
leading to a resolution of its indeterminacies. As
we will see, this indeterminacy locus is a non-reduced and
reducible scheme.  The resolution of indeterminacies is done by a
careful analysis of the irreducible 
components of the  indeterminacy locus, blown up one at a time,

\medskip
$$\xymatrix%@R1pc@C1pc
{\bb Y_4\,\ar@/^2.0pc/[rrrrr]\,\ar@{->}[r]
&
\bb Y_3\,\ar@{->}[r]&\bb Y_2\,\ar@{->}[r]&\bb Y_1\,\ar@{->}[r]&\bb Y\,\ar@{-->}[r]&E(3)\subset\p{44}.
}$$

Once accomplished the resolution of the map, 
we obtain a line bundle $\cl W$ over $\bb Y_4$, pullback
of $\cl O_{\p{44}}(-1)$,
 with each fiber spanned by
a differential 1-form possibly in the boundary of $E(3)$.
% over the compactified space (simple by taking the
% well-defined $\omega$ with its multiple writings up to
% a scalar). 
Since the dimension of the exceptional component is 13, the required degree is (see Theorem \ref{integral})
$$\int -c_1^{13}(\W).$$

This will be done in Chapter \ref{chapdescricao}. In Chapter \ref{chapgrau}, we proceed to compute this integral. In order to do this, we will apply Bott's formula
$$\int -c_1^{13}(\mathcal{W}) =\sum\limits_{F}\dfrac{-c_1^T(\mathcal{W})^{13} \cap [F]_T}{c_{top}^T(\mathcal{N}_{F})},$$
where the sum runs through all  the fixed components of
a suitable action of the torus $T:=\mathbb{C}^*$ on $\bb Y_4$, induced
from an action on $\mathbb{P}^3$.
%de fato no espa\c co de par\^ametros constru\'\i do}, 
The $\mathcal{N}_F$ appearing in the denominator denotes
the normal bundle of a fixed component $F$ in $\bb Y_4$.

Putting together the contributions from
%After going through 
a total of 1728 fixed isolated points and 120 fixed lines we  find in Theorem \ref{GRAUE3} the number $\mathbf{168208}$ as the desired degree.

In the appendix we include the scripts for
\textit{Macaulay2} used
to study the resolution of the indeterminacies and to
compute the sum of all contributions in Bott's formula. 

\section{Notation and conventions}

 Let $x_0,x_1,x_2,x_3 $ denote homogeneous
 coordinates of $\mathbb{P}^3$.  

Write  
\\\centerline{$S_d=\id{x_0^d,x_0^{d-1}x_1,\dots,x_3^d} $}
for the vector space of  homogeneous polynomials of
 degree $d$. If no
confusion arises, we abuse and write $S_d$ for the trivial
vector bundle $S_d\vez X$ over some base variety
$X$. 

The
dual projective space $\pd3=\p{}(S_1)$ comes equipped
with a tautological line subbundle, 
\\\centerline{$\cl O_{\pd3}(-1)
{\hooklongrightarrow}S_1$.}

Given a vector bundle \cl V, we write $\p{}(\cl V)$
for the projective bundle of \rk\ one subspaces of the
fibers of \cl V. It comes endowed with a tautological
\rk\ one subline bundle
\\\centerline{$
\cl O_{\cl V}(-1){\hooklongrightarrow}\cl V,
$}
where we omit the pullback $\cl V_{|\p{}(\cl V)}$.

For an action of $\mathbb{C}^*$ in $X$, we denote
\  $t \circ v$ \
for the image of the point $(t,v) \in \mathbb{C}^* \times X$ by the action.

For clarity, we will use a left square ($\square$) also at the end of examples.%, like example \ref{calcl1}.

\chapter{Description of the Exceptional 
Component}  \label{chapdescricao}

\section{Associated complete flag}

 In \cite{alcides} the authors  describe
the exceptional component of foliations in $\mathbb{P}^3$
 starting with the two  homogeneous polynomials
\be\label{fg}
f_0:=x_0^2x_3-x_0x_1x_2+\dfrac{x_1^3}{3} \ \ 
\text{ and } \ \ 
g_0:=x_0x_2-\dfrac{x_1^2}{2}
\ee
from which arises the integrable differential 1-form
$$\omega_0 := \dfrac{2g_0df_0 - 3f_0dg_0}{x_0}.$$
Explicitly,
$$\ba{rcrcr}
\omega_0&=&
(x_1x_2^2-2x_1^2x_3+x_0x_2x_3)dx_0& +&
x_0(3x_1x_3-2x_2^2)dx_1  \\
&&\! +x_0(x_1x_2-3x_0x_3)dx_2& +&\!\!
x_0(2x_0x_2-x_1^2)dx_3.
\ea$$

The singular locus of  such  a  foliation consists
of a union of three curves\,(cf.\,appendix\,\ref{primdec}):

\begin{enumerate}
\item the conic given by the ideal     
$\langle x_0,x_2^2-2x_1x_3 \rangle$
(it lies in  the plane $x_0=0$);
\item the line $\ell_0$ defined by
          $ x_0=x_1=0$;
\item the twisted cubic given by $ \langle 2x_2^2-3x_1x_3, \ x_1x_2-3x_0x_3, \ x_1^2-2x_0x_2 \rangle$. 
\end{enumerate}  
Notice that these three components meet at the
point $p_0:=(0:0:0:1) \in \mathbb{P}^3$. The line
$\ell_0$ is tangent to the twisted cubic, the
plane is the osculating plane at $p_0$
and the conic is the
osculating conic.  The component $E(3)$ consists
of the orbit of $\omega_0$ under the group Aut\p3
of automorphism of $\mathbb{P}^3$; its dimension
is equal to 13 (see \cite{jinvariante}).

Let us examine the geometry of the  surface defined
by our cubic form (\ref{fg}) 
$f_0$.

We see  that it is an irreducible cubic, singular along the
line \,$\ell_0$.
Moreover, given a point $ p_t = (0: 0: t: 1) \in \ell_0
$, the tangent cone to our cubic surface at the point $
p_t $ has equation
$$ x_0^2 - tx_0x_1 = x_0 (x_0-tx_1). 
$$
Thus,  the tangent cone is a pair of planes containing the double
line \,$\ell_0$, one of which is fixed and the other varies
with the point $ p_t$.

We have also perceived that at  the special
point $ p_0= (0: 0: 0: 1) $ the tangent
cone to the cubic is the double plane $ x_0^2 = 0$.

Therefore, our cubic $f_0$ comes endowed with
a companion complete flag 
\begin{equation} \label{pdv}
  \varphi_0:\ 
p_0=\{x_0=x_1=x_2=0\} \in \ell_0 = \{x_0=x_1=0\} 
\subset v_0 =
\{x_0=0\} .
\ee Such a cubic  belongs to one of the strata of  cubics
singular along a  line, as described in
\cite{Abhyankar} and revisited in \cite{Coray}.

As for the quadric
$$ g_0 = x_0x_2 - \frac{x_1^2}{2},$$
notice we have a cone containing the line $\ell_0$
and vertex 
$p_0$. Moreover, the tangent plane to the cone at the
smooth points on $\ell_0$ is precisely the same plane $v_0$.

\begin{remk}\normalfont
\label{recupflag}  
The flag can move under the action of an
automorphism of $\mathbb{P}^3$, but it can be
recovered directly from the 1-form $\omega $
defining the exceptional foliation, just by
looking at the singular locus. For any such
$\omega$, the singular locus has three components
- a conic living on the new plane, the new line
and a twisted cubic that meets the other two
components just at the new point; the line is
tangent to the twisted cubic, the plane is the
osculating plane and the conic is the osculating
conic. By the way, the dimension 13 mentioned just
above is the dimension of the family of pointed
twisted cubics.
\end{remk}

These considerations about the cubic/quadric pairs
as in (\ref{fg}) lead us to consider the construction of a
parameter space for the family of exceptional foliations by starting with the variety $ {\bb F}$  of complete flags on $ \mathbb{P}^3 $, 
$$
p \textrm{(point)} \in \ell \textrm{(line)} \subset \pi
\textrm{(plane)}.
$$
Given such a flag we take: \vspace{1cm}

$\bullet$ A cubic surface $f$ with the properties:
\begin{enumerate}
	\item $f$ is singular along $\ell$;   
	\item the tangent cone to $f$ at
          each point $q \in \ell$ is
          the union of two planes $\pi \cup \alpha_q$.  
	\item $\alpha_p = \pi$, that is, at the point $p$ the tangent cone is the double plane.
\end{enumerate} \vspace{0,5cm}

$\bullet$ A quadratic cone $g$ with the properties:
\begin{enumerate}
	\item the line $\ell$ is contained in the cone $g$;
	\item the plane $\pi$ is tangent to $g$ along
the line $\ell$;
	\item the point $p$ is in the vertex of
          the cone.  
\end{enumerate} \vspace{0,5cm}

\begin{defi}\normalfont
	For simplicity we will call here a cubic or a quadric that meet the \linebreak requirements above as \textit{special}.
\end{defi}

An exceptional foliation is then given by the differential form
$$\omega = \dfrac{3fdg - 2gdf}{h},$$
where $h$ is an equation of the plane $\pi$. 

At this point, we realize that after the
construction of  special pairs
$(f, g)$, it is still necessary to
impose the condition that the differential form $
3fdg-2gdf $ be divisible by the equation of the plane. We
refer to this as  the \textit{divisibility condition}.
It deserves a closer look. 
In order to better understand this divisibility
condition, let us fix a complete flag as in
(\ref{pdv}). 

The conditions about the cubic surface show that its equation can be expressed as
\begin{equation}\label{ff}\left.\ba l
f=\left(a_0x_0 + a_1x_1 + a_2x_2+a_3x_3\right)x_0^2 +
a_4x_0x_1^2+a_5x_0x_1x_2 + a_6x_1^3,\\\na9
\text{\!\hskip-1.39 cm
	while the conditions about the quadric implies
}\\\na9
g=b_0x_0^2 + b_1x_0x_1 + b_2x_0x_2 + b_3x_1^2.
\ea\right\}
\end{equation}
We register for later use the following invariant description.

\begin{lema} \label{monsfg} 
	The $7$ monomials appearing in \ $f$
	\\\centerline{$x_{0}^{3},x_{0}^{2}x_{1},
		x_{0}^{2}x_{2}, x_{0}^{2}x_{3}, x_{0}x_{1}^{2},   x_{0}x_{1}x_{2}, x_{1}^{3} $}
 are the monomial generators of the
	subspace  of cubics,
	\\\centerline{$x_0 ^2\id{x_0,\dots,x_3}+x_0
		\id{x_0,x_1}\id{x_0,x_1,x_2}+\id{x_0,x_1}^3\subset S_3$.}
	Likewise, the $4$ monomials \ $x_0^2,x_0x_1,x_0x_2,x_1^2$ \ 
	in the quadric \ $g$ \ generate the subspace
	\\\centerline{$x_0\id{x_0,x_1,x_2} +  \id{x_0,x_1}^2\subset S_2
		$.}
	\qed
\end{lema}

\noindent
Varying the flag, we obtain equivariant vector
subbundles with respective ranks 7 and 4,
\be\label{AB}
\cl A\subset S_3 \ \text{ and } \ \cl B\subset S_2.
\ee

Projectively we obtain, at this stage, the variety
$$
\p{}(\cl A)\vez_{\bb F}\p{}(\cl B),
$$ a $\mathbb{P}^6 \times \mathbb{P}^3$ bundle
of pairs $(f,g)$ over a fixed flag. 

Continuing the discussion about the divisibility
condition, we may write
$$3fdg-2gdf = x_0w_0 + \left[(3a_6b_1-2a_4b_3)x_1^4 + (3a_6b_2 - 2a_5b_3)x_1^3x_2\right]dx_0,$$
where $w_0$ is a 1-form. From this it follows that the
divisibility condition is given, on the fiber over
\ffi,\ by the equations 
$$\left\{
\ba c
3a_6b_1 = 2a_4b_3, 
\\3a_6b_2 = 2a_5b_3 .
\ea\right.
$$
This locus consists of two irreducible components inside
$\mathbb{P}^6 \times \mathbb{P}^3$, both of
codimension two:
\begin{equation}
  \label{divcond1}
\ba{lrc}
&  (\star)& \ \ \ \
\left\{\begin{array}{lcr}
3a_6b_1 &=& 2a_4b_3\\
3a_6b_2 &=& 2a_5b_3\\
a_4b_2 &=& a_5b_1
\end{array}\right.
\\
\text{and} & (\star \star )& a_6=b_3=0.
\ea
\end{equation}
For a pair $(g,f)$ satisfying $(\star \star )$,
we actually get
$$\left\{\ba l
f=x_0\left(a_0x_0^2 + a_1x_0x_1 + a_2x_0x_2+a_3x_0x_3 +
a_4x_1^2+a_5x_1x_2\right),\\\na7
g=x_0\left(b_0x_0 + b_1x_1 + b_2x_2\right).
\ea\right.$$
It means that a general element in the second component $(\star \star)$
consists of a cubic and a quadric both
divisible by the equation of the plane,
certainly not interesting for our study of the 
exceptional component.

Henceforth, we refer to the divisibility condition
as  the three
equations  (\ref{divcond1}).

Recalling the dimension,  13, of the exceptional
component, it 
is reassuring to find 
$$\!\!
\ba{*8c}
{\text{\large{6}}}& \!\mbox{\large+}\! &{\text{\large{6}}}
&\!\mbox{\large+}\!&{\text{\large{3}}}
&\! \mbox{\large$-$}\! &{\text{\large{2}}}&\!\!=13.\\
\overbrace{{\mbox{dimension of }}}&&
\overbrace{{\mbox{$\mathbb{P}^6$ of special cubics}}}&&
% \ba c 
\overbrace{
	\p3 \text{ of special}}
&&
\overbrace{{\mbox{divisibility}}}&
\\\na{-4.2}
\text{complete flags}&&&&\text{quadric cones}&
&\text{condition}&
\ea$$

\section {Construction of a parameter space}
Let $\mathbb{G}=\mathbb{G}(1,3)$ be the Grassmann variety of lines in $\mathbb{P}^3$, with tautological sequence
$$\mathcal{S} \hooklongrightarrow \mathbb{G} \times \mathbb{C}^4 \twoheadlongrightarrow \mathcal{Q},$$
where \rk (\cl S)=2.
Denote by  $\mathcal{Q}^{\vee}$ the dual
of $\mathcal{Q}$. We also have the tautological sequence
over \p3,
$$\ba {*5c}
\cl O_{\p3}(-1)&
\hooklongrightarrow& \p3
\times \mathbb{C}^4& \twoheadlongrightarrow&\cl P
%\\\na8
%\cl O_{\pd3}(-1)&
%\hooklongrightarrow& \p3
%\times (\mathbb{C}^4)\ve &\twoheadlongrightarrow&\cl R
\ea
$$
Then,
$$
\ba {*5c}
\mathbb{P}(\mathcal{S})&=&\{(p,\ell) \,|\, p \in \ell
\} &\subset& \mathbb{P}^3 \times \mathbb{G}
\\\na6
\mathbb{P}(\mathcal{Q}^{\vee})&=&
\{(\ell, \pi) \,|\, \ell \subset \pi \}&
\subset& \mathbb{G} \times \pd3
.\ea$$
The variety of complete flags  is just the fiber product
$$\mathbb{F}=\mathbb{P}(\mathcal{S}) \times _\mathbb{G}
\mathbb{P}(\mathcal{Q}^{\vee}) = \{(p,\ell,\pi) | \ p \in
\ell \subset \pi \}
\subset\p3\vez\bb G\vez\pd3.
$$

The tautological bundles of these spaces lift
to bundles over $\mathbb{F}$ still denoted by the
same letters. They
fit together into the diagram (pullbacks omitted),
$$\ba{*9c}
\cl O_{\cl Q\ve}(-1)&=&
\cl O_{\pd3}(-1) & \hooklar & \cl Q\ve & 
\hooklar & \cl P\ve & 
\hooklar & S_1:=(\mathbb{C}^4)^{\vee}\\\na6&\ffi:&
\id{x_0}&\subset&\id{x_0,x_1}&\subset&\id{x_0,x_1,x_2}
&\subset&\id{x_0,\dots,x_3}
\ea$$
where the bottom row indicates the corresponding
fibers over our favorite flag (\ref{pdv}) $\varphi_0$.

Now, recall from (\ref{AB}) the rank 7 subbundle 
$\cl A\subset S_3$ of special cubics
and the rank 4 vector subbundle 
$\cl B\subset S_2$ of special quadric cones 
over the variety ${\bb F}$ of complete
flags.
In view of Lemma\,\ref{monsfg},
we have  the surjections$$\ba c
\left(\cl O_{\cl Q\ve}(-2)\otimes S_1\right)\oplus\left(
\cl O_{\cl Q\ve}(-1)\otimes{\cl Q\ve}\otimes\cl
P\ve\right)
\oplus{}{\sym}_3\cl Q\ve\twoheadlongrightarrow
\cl A\subset S_3
\\\na6
(\cl O_{\cl Q\ve}(-1)\otimes\cl P\ve)
\oplus{}{\sym}_2\cl Q\ve
\twoheadlongrightarrow
\cl B\subset S_2.
\ea
$$

By construction, the vector bundles \,\cl A, \cl B \,fit into the
exact sequences 
\be\label{ovA}\ba {*4cl}
S_1\otimes\cl{O}_{\pd3}(-2)& \hooklar& \cl A&
\twoheadlar&
\ov{\cl A}:=\cl A\big/\left(S_1\otimes\cl{O}_{\pd3}(-2)
\right)\\\na6
\cl O_{\pd3}(-2)&\hooklar&
\cl B &\twoheadlar&\ov{\cl B}.
\ea\ee
%
%$$
%\ba c
%\xymatrix%@R1pc@C1pc
%{0 \ \ar@{->}[r]& 
% S_1\otimes\cl{O}_{\pd3}(-2)
% \ \ar@{->}[r]& 
% \cl A
%\ \ar@{->}[r]&
%\ov{\cl A}:=\cl A\big/\left(S_1\otimes\cl{O}_{\pd3}(-2)
%\right)
%\ \ar@{->}[r]& 0}
%\\
%\xymatrix%@R1pc@C1pc
%{
%0 \ \ar@{->}[r]&
%\cl O_{\pd3}(-2) \ar@{->}[r]&
%\cl B  \ar@{->}[r]&\ov{\cl B}\ar@{->}[r]&0.
%}
%\ea$$

\begin{remk}\normalfont
	For later use, remember that knowing the fibers of an
	equivariant, say $\bb C^\star$-vector bundle over 
	the variety \bb F\ of complete flags suffices to get our
	hands on the $\bb C^\star$-equivariant Chern classes.
\end{remk}

We recall the fiber $\ov{\cl A}_{\varphi_0}$
is spanned by the classes  
$$\ov{x_{0}x_{1}^{2}},\ov{x_{0}x_{1}x_{2}},
\ov{x_{1}^{3}} \ \ 
\text{mod}.\,x_0^2S_1.
$$

We have likewise the subbundle $\cl O_{\pd3}(-2)\subset
\cl B$ and the corresponding quotient \ov{\cl B}. 

The fiber of \ov{\cl B} over \ffi \ has the generators
\\\centerline{$\ov{x_0x_1},\ov{x_0x_2},\ov{x_1^2}$
	mod.\,\id{x_0^2}.} 
This shows the following:
\begin{lema}\label{ovAovB}
	\ov{\cl A} is isomorphic to 
	$\left(\cl Q\ve/\cl
	O_{\pd3}(-1)\right)\otimes\ov{\cl B}$. 
\end{lema}

\begin{proof}
	Indeed,
	on the fiber over \ffi\ we have
	$$(
	\cl Q\ve/\cl O_{\pd3}(-1))_{\ffi}=\id{x_0,x_1}/\id{x_0}
	= \id{\ov{x_1}} 
	$$
	hence
	$$
	\left((
	\cl Q\ve/\cl O_{\pd3}(-1))\otimes\ov{\cl B}\right)_{\ffi}=
	\id{\ov{x_1}} \otimes
	\id{\ov{x_0x_1},\ov{x_0x_2},\ov{x_1^2} }
	=\ov{\cl A}_{\ffi},
	$$
	where the = signs mean isomorphisms of
	representations of the stabilizer of the flag
	\ffi.
	
\end{proof}
We define
\be\label{X}
\mathbb{X}=\bb P(\cl B)
\ee 
the corresponding projective bundle of special quadrics.
The fiber $\mathbb{X}_{\varphi_0}%=\{ g \}
$ is the $\mathbb{P}^3$ of special quadrics over the flag
\ffi $=(p_{0},\ell_0,v_0)$ fixed in (\ref{pdv}). 

\medskip
Next, we rewrite the
divisibility condition (\ref{divcond1}) as the following 
system of linear \linebreak equations 
in the variables $(\underline{a})=(a_4,a_5,a_6)$ with
coefficients $(\underline{b})=(b_0,b_1,b_2,b_3)\in\p3$,
\begin{equation}\label{syslin}
\left[\begin{array}{ccc}
2b_3 & 0     & -3b_1\\
0    & -2b_3 & 3b_2\\
-b_2 & b_1   & 0
\end{array}\right] \cdot 
\left[\begin{array}{c}
a_4\\
a_5\\
a_6
\end{array}\right] =
\left[\begin{array}{c}
0\\
0\\
0
\end{array}\right].
\end{equation}

The matrix of the coefficients has determinant zero and 
generic rank two. 
Computing the 2\vez2-minors, we find that the rank of
the coefficients matrix will be less than two when
$b_1=b_2=b_3=0\leftrightarrow x_0^2$. Off the
point $x_0^2$  the solution space is spanned by
the vector product of two rows, 
\\\centerline{$
	%\ba c
	%(-6b_1b_3:-6b_2b_3:-4b_3^2)=( 3b_1: 3b_2: 2b_3)
	%\\
	(-3b_1b_2:-3b_2^2:-2b_2b_3)=( 3b_1:3b_2: 2b_3)
	%\ea
	$}

Our idea now is to describe the locus of special cubic
forms with the divisibility condition as a
$\mathbb{P}^4$-bundle over $\mathbb{X}$ (since it has
codimension two in the $\mathbb{P}^6$ of special
cubics). But to make it possible, we need to replace
$\mathbb{X}$ by a new space, $\mathbb{X}'$,
for which the rank of the
coefficients matrix is two everywhere.

Precisely, we think of the fiberwise solution space to
(\ref{syslin}) as defining a rational map 
$\psi:\bb X\rar \p{}(\ov{\cl A})$, notation as in 
(\ref{ovA}), which on the fiber over the standard
flag $\varphi_0$ (\ref{pdv}) reads
\be\label{b2a}
\psi:\ 
(\un b)\mapsto
(a_4:a_5:a_6)=
(3b_1:3b_2:2b_3).
\ee
Look at the closure $\mathbb{X}'$ of the graph of
$\psi$. On the fiber over $\varphi_0$,
this is the blowup of $\p3=\p{}(\cl B_{\varphi_0})$ 
at the point $x_0^2$. 
As a matter of fact, this turns out to be the
restriction to the fiber over $ \varphi_0$  
of the blowup of $\mathbb{X}$ along the section
$\mathbb{P}\left(\mathcal{O}_{\pd3}(-2)\right)
\hookrightarrow \mathbb{P}(\mathcal{B})$ over \bb
F. We have the diagram
$$%\be\label{PA=PB}
\xymatrix%@R1pc@C1pc
@R2.5pc
{&&
\,  \X'
\ar@{->}[d]_{\text{\small (blowup)}}
\ar@{->}[drr]^{\ \ \ \ \psi':\text{\small\,\,(\p1-bundle)}}
 &&&
\\ 
\mathbb{P}\left(\mathcal{O}_{\pd3}(-2)\right)
\ \ \ar@{=}[drrr] 
\!\! \ar@{^(-}[r]& 
%\subset
\p{}(\cl B) \ar@{=}[r]& \X
\ar@{->}[dr]
\ar@{-->}[rr]^{\!\hskip-.31cm^\psi}&&\p{}(\ov{\cl A})
\ \ar@{->}[dl] \ar@{=}[r]& \p{}(\ov{\cl B})
\\&&&\ \ \ \bb F \ \ \ & &
}
$$%    \ee
where the leftmost inclusion is defined by
taking the square of the equation of the plane in
the flag.
The equality $\p{}(\ov{\cl A})=\p{}(\ov{\cl
	B})$ comes from Lemma\,\ref{ovAovB}; under this
identification, we have 
$$\cl O_{\ov{\cl A}}(-1)=
	\cl O_{\ov{\cl B}}(-1)\otimes (
	\cl Q\ve/\cl O_{\pd3}(-1)).
$$
Pulling back the tautological line
subbundle via $\psi'$ , we get the diagrams
\be\label{ovAB}
\ba {*5c}
S_1\otimes\cl{O}_{\pd3}(-2)
&\hooklar&\cl A'&
\twoheadlar&
\cl O_{\ov{\cl A}}(-1)\\
||&&\hookdownarrow&&\hookdownarrow\\
S_1\otimes\cl{O}_{\pd3}(-2) &\hooklar& {\cl A}&
\twoheadlar&
\ov{\cl A}\
\ea
\ee
and
\be\label{ovAB1}
\ba {*5c}
\cl{O}_{\pd3}(-2)
&\hooklar&\cl B'&
\twoheadlar&
\cl O_{\ov{\cl B}}(-1)\\
||&&\hookdownarrow&&\hookdownarrow\\
\cl{O}_{\pd3}(-2) &\hooklar& {\cl B}&
\twoheadlar&
\ov{\cl B}\
\ea
\ee
where $\rk\cl A'=5,\rk\cl B'=2.$ We have 
\be\label{defX'}
\bb X'=\p{}(\cl B'),
\ee a \p1
bundle over the \p2 bundle \p{}(\ov{\cl B}) over
\bb F. 
We define
\be\label{Y}
\bb Y=\p{}(\cl A'),
\ee
an equivariant \p4--bundle over $\bb X'$.
\begin{prop}
A
general point in \bb Y corresponds to a pair
$(g,f)$ of special quadric, cubic satisfying the
divisibility condition.
\end{prop}
\begin{proof}
  The assertion follows from the previous
  considerations. 
\end{proof}

\medskip
\noindent{{\bfseries Remark.}}
Since all constructions performed so far, as well
as in the sequel are equivariant, 
henceforth we drop the reference to the fiber over
our favorite flag
\ffi, and simplify notation writing, until further
notice,
$$\bb X=\bb
X_{\ffi},\ \cl A=\cl A_{\ffi},\dots,\  etc.
$$

In view of Lemma \ref{ovAovB}, the rational map (\ref{b2a}) can also be written as the rational linear projection map
$$\begin{array}{cccc}
  \psi_{\overline{\mathcal{B}}}: & \mathbb{X} &
\dashrightarrow   & \mathbb{P}(\overline{\mathcal{B}})\\\na8
& g:=b_0x_0^2 + b_1x_0x_1 + b_2x_0x_2 + b_3x_1^2 &
 \longmapsto
& g':=b_1\ov{x_0x_1} + b_2\ov{x_0x_2} + b_3\ov{x_1^2}
\end{array}$$
%$$
%\xymatrix%@R1pc@C1pc
%{
%  \psi_{\overline{\mathcal{B}}}: & \mathbb{X}
%  \  \  \  \  \  \ \ar@{-->}[r] &
%  \ \ 
%  \mathbb{P}(\overline{\mathcal{B}})\\
%  & g:=b_0x_0^2 + b_1x_0x_1 + b_2x_0x_2 + b_3x_1^2 
% \ \ar@{|->}[r] & g':=b_1x_0x_1 + b_2x_0x_2 + b_3x_1^2
%}
%$$
(class mod.\,$x_0^2$). Write for short 
\begin{equation} \label{X'}
  \mathbb{X}' = \{ (g,g') =
  \left((b_0:b_1:b_2:b_3),(u_1:u_2:u_3)\right)
  | \ b_iu_j=b_ju_i  \},
\end{equation}
where
\begin{equation} \label{defg'}
g':=u_1\ov{x_0x_1} + u_2\ov{x_0x_2} + u_3\ov{x_1^2}.
\end{equation}

Here, the $\underline{u}$ are 
homogeneous coordinates in \p2 = fiber of $\p{}(
\ov{\cl B})$ over \ffi.
The map $\psi'$ is (see  (\ref{b2a}))
$$\begin{array}{cccc}
\psi': & \mathbb{X'} & \longrightarrow    & \mathbb{P}(\overline{\mathcal{A}})\\
& (g,g')%=u_1x_0x_1 + u_2x_0x_2 + u_3x_1^2 
& \mapsto & 3u_1\ov{x_0x_1^2} + 3u_2\ov{x_0x_1x_2}
+ 2u_3\ov{x_1^3} \\
\end{array}$$

Notice that the divisibility condition
(\ref{divcond1})
has now the expression

 \begin{equation}\label{divcond}
   \left\{\begin{array}{lcl}
3a_6u_1 &=& 2a_4u_3\\
3a_6u_2 &=& 2a_5u_3\\
a_4u_2 &=& a_5u_1
\end{array}\right.\end{equation}

Set
$$\bb U_0:=
	\bb A^3\subset \bb X=\p3\ni(1:0:0:0)
	\leftrightarrow x_0^2,$$
{\em i.e,} we put $b_0=1$ and $b_1,b_2,b_3$ are affine
coordinates 
on $\bb U_0 $. We write
$$
	g=x_0^2+b_1x_0x_1+b_2x_0x_2+b_3x_1^2\in \bb U_0.
	$$
Let $\bb U_0'\subset\bb U_0\vez\p2$ 
be defined by $u_1=1.$
So $b_2=b_1u_2,b_3=b_1u_3$. Thus the
affine coordinates on $\bb U_0'$ are
$b_1,u_2,u_3$. We have
$$g'=x_0x_1+u_2x_0x_2+u_3x_1^2.$$

The rational map (\ref{b2a}) gives 
$$(a_4:a_5:a_6)=(3:3u_2,2u_3) , \ \text{ so }
\ \frac{a_4}3=\frac{a_5}{3u_2}=\frac{a_6}{2u_3}
\  \text{ and } \ a_5=u_2a_4,a_6=\frac23u_3a_4.
$$
Recall $\ov{\cl A}=\id{\,\ov{x_{0}x_{1}^{2}},\ov{x_{0}x_{1}x_{2}},
	\ov{x_{1}^{3}} \ \ 
	\text{mod}.\,x_0^2S_1\,}%=\ov{x_1}\ov{\cl B}.
$. We have the trivialization,
$$
\cl O_{\ov{\cl A}}(-1)_{|\bb U_0'}=
\id{\,\ov{x_{0}x_{1}^{2}}+u_2\ov{x_{0}x_{1}x_{2}}+\frac23
	u_3  \ov{x_{1}^{3}}\,}
$$
and
$$\cl A'_{|\bb U_0'}=x_0^2S_1+
\id{\,{x_{0}x_{1}^{2}}+u_2{x_{0}x_{1}x_{2}}+\frac23
	u_3  {x_{1}^{3}}\,}.
$$
We have $\bb Y_{|\bb U_0'}=  \bb U_0' \vez\p4$. We
redefine homogeneous coordinates $(a_i)$ on 
the \p4 factor and write
\be\label{fP4}\left\{\ba l
f=(a_0x_0 + a_1x_1 + a_2x_2+a_3x_3)x_0^2 +
a_4({x_{0}x_{1}^{2}}+u_2{x_{0}x_{1}x_{2}}
+\frac23u_3{x_{1}^{3}}),\\\na5
g=x_0^2+b_1(x_0x_1+u_2x_0x_2+u_3x_1^2).
\ea\right.
\ee
For example, choose  $g=x_0^2 ,\ \  g'=x_0x_1+
u_2x_0x_2+u_3x_1^2 \ \in \mathbb{X}', \ (u_2,u_3) \in \mathbb{A}^2$ .  We may pick $f=
{x_{0}x_{1}^{2}}+u_2{x_{0}x_{1}x_{2}} 
+\frac23u_3{x_{1}^{3}} $ in the fiber of \bb Y
over $(g,g')$ and evaluate 
\begin{equation} \label{wcomus}
\omega=\dfrac{3fdg - 2gdf}{x_0} = 
\end{equation}
$$
=(4x_0x_1^2+4u_2x_0x_1x_2+4u_3x_1^3)dx_0 -
(4x_0^2x_1+2u_2x_0^2x_2+4u_3x_0x_1^2)dx_1 -
2u_2x_0^2x_1dx_2.
$$
Thus, $\omega$ spans a linear space of dimension
one over the chosen points.
If this space were well defined
for all points of $\mathbb{Y}$, we would have a line
subbundle $\mathcal{W}$ of the trivial bundle of twisted
differential forms, pullback of the $\cl O_{\p{44}}(-1)$ of the
corresponding projective space where $E(3)$ lives.
%and the projective bundle $\mathbb{P}(W)$ should give an
%exceptional foliation for any element in $\mathbb{Y}$,
%and 
So $\mathbb{Y}$ would be an adequate parameter space for
the exceptional component. 

However, there are indeterminacies in \bb
Y\ corresponding to the vanishing of the form $w=3fdg-2gdf$,
{\em e.g.,} when $f=x_0^3$, $g=x_0^2$,
just to give an example.

\medskip
The next step is to solve the indeterminacies of the rational map
\be\label{mapomega}
\begin{array}{ccc}
\mathbb{Y} & \makebox{\Large$\dashrightarrow$}    & \mathbb{P}(H^0(\Omega^1_{\mathbb{P}^3}(4)))\\\na8
(g,g',f) & \longmapsto & \omega = \dfrac{3fdg - 2gdf}{x_0}\\
\end{array}
\ee

Note  that $g'$ is implicit in this formula: its
coefficients appear in $f$, cf.\,(\ref{fP4}).

We will see below (cf.\,\S\ref{indtw}) that the indeterminacy locus of this map is a
non-reduced and reducible scheme. But its reduction
 consists of two beautiful components. 

One
of them corresponds to the points $(g,g',f)$ parametrized
by $(s:t) \in \mathbb{P}^1$, where 
$$f=(sx_0+tx_1)^3, \ \ g=(sx_0+tx_1)^2, \ \ g'=x_1(2sx_0
+ tx_1).
$$

The other  component corresponds to a section over the whole
$\mathbb{P}^2$ fiber of the \linebreak exceptional divisor of
$\mathbb{X}'$ over  $g=x_0^2$, with $f=x_0^3$, {\em i.e.,}
the points of the form $(x_0^2,g',x_0^3)$, 
with $g'$  as in (\ref{defg'}), cf. the diagram below:
%a component produced by our replacement of $\mathbb{X}$
%by $\mathbb{X}'$.    
%\begin{equation}\label{diagE}
%\ba {*5c}
%\bb Y&
%\stackrel{\p4-\text{bundle}}{
%	\lar}&\bb X'&
%\stackrel{\text{blowup}}
%\lar&\bb X=\p3
%\\\cup&&\cup&&
%\\
%\{(x_0^2,g',x_0^3)\}
%&\leftrightarrow&\{(x_0^2,g')\}=\p2&\lar&x_0^2
%\ea
%\end{equation}
\begin{equation}\label{diagE}
\xymatrix%
@R.781 pc@C1.951 pc
{
  \bb Y
  \ \ar@{->}[r]^{\p4-\text{bundle}}&
  \bb X'
  \ \ar@{->}[r]^{\text{\!\!\!blowup}} & \bb X=\p3
  \\\bigcup&\bigcup\phantom'&\bigcup\\
  \{(x_0^2,g',x_0^3)\}
  \ \ar@{<->}[r]&\{(x_0^2,g')\}=\p2
 \ \ar@{->}[r]& \left\{x_0^2\right\}
}
\end{equation}

This second component is then the $\mathbb{P}^2$
%(see (\ref{defg'}))
projectivization of the space
\begin{equation}\label{geraE}
\langle \ov{x_0x_1}, \ov{x_0x_2}, \ov{x_1^2} \rangle .
\end{equation}

%Notice that the second component (the $\mathbb{P}^2$) is
%a nice complete intersection, but the first described component
% is not. However, it is a locally complete
% intersection, as we will see in the next section.

\medskip

The (global) equations of the indeterminacy locus
of the map\,(\ref{mapomega})
are too hard to
treat. However, %when we see it
locally they look rather
friendly.
Let us %, so we will 
work over an affine cover of $\mathbb{X}'$ to deal with it.

\section {Indeterminacies of $\omega$}
\label{indtw}
Let's sketch our plan to solve the indeterminacies of the map $(g,g',f) \mapsto \omega$\,(\ref{mapomega}).

The links to appendix\,\ref{apA} refer to scripts for the software \textit{Macaulay2} that perform the calculations.

%The (global) equations of the indeterminacy locus are too hard to
%treat. However, when we see it locally they look very
%friendly. Now, we start the study over all standard
%neighborhoods that cover $\mathbb{X}'$. \\ 

Notation as in (\ref{X'}), we list 
 the typical neighborhoods to consider.

\subsection{$\bullet \ \ b_3=1$}\label{a6b3=1txt}

In this case, $u_1=b_1, u_2=b_2, u_3=1$ and the
divisibility condition (\ref{divcond1}),
(\ref{divcond}) is
$$\left\{\begin{array}{l}
3a_6b_1 = 2a_4\\
3a_6b_2 = 2a_5\\
a_4b_2 = a_5b_1
\end{array}\right.$$
Observe that the third equation is irrelevant, and we can
write 
\begin{equation}\label{fga6b3=1}
\left\{
\ba l
f=\left(a_0x_0 + a_1x_1 + a_2x_2+a_3x_3\right)x_0^2 +
a_6\left(\frac{3}{2}b_1x_0x_1^2+\frac{3}{2}
b_2x_0x_1x_2 +   x_1^3\right),
\\\na6
g=b_0x_0^2 + b_1x_0x_1 + b_2x_0x_2 + x_1^2.
\ea\right.
\end{equation}
Presently our \p4 of special cubics 
has homogeneous coordinates
$a_0,a_1,a_2,a_3,a_6$. 
Compute $\omega$ as in \ref{compw1}. We find
$$\omega=\left(\textcolor{blue}{(3a_0b_1-2a_1b_0)}x_0^2x_1
  + \cdots
  + \textcolor{blue}{a_3b_2}x_0x_2x_3\right)
  dx_0 + \cdots +
 \blu{2a_3}
\left(\textcolor{blue}{b_0}x_0^3 + \cdots
+ \textcolor{blue}{b_2}x_0^2x_2\right)dx_3.
$$
%
%$$\omega=[\textcolor{blue}{(-2a_1b_0+3a_0b_1)}x_0^2x_1 + \textcolor{blue}{(-6a_6b_0b_1+a_1b_1+6a_0)}x_0x_1^2 + $$ 
%$$ + \textcolor{blue}{(-(3/2)a_6b_1^2-6a_6b_0+4a_1)}x_1^3 + \textcolor{blue}{(-2a_2b_0+3a_0b_2)}x_0^2x_2 + $$ $$ + \textcolor{blue}{(-6a_6b_0b_2+a_2b_1+a_1b_2)}x_0x_1x_2 + 
%\textcolor{blue}{(-3a_6b_1b_2+4a_2)}x_1^2x_2 + \textcolor{blue}{a_2b_2}x_0x_2^2 + $$  $$ \textcolor{blue}{- (3/2)a_6b_2^2}x_1x_2^2 - \textcolor{blue}{2a_3b_0}x_0^2x_3 + \textcolor{blue}{a_3b_1}x_0x_1x_3 + 
%\textcolor{blue}{4a_3}x_1^2x_3 + \textcolor{blue}{a_3b_2}x_0x_2x_3]dx_0 + $$ $$
%[\textcolor{blue}{(2a_1b_0-3a_0b_1)}x_0^3 + \textcolor{blue}{(6a_6b_0b_1-a_1b_1-6a_0)}x_0^2x_1 + $$ $$ + 
%\textcolor{blue}{((3/2)a_6b_1^2+6a_6b_0-4a_1)}x_0x_1^2 + 
%\textcolor{blue}{(3a_6b_0b_2-3a_2b_1+2a_1b_2)}x_0^2x_2 + $$ $$ + 
%\textcolor{blue}{((9/2)a_6b_1b_2-6a_2)}x_0x_1x_2 + \textcolor{blue}{3a_6b_2}^2x_0x_2^2 \textcolor{blue}{- 
%	3a_3b_1}x_0^2x_3\textcolor{blue}{-6a_3}x_0x_1x_3]dx_1 + $$  $$ + 
%[\textcolor{blue}{(2a_2b_0-3a_0b_2)}x_0^3 + 
%\textcolor{blue}{(3a_6b_0b_2+2a_2b_1-3a_1b_2)}x_0^2x_1 + $$  $$ + \textcolor{blue}{(-(3/2)a_6b_1b_2+2a_2)}x_0x_1^2 \textcolor{blue}{- a_2b_2}x_0^2x_2 \textcolor{blue}{- (3/2)a_6b_2^2}x_0x_1x_2 \textcolor{blue}{- 3a_3b_2}x_0^2x_3]dx_2 + $$  $$ + 
%[\textcolor{blue}{2a_3b_0}x_0^3 + \textcolor{blue}{2a_3b_1}x_0^2x_1 + \textcolor{blue}{2a_3}x_0x_1^2 + \textcolor{blue}{2a_3b_2}x_0^2x_2]dx_3. $$

Notice that if we set $a_6=0$ the expression for 
$\omega$ acquires the terms
$$[\cdots + \textcolor{blue}{(a_1b_1+6a_0)}x_0x_1^2 + \textcolor{blue}{4a_1}x_1^3+\textcolor{blue}{4a_2}x_1^2x_2+\textcolor{blue}{4a_3}x_1^2x_3 + \cdots ]dx_0 + \cdots$$

Hence, if this is zero then 
\begin{equation}\label{as=0}
a_0=a_1=a_2=a_3=a_6=0,
\end{equation}
which is impossible since $(a_0:a_1:a_2:a_3:a_6) \in \mathbb{P}^4$.

It means that over the neighborhood $b_3=1$ we do not
have any indeterminacies when $a_6=0$. Thus in order to study the indeterminacy locus we can take $a_6=1$.

Setting $a_6=1$ in the expression of $\omega$, and collecting the coefficients in the resulting 1-form, we see that the indeterminacy locus is given by the ideal (cf.\,\ref{obtemJ1})
$$J=\langle
a_3,3b_1b_2-4a_2,6a_0-6b_0b_1+a_1b_1,3b_1^2-8a_1+12b_0,
b_2^2, b_2(a_1-3b_0),(a_1-3b_0)^2 \rangle.
$$
Note that it contains the square $\langle b_2,a_1-3b_0 \rangle ^2$.\\

The radical  is 
\begin{equation} \label{eqC}
%K
J_{red}
:=rad(J)=\langle 8a_0-b_1^3,a_2,b_2,a_3,4b_0-b_1^2,4a_1-3b_1^2 \rangle.
\end{equation}
It  is generated by a regular sequence (the affine
coordinates here are $a_0,a_1,a_2,a_3,b_0,b_1,b_2$),
therefore we have a locally complete intersection. 
%Moreover, the equations in $J_{red}$ show that 
Plugging (\ref{eqC}) into (\ref{fga6b3=1}) we find
\begin{equation} \label{expC}
f=\left(\frac{b_1}{2}x_0+x_1\right)^3, \ \ g=\left(\frac{b_1}{2}x_0 + x_1\right)^2.
\end{equation}

\begin{defi}\normalfont
We call $C$  the curve defined locally by the
above ideal $J_{red}$.
\end{defi} Blowing up $\mathbb{Y}$ along this one
dimensional variety $C$ does not solve the
indeterminacies yet, but makes the resulting ones very nice.

Indeed, let $\mathbb{A}^7$ the open dense set defined by $b_3=a_6=1$. The blowup of this $\mathbb{A}^7$ along $C$ \,is
\begin{equation}\label{Y1a6b3}
\mathbb{Y}_1 = \{\left((\underline{a},\underline{b}),(s_0:\ldots:s_5)\right)  \ | \ s_i \cdot e_j = s_j \cdot e_i\} \subset \mathbb{A}^7 \times \mathbb{P}^5,
\end{equation}
where $e_i, \ 0\leq i \leq 5$, are the equations
of $J_{red}$, ordered as displayed in (\ref{eqC}). 

We will write the equations of $\mathbb{Y}_1$ as
\begin{equation}\label{eqbup1}
\dfrac{8a_0-b_1^3}{s_0}=\dfrac{a_2}{s_1}=\dfrac{b_2}{s_2}=\dfrac{a_3}{s_3}=\dfrac{4b_0-b_1^2}{s_4}=\dfrac{4a_1-3b_1^2}{s_5}\  .
\end{equation}

The exceptional divisor \
$\mathbb{P}(\mathcal{N}_{C|\mathbb{Y}})$ \ is a  $\mathbb{P}^5$-bundle over $C$, where $\mathcal{N}_{C|\mathbb{Y}}$ denotes the normal bundle of $C$ in $\mathbb{Y}$.

Choose $a_2$ for the local equation of the exceptional
divisor, so take $s_1=1$. The equations (\ref{eqbup1}) become (cf.\,\ref{blowC1}),
$$\left\{\begin{array}{rcl}
8a_0-b_1^3 & = & s_0a_2 \\
b_2        & = & s_2a_2 \\
a_3        & = & s_3a_2 \\
4b_0-b_1^2 & = & s_4a_2 \\
4a_1-b_1^2 & = & s_5a_2.
\end{array}
\right.$$
Substituting into the expression of $\omega$, we find
(see the expression for
wnt$\underline{\ \ }$1 in \ref{blowC1}), 
yes,  a horrible expression. But is easy to see that
the equation $a_2$ of the exceptional divisor
divides this new $\omega$.
Then, we can perform this division to
obtain the new expression of $\omega$ in the present coordinates $b_1,a_2,s_0,s_2,s_3,s_4,s_5$. And now we are able to look at the indeterminacy locus of this new $\omega$, finding the generators (ideal Jn$\underline{\ \ }$1 in \ref{blowC1})
\begin{equation}\label{j''a2}
J''=\langle 3s_4-2s_5,s_3,s_0-b_1s_5,4-3b_1s_2,a_2 \rangle
\end{equation}
(remember that we are doing this on the neighborhood $s_1=1$ of $\mathbb{P}^5$).

We can perform this analysis on the other 5 standard neighborhoods of $\mathbb{P}^5$ to see
that the ideal of the new indeterminacy locus is, in homogeneous
coordinates of $\mathbb{P}^5$, given by the
dehomogenization of (\ref{JR}) w.r.t. the variable $s_i$,
cf. the ideal wn$\underline{\ \ }$i in \ref{othviz}
\begin{equation} \label{JR}
J''=\langle 3s_4-2s_5,s_3,s_0-b_1s_5,4s_1-3b_1s_2,\textrm{exc} \rangle,
\end{equation}
where the exc in ideal $J''$ represents the chosen local
equation of the exceptional divisor (in
(\ref{j''a2})
this is the equation $a_2$).

This ideal $J''$ corresponds to a ruled surface with basis $C$
(due to the presence of the exc equation), which is a
$\mathbb{P}^1$--subbundle of the exceptional divisor: note that the other four equations are homogeneous linear equations in the variables $\underline{s}$.

Hence, we have now a reduced and irreducible 
local complete intersection of  dimension two.
Since it is the full indeterminacy locus,
a blowup along this
subscheme will solve the indeterminacies over the
neighborhood $b_3=1$. Thus, the map to $\omega$
becomes a morphism after blowing up this neighborhood. (cf. \cite{hartshorne2013algebraic}, ex. 7.17.3, page 168).

\begin{defi}\label{defR}\normalfont
We call $R$ the ruled surface defined by the ideal $J''\  (\ref{JR})$.
\end{defi}

Let us register for future use (see Table \ref{x13x12s2}, page \pageref{x13x12s2}) the blowup of $\bb{Y}_1$ along $R$ as
\begin{equation}\label{RRa6b3}
\{\left((\underline{b},\underline{s}),(v_0:\ldots:v_4) \right) \ | \ v_i \cdot e_j = v_j \cdot e_i\} \subset \mathbb{A}^7 \times \mathbb{P}^4,
\end{equation}
where $e_i, \ 0\leq i \leq 4$, are the equations of $J''$, ordered as in (\ref{JR}), 

As in (\ref{eqbup1}), the equations of this blowup are
\begin{equation}\label{eqbup7}
\dfrac{3s_4-2s_5}{v_0}=\dfrac{s_3}{v_1}=\dfrac{s_0-b_1s_5}{v_2}=\dfrac{4s_1-3b_1s_2}{v_3}=\dfrac{\textrm{exc}}{v_4}\  .
\end{equation}

\subsection{$\bullet \ \ b_0=1$ and $u_1=1$.}\label{solvb0u1}
Keep the notation as in (\ref{X'}).

%We will study the indeterminacies over the exceptional
%divisor of $\mathbb{X}'$. Remember that 
%$$\mathbb{X}' = \{ (g,g') =
%(b_0:b_1:b_2:b_3),(u_1:u_2:u_3) | \ b_iu_j=b_ju_i \ , \
%i,j \in \{1,2,3 \} \}.$$ 

Taking $u_1=1$, the interesting equations are
$$ \left\{\begin{array}{l}
b_2 = b_1u_2\\
b_3 = b_1u_3
\end{array}\right\} \ \textrm{ and } \ \left\{\begin{array}{l}
3a_6 = 2a_4u_3\\
a_5 = a_4u_2
\end{array}\right\},$$
where the former corresponds to the blowup of our \p3 of
special quadrics at $x_0^2$ and the latter 
to divisibility.
We have
$$\left\{\ba l
f=\left(a_0x_0 + a_1x_1 + a_2x_2
+a_3x_3\right)x_0^2 +a_4\left(
x_0x_1^2+u_2x_0x_1x_2 + \frac{2}{3}u_3x_1^3\right)
\\\na8
g=x_0^2 + b_1x_0x_1 + b_1u_2x_0x_2 + b_1u_3x_1^2.
\ea\right.
$$

After computing $\omega$ and collecting coefficients, we
realize  that if $a_0=0$ then there is no indeterminacy
(similar to (\ref{as=0})). Hence, take $a_0=1$ (cf.\,\ref{compw2}).

The indeterminacy locus is again a non-reduced scheme. The
reduction now presents a surprise in relation to the previous neighborhood: its ideal can be written  as  \ \ref{compCE}
\begin{equation}\label{idealK}
J_{red}=rad(J)=\langle 3b_1^2-4a_4,a_3,a_2,2a_1-3b_1,b_1u_2,b_1(b_1-4u_3) \rangle,
\end{equation}
so it is reducible. There are two components: (see primaryDecomposition K in \ref{compCE})

\begin{enumerate}
	\item A component $C$, 
with ideal
	\begin{equation} \label{JC}
	J_C=\langle b_1-4u_3,a_1-6u_3,a_2,a_3,a_4-12u_3^2,u_2 \rangle
	\end{equation}
	which represents (compare with (\ref{expC}))
	$$f=(x_0+2u_3x_1)^3, \ \ g=(x_0+2u_3x_1)^2,$$
	the same curve $C$ as before, viewed in
        the present neighborhood.

	\item A component with ideal 
	\begin{equation} \label{JEE}
	J_E=\langle a_1,a_2,a_3,a_4,b_1 \rangle,    
	\end{equation}
	which is the whole $\mathbb{P}^2$--fiber of the
        exceptional divisor of $\mathbb{X}'$ over
        $g=x_0^2$, and $f=x_0^3$. See (\ref{diagE}).  
\end{enumerate}
Notice that there is precisely one  point in the
intersection of these two components, since
$$J_C + J_E = \langle a_1,a_2,a_3,a_4,b_1,u_2,u_3 \rangle .$$
This ideal represents the single point $(g,g',f)=(x_0^2,x_0x_1,x_0^3) \in \mathbb{Y} $.

\begin{defi}\label{defE} \normalfont
	We call $E$ the component given by the ideal $J_E$ (\ref{JEE}).
\end{defi}

Let's see what happens when we make the blowup of
$\mathbb{Y}$ first along $C$, followed by a blowup on
$E'$, the strict transform of $E$.

\subsubsection{Blowing up along $C$}

Now, let $\mathbb{A}^7$ be the affine neighborhood defined by $b_0=a_0=u_1=1$. The blowup of this $\mathbb{A}^7$ along $C$ is
$$\mathbb{Y}_1 = \{\left(
(\underline{a},\underline{b}),(s_0:\ldots:s_5)\right)
  \ | \ s_i \cdot e_j = s_j \cdot e_i\} \subset \mathbb{A}^7 \times \mathbb{P}^5,$$
where $e_i, \ 0\leq i \leq 5$, are the equations of $J_C$, ordered as in (\ref{JC}). 

We will write the equations of $\mathbb{Y}_1$ as
\begin{equation}\label{eqbup3}
\dfrac{b_1-4u_3}{s_0}=\dfrac{a_1-6u_3}{s_1}=\dfrac{a_2}{s_2}=\dfrac{a_3}{s_3}=\dfrac{a_4-12u_3^2}{s_4}=\dfrac{u_2}{s_5}\  .
\end{equation}

There are many choices for the local equation of the first
exceptional divisor, as we see in (\ref{JC}). Let us pick, say,
exc$_1=u_2$. 
Then, we make the following substitutions on the
expression of $\omega$, which corresponds to $s_5=1$ in
(\ref{eqbup3}),
(cf.\,\ref{primC}):
$$\left\{\begin{array}{rcl}
b_1-4u_3    & = & s_0u_2 \\
a_1-6u_3    & = & s_1u_2 \\
a_2         & = & s_2u_2 \\
a_3         & = & s_3u_2 \\
a_4-12u_3^2 & = & s_4u_2
\end{array}\right.
$$

The new 7 affine coordinates on the blowup are 
\begin{equation}\label{affineY1}
s_0,s_1,s_2,s_3,s_4,u_2,u_3.
\end{equation}

%($``="(b_1-4u_3:a_1-6u_3:a_2:a_3:a_4-12u_3^2:u_2)$). 
After
these substitutions, one can check that the new
expression 
for $\omega$ 
is divisible by $u_2$. Dividing, we obtain the transform of $\omega$ over the first blowup (cf. wn$\underline{\ \ }$5 in \ref{primC}).

\begin{defi}\normalfont
Denote by $\mathbb{Y}_1 \rightarrow \mathbb{Y}$ the blowup of $\mathbb{Y}$ along $C$.
\end{defi}

Recall (\ref{idealK}) the radical $J_{red}=rad(J)$ represents
the union $C \cup E$. After the blowup of $\mathbb{Y}$
along $C$, calculations as in \ref{primC} reveal that 
%the new
%indeterminacies coincide with the strict
%transform of $J_{red}$. The second step is to blowup along this component. 
%
the strict transform of $J_{red}$ is given by 
%now locally
\begin{equation} \label{KE'}
J_E'=\langle s_3,s_2,3s_0-2s_1,6u_3+s_1u_2,3s_4+s_1^2u_2 \rangle.
\end{equation}
It coincides with the strict transform $E'$ of $E$ under
the first blowup. 

\subsubsection{Blowup along $E'$}

We take $E'$ as our new blowup center.

\begin{defi}\normalfont
Denote by $\mathbb{Y}_2 \rightarrow \mathbb{Y}_1$ the blowup of $\mathbb{Y}_1$ along $E'$.
\end{defi}

After the blowup along $C$, the reduction of the
new indeterminacy locus in $\mathbb{Y}_1$ consists
again of two components: the component $E'$
(\ref{KE'}), and the ruled surface $R$ (Definition
\ref{defR}), with ideal (cf.\,\ref{secE})
$$\langle u_2,\ s_3,\ s_2-6u_3,\ 3s_0-2s_1,\ 3s_1u_3 - s_4 \rangle.$$
Note the presence of the local equation
exc$_1=u_2$, representing the base $C$, and the
other four linear equations in the ``fiber'' variables $\underline{s}$.

These two components intersect in codimension one in $R$, 
since the sum of their ideals is (cf.\,\ref{secE})
$$\langle u_3,\ u_2,\ s_4,\ s_3,\ s_2,\ 3s_0  - 2s_1
\rangle.$$ Thus, the blowup along $E'$ does not change
the structure of $R$, that is, its strict transform $R'
\subset \mathbb{Y}_2$ is a ruled surface, a
$\mathbb{P}^1$--bundle over the (transform of the) curve $C$.

%We will simplify the notations from now on writing
%$E'=E, \ R'=R$ for the strict transforms, without
%danger of confusion.

The blowup of the $\mathbb{A}^7$ (\ref{affineY1}) along $E'$ is
$$\mathbb{Y}_2 = \{\left(
(\underline{u},\underline{s}),(t_0:\ldots:t_4)\right)
  \ | \ t_i \cdot e_j = t_j \cdot e_i\} \subset \mathbb{A}^7 \times \mathbb{P}^4,$$
where $e_i, \ 0\leq i \leq 4$, are the equations of
$J_{E'}$, ordered as in (\ref{KE'}), and
$(\underline{u},\underline{s})$ as in (\ref{affineY1}).

We will write the equations of $\mathbb{Y}_2\subset
\mathbb{A}^7 \times \mathbb{P}^4
$ as
\begin{equation}\label{eqbup4}
\dfrac{s_3}{t_0}=\dfrac{s_2}{t_1}=\dfrac{3s_0-2s_1}{t_2}=\dfrac{6u_3+s_1u_2}{t_3}=\dfrac{3s_4+s_1^2u_2}{t_4}\  .
\end{equation}

Choose now the equation exc$_2=s_2$ in (\ref{KE'})  as
the new local exceptional equation, {\em i.e,} take $t_1=1$. Equations (\ref{eqbup4}) become
$$\left\{
\begin{array}{rcl}
s_3           & = & t_0s_2 \\
3s_0-2s_1     & = & t_2s_2 \\
6u_3+s_1u_2   & = & t_3s_2 \\
3s_4+s_1^2u_2 & = & t_4s_2
\end{array}\right.
$$
The new 7 affine coordinates on $\mathbb{Y}_2$ are 
\begin{equation}\label{affineY2}
t_0,t_2,t_3,t_4,u_2,s_1,s_2.
\end{equation}

Again, the expression of the renewed $\omega$, obtained
by these substitutions, is divisible by the equation of
the exceptional divisor ($s_2$ in this case). 

However, after this division, the map to $\omega$ is not
solved yet, (cf.\,\ref{secE}).

There are now two components in the new indeterminacy locus. One of them is 
\begin{equation} \label{JR2}
J_{R'}=\langle u_2,t_3-1,t_2,t_0,3s_1-2t_4 \rangle.
\end{equation}
The ideal $J_{R'}$ represents the ruled surface $R'$ which is a
\p1--subbundle of the (transform of the) exceptional divisor $\p{}(\cl
N_{C|\bb Y})$ (the same as (\ref{JR}), viewed in other neighborhood).

\subsubsection{Blowup along $R'$}

\begin{defi}\normalfont
Denote by $\mathbb{Y}_3 \rightarrow \mathbb{Y}_2$ the blowup of $\mathbb{Y}_2$ along $R'$.
\end{defi}

This blowup, in the affine chart $s_5=1,\ t_1=1$, is 
 $$\mathbb{Y}_3 = \{\left(
(\underline{u},\underline{s},\underline{t}),(v_0:\ldots:v_4)  \right)
\ | \ v_i \cdot e_j = v_j \cdot e_i\} \subset \mathbb{A}^7 \times \mathbb{P}^4,$$
 where $e_i, \ 0\leq i \leq 4$, are the equations of
 $J_{R'}$, ordered as in (\ref{JR2}), and
 $(\underline{u},\underline{s},\underline{t})$ as in (\ref{affineY2}). 

We will write the equations of $\mathbb{Y}_3$ as
\begin{equation}\label{eqbup5}
\dfrac{u_2}{v_0}=\dfrac{t_3-1}{v_1}=\dfrac{t_2}{v_2}=\dfrac{t_0}{v_3}=\dfrac{3s_1-2t_4}{v_4}\  .
\end{equation}

Choose now the equation exc$_3=u_2$ in (\ref{JR2})  as
the new local exceptional equation, {\em i.e,} take $v_0=1$. Equations (\ref{eqbup5}) become
$$\left\{
\begin{array}{rcl}
t_3-1     & = & v_1u_2 \\
t_2       & = & v_2u_2 \\
t_0       & = & v_3u_2 \\
3s_1-2t_4 & = & v_4u_2
\end{array}\right.
$$

The new 7 affine coordinates on $\mathbb{Y}_3$ are 
\begin{equation}\label{affineY3}
v_1,v_2,v_3,v_4,u_2,s_2,t_4.
\end{equation}

In $\bb{Y}_3$ the indeterminacy locus reduces to (Jn$\underline{\ \ }$11 in \ref{tercR})
\begin{equation}\label{J'''L}
J_L=\langle s_2,v_1,v_2,v_3,v_4,t_4 \rangle.    
\end{equation}

\subsubsection{Blowup along $L$}

The component represented by the ideal $J_L$ (\ref{J'''L}) is easy to
describe. Since the affine variables are as listed in
(\ref{affineY3}),
the six linear equations represent a
line inside the (transform of) the second blowup center $E$ (note the presence of the equation exc$_2 = s_2$ in $J_L$), a line parametrized by the variable
$u_2$. 
\begin{defi}\normalfont
We call $L$ the component described by the ideal $J_L \  
(\ref{J'''L})$.
\end{defi}

Finally, we solve the indeterminacies by blowing up along the still remaining indeterminacy locus, the component $L$. 

\begin{defi}\normalfont
We will denote by $\mathbb{Y}_4 \rightarrow \mathbb{Y}_3$ the blowup of $\mathbb{Y}_3$ along $L$.
\end{defi}

To put coordinates for later use (see Table \ref{x01s5pl3v0}, page \pageref{x01s5pl3v0}), set
$$\mathbb{Y}_4 = \{\left(
(\underline{u},\underline{s},\underline{t},\underline{v}),(z_0:\ldots:z_5)\right)  \ | \ z_i \cdot e_j = z_j \cdot e_i\} \subset \mathbb{A}^7 \times \mathbb{P}^5,$$
where $e_i, \ 0\leq i \leq 5$, are the equations of $J_L$, ordered as in (\ref{J'''L}), and $(\underline{u},\underline{s},\underline{t},\underline{v})$ as 
in (\ref{affineY3}).
The equations of $\mathbb{Y}_4
\subset \mathbb{A}^7 \times \mathbb{P}^5$ are
\begin{equation}\label{eqbup6}
\dfrac{s_2}{z_0}=\dfrac{v_1}{z_1}=\dfrac{v_2}{z_2}=\dfrac{v_3}{z_3}=\dfrac{v_4}{z_4}=\dfrac{t_4}{z_5}\  .
\end{equation}

Over $\mathbb{Y}_4$ the map is solved, at least in the
affine chart [$s_5=t_1=v_0=1$], (cf.\,\ref{quartL}). 

The calculations in all other standard neighborhoods
reveal that this sequence of blowups still work to solve
the map. 
(See\,\ref{othviz2}).
 
 \begin{teo}
The four blowups of $\mathbb{Y}$ along $C$, $E$, $R$ and $L$ will solve the map $(g,f) \mapsto \omega$ on the neighboorhod $b_0=1, \ u_1=1$.    
 \end{teo}

There are other four standard neighborhoods of
$\mathbb{Y}$ to be checked,  corresponding to
$[b_2=1], \ [b_1=1], \ [b_0=1 \textrm{ and }
  u_2=1], \ [b_0=1 \textrm{ and } u_3=1].$ The
four blowups \linebreak described will solve the
maps over each of these neighborhoods, cf.\,\ref{a0b0u2=1}, \ \ref{a0b0u3=1}, \ \ref{b2=1}

\section{Summary}%ized strategy} 

The whole discussion of the indeterminacy loci was made over our fixed flag (\ref{pdv}).
%$$\ffi:\,p_0=\{x_0=x_1=x_2=0 \} \in \ell_0=\{x_0=x_1=0 \} \subset
%v_0=\{x_0=0 \}.
%$$
Since the blowup centers as described were actually fibers of
 bundles over the variety of complete flags  ${\bb F}$, 
%But over other flags it is just . Hence running through the  
we have at the end a bundle which constitutes the
desired parameter space:
%(notice that $\mathbb{X}$ is a $\mathbb{P}^3$-bundle
%over the complete flags). 
\begin{equation} \label{summarized}
\ba c
\mathbb{Y}_4   \xrightarrow{\textrm{blowup } L} \mathbb{Y}_3 \xrightarrow{\textrm{blowup }R} \mathbb{Y}_2  \xrightarrow{\textrm{blowup }E} \mathbb{Y}_1  \xrightarrow{\textrm{blowup }C} \mathbb{Y} \\\\
\mathbb{Y} \xrightarrow{\mathbb{P}^4 \textrm{ bundle}} \mathbb{X}' \xrightarrow{\textrm{blowup }\mathbb{P}(\mathcal{O}_{\pd3}(-2))} \mathbb{X}
\xrightarrow{\mathbb{P}^3 \textrm{ bundle } } {\bb F}=\textrm{ all flags}
\ea
\end{equation}

This can be summarized as follows:

\begin{teo}\label{Y4}
Let $\mathbb{Y}_4$ be the variety obtained by the four
blowups as described above. Then $\mathbb{Y}_4$ 
is
equipped with a morphism 
onto the exceptional component $E(3)$.
\qed\end{teo}

By construction, $\mathbb{Y}_4$ is equipped with a line bundle $\mathcal{W}$, pullback of the line bundle $\mathcal{O}_{\p{44}}(-1)$ over 
${\p{44}}=\mathbb{P}(H^0(\Omega^1_{\mathbb{P}^3}(4)))$. The fiber of $\mathcal{W}$ over a point in $\mathbb{Y}_4$ is the rank one space spanned by the computed 1-form $\omega$.

$$\begin{array}{ccccc}
\mathcal{W}=\Phi^*(\mathcal{O}(-1))&   & \mathcal{O}(-1)  & &  \\
\downarrow &  & \downarrow  & & \\
\mathbb{Y}_4 & \stackrel{\Phi}\longrightarrow    & \mathbb{P}(H^0(\Omega^1_{\mathbb{P}^3}(4))) & \supset & E(3)=\Phi(\mathbb{Y}_4)\\
\end{array}$$

\begin{prop}
The map $\Phi: \mathbb{Y}_4 \longrightarrow E(3)$ is generically injective.
\begin{proof}
	For a given $\omega \in E(3)$ (off the boundary) we already saw that the flag can be recovered (see Remark \ref{recupflag}, page \pageref{recupflag}). Hence, we can look at a fiber over a fixed flag. 
	
	So fix the flag $\ffi=(p_{0},\ell_{0},v_0) (\ref{pdv})$. The fiber of $\mathbb{X}$ at this flag is the $\mathbb{P}^3$ of special quadrics
$$g=b_0x_0^2 + b_1x_0x_1 + b_2x_0x_2 + b_3x_1^2.$$
Now, we will pay attention to the dense
open set where $b_2=1$, so 
		\begin{equation} \label{gfibra}
	g=b_0x_0^2 + b_1x_0x_1 + x_0x_2 + b_3x_1^2.
	\end{equation}
The divisibility condition (\ref{divcond1}) shows that the fiber of $\mathbb{Y}$ over a quadric as in (\ref{gfibra}) is the $\mathbb{P}^4=\{(a_0:a_1:a_2:a_3:a_5) \}$ of special cubics
$$f=(a_0x_0+a_1x_1 +a_2x_2+a_3x_3)x_0^2
+ a_5\left(b_1x_0x_1^2 +    x_0x_1x_2 +
\frac{2}{3}   b_3x_1^3\right)
.$$
Computing $\omega = \dfrac{3fdg-2gdf}{x_0}$ one can see that there are no indeterminacies on this neighborhood $[b_2=1]$ (cf.\,\ref{b2=1}). Moreover,
	$$2g\dfrac{\partial f}{\partial x_3} - 3f\dfrac{\partial g}{\partial x_3} = 2a_3gx_0^2,$$ 
	and the quadric $g$ can be recovered from
        $\omega $ just by looking at the coefficient of $dx_3$ (at least on the dense open set $a_3=1$).
	
	The coefficient of $dx_2$ in $\omega$ is
	
	$$\ba c
\dfrac{1}{x_0}\left(2g\dfrac{\partial f}{\partial x_2} -
  3f\dfrac{\partial g}{\partial x_2}\right) = 
\\\textcolor{blue}{(2a_2b_0-3a_0)}x_0^3+\textcolor{blue}{(2a_5b_0+2a_2b_1-3a_1)}x_0^2x_1+\textcolor{blue}{(-a_5b_1+2a_2b_3)}x_0x_1^2
\\\textcolor{blue}{-a_2}x_0^2x_2\textcolor{blue}{-a_5}x_0x_1x_2\textcolor{blue}{-3a_3}x_0^2x_3.
\ea$$
		
			Hence the coefficients $a_2$ and $a_5$ can be recovered from $x_0^2x_2dx_2$ and $x_0x_1x_2dx_2$, \linebreak
			respectively. And since the coefficients $b_0$ and $b_1$ are also known one can obtain from $\omega$ the coefficients $a_0$ and $a_1$ from $x_0^3dx_2$ and $x_0^2x_1dx_2$ 
			respectively. 
\end{proof} 
\end{prop}

\begin{teo} \label{integral}
	The degree of the exceptional component of
        codimension one and degree two foliations
        in $\mathbb{P}^3$ is given by
	$$\int\limits_{\mathbb{Y}_4}-c_1^{13}(\mathcal{W}) \cap [\mathbb{Y}_4].$$
	\begin{proof}
		We have $\dim (\mathbb{Y}_4) =13$. Since the map $\Phi$ is generically injective the required degree is (cf. definition of deg$_f\tilde{X}$ in \cite{Fulton}, page 83) $$\int\limits_{\mathbb{Y}_4}c_1(\Phi^*\mathcal{O}(1))^{13} \cap [\mathbb{Y}_4] = \int\limits_{\mathbb{Y}_4}-c_1(\Phi^* \mathcal{O}(-1))^{13} \cap [\mathbb{Y}_4]= \int\limits_{\mathbb{Y}_4}-c_1(\mathcal{W})^{13} \cap [\mathbb{Y}_4].$$
	\end{proof}
\end{teo}

\chapter{Calculation of the Degree} \label{chapgrau}

\section{Action of a Torus}

In order to obtain the value of the integral in Theorem \ref{integral}, we will apply Bott's formula
\begin{equation}  \label{BOTT}
\int\limits_{\mathbb{Y}_4}-c_1^{13}(\mathcal{W}) \cap [\mathbb{Y}_4]=\sum\limits_{F}\dfrac{-c_1^T(\mathcal{W})^{13} \cap [F]_T}{c_{top}^T(\mathcal{N}_{F|\mathbb{Y}_4})},\end{equation}
where the sum runs through all fixed components $F$ under
a convenient action of the torus $\mathbb{C}^*$. 
Fix  an action 
\begin{equation}\label{actionp3}
\begin{array}{ccc}
\mathbb{C}^* \times \pd3 & \longrightarrow & \pd3 \\
(t,x_i) & \mapsto & t^{w_i}x_i
\end{array}
\end{equation}
of $\mathbb{C}^*$ in $\pd3$ with distinct weights $w_i$, $i \in \{0,1,2,3 \}$.

For a good choice of these weights, the only fixed flags are the standard ones (cf. \ref{acaoquadcub}),
$$p_{ijk}=\{x_i=x_j=x_k=0\} \in \ell_{ij}=\{x_i=x_j=0\} \subset v_i=\{x_i=0\}.$$
Hence, we have 24 fixed flags.

\noindent
The reader not familiarized with Bott's formula
may consult 
%, we strongly recommend 
the  reference \cite{meirelles}, in special the examples
with a fixed component of positive dimension. See also
the appendix \ref{3planos}.

\section{Fixed Points over a Flag}

Now, consider our favourite fixed flag $\varphi _0=(p_{0},\ell_{0},v_0)$ (\ref{pdv}).

The induced action on $\mathbb{X}'$ is given by
\begin{equation}\label{acaoemg}
\left\{
\ba l
g=b_0x_0^2 + b_1x_0x_1 + b_2x_0x_2 + b_3x_1^2 
\ \mbox{\large$\mapsto$}  
\\\na7
t\circ g=
t^{2w_0}b_0x_0^2 + t^{w_0+w_1}b_1x_0x_1 + t^{w_0+w_2}b_2x_0x_2 + t^{2w_1}b_3x_1^2 \ ;\\\na{12}
g'=u_1x_0x_1 + u_2x_0x_2 + u_3x_1^2 \ \mbox{\large$\mapsto$}  
\\\na7
t\circ g'=
t^{w_0+w_1}u_1x_0x_1 + t^{w_0+w_2}u_2x_0x_2 +
t^{2w_1}u_3x_1^2\ .
\ea\right.
\end{equation}

We recognize at once the six isolated fixed points
$$
\mb x_1':=(x_0x_1), \
\mb x_2':=(x_0x_2), \
\mb x_3':=(x_1^2), \
\mb x_4':=(x_0^2,x_0x_1), \
\mb x_5':=(x_0^2,x_0x_2), \
\mb x_6':=(x_0^2,x_1^2).
$$
The last three points lie in the exceptional \p2 of
the blowup of $\p3=\{(b_0:\cdots:b_3)\}$
at the point $x_0^2\leftrightarrow(1:0:0:0)$.

Over each of these six points, in the fiber $\bb Y_{\mb
  x_i'}$ 
of $\mathbb{Y}$ we have 5 other  isolated fixed points, since this fiber is a $\mathbb{P}^4$ (again, for a good choice of the weights, cf. \ref{acaoqc}).

For example, take the fixed point $\mb x_1'=x_0x_1$. The fiber $\bb{Y}_{\mb x_1'}$ is the $\mathbb{P}^4$ of special cubics
$$f=a_0x_0^3 + a_1x_0^2x_1 + a_2x_0^2x_2+a_3x_0^2x_3 +
a_4x_0x_1^2,$$
by the divisibility condition (\ref{divcond1}) applied to (\ref{ff}). Thus, the five fixed cubics here are
$$
\mb y_1':=x_0^3, \
\mb y_2':=x_0^2x_1, \
\mb y_3':=x_0^2x_2, \
\mb y_4':=x_0^2x_3, \
\mb y_5':=x_0x_1^2. 
$$

For the pair $(\mb x_1', \mb y_3')=(x_0x_1,x_0^2x_2) \in \bb Y$, one can compute 
$$\omega = \dfrac{3\mb x_1'd\mb y_3' - 2\mb y_3' d\mb x_1'}{x_0} = -x_0x_1x_2dx_0 + 3x_0^2x_2dx_1 - 2x_0^2x_1dx_2.$$

The fiber of the line bundle $\cl W$ over the fixed point $(x_0x_1,x_0^2x_2)$ is the $\mathbb{C}$--linear space spanned by
$$\omega = x_0x_1x_2dx_0 - 3x_0^2x_2dx_1 + 2x_0^2x_1dx_2.$$ 

Notice that presently there are many fixed
points outside of the indeterminacy locus of the map $(g,f) \mapsto
\omega$ (see Table \ref{Fixosem} below). Their
contribution for Bott's formula can be immediately
computed.

We present a table of the 30 fixed points described until
here. Computations in \ref{fibrasY}.

\newpage

\begin{table}[h!]
	\centering
	\caption{Fixed points over $\mathbb{Y}$}
	\label{Fixosem} 
$\begin{array}{ccccc} \hline \hline \\
\textrm{fixed point at} & & \textrm{fixed point at} & & \textrm{generator of fiber} \\
\mathbb{X}':\mb{x}_i' & & \mathbb{Y}_{\mb{x}_i}:\mb{y}_j' & & \textrm{in } \mathcal{W} \\
\hline
&   &  &  & \\
& \vline & x_0^3    & \leadsto & \omega = x_0^2x_1dx_0 - x_0^3dx_1 \\
& \vline & x_0^2x_1 & \leadsto & \omega = x_0x_1^2dx_0 - x_0^2x_1dx_1  \\
x_0x_1 & \vline & x_0^2x_2 & \leadsto & \omega = x_0x_1x_2dx_0 - 3x_0^2x_2dx_1 + 2x_0^2x_1dx_2  \\
& \vline & x_0^2x_3 & \leadsto & \omega = x_0x_1x_3dx_0 - 3x_0^2x_3dx_1 + 2x_0^2x_1dx_3   \\
& \vline & x_0x_1^2 & \leadsto & \omega = x_1^3dx_0 - x_0x_1^2dx_1 \\
&        &     &        &   \\
& \vline & x_0^3    & \leadsto & \omega = x_0^2x_2dx_0 - x_0^3dx_2 \\
& \vline & x_0^2x_1 & \leadsto & \omega = x_0x_1x_2dx_0 + 2x_0^2x_2dx_1 - 3x_0^2x_1dx_2  \\
x_0x_2 & \vline & x_0^2x_2 & \leadsto & \omega = x_0x_2^2dx_0 -          x_0^2x_2dx_2 \\
& \vline & x_0^2x_3 & \leadsto & \omega = x_0x_2x_3dx_0 - 3x_0^2x_3dx_2 + 2x_0^2x_2dx_3   \\
& \vline & x_0x_1x_2 & \leadsto & \omega = -x_1x_2^2dx_0 + 2x_0x_2^2dx_1 - x_0x_1x_2dx_2 \\
&        &     &        &   \\
& \vline & x_0^3    & \leadsto & \omega = x_0x_1^2dx_0 - x_0^2x_1dx_1 \\
& \vline & x_0^2x_1 & \leadsto & \omega = x_1^3dx_0 - x_0x_1^2dx_1 \\
x_1^2  & \vline & x_0^2x_2 & \leadsto & \omega = 2x_1^2x_2dx_0 -          3x_0x_1x_2dx_1 + x_0x_1^2dx_2 \\
& \vline & x_0^2x_3 & \leadsto & \omega = 2x_1^2x_3dx_0 - 3x_0x_1x_3dx_1 + x_0x_1^2dx_3   \\
& \vline & x_1^3 & \leadsto & \textcolor{red}{\textrm{not defined}} \\
&        &     &        &   \\ 
& \vline & x_0^3    & \leadsto & \textcolor{red}{\textrm{not defined}} \\
& \vline & x_0^2x_1 & \leadsto & \omega = x_0^2x_1dx_0 - x_0^3dx_1 \\
x_0^2, x_0x_1  & \vline & x_0^2x_2 & \leadsto & \omega =     x_0^2x_2dx_0 - x_0^3dx_2 \\
& \vline & x_0^2x_3 & \leadsto & \omega = x_0^2x_3dx_0 - x_0^3dx_3   \\
& \vline & x_0x_1^2 & \leadsto & \omega = x_0x_1^2dx_0 - x_0^2x_1dx_1 \\
&        &     &        &   \\
& \vline & x_0^3    & \leadsto & \textcolor{red}{\textrm{not defined}} \\
& \vline & x_0^2x_1 & \leadsto & \omega = x_0^2x_1dx_0 - x_0^3dx_1 \\
x_0^2, x_0x_2  & \vline & x_0^2x_2 & \leadsto & \omega =         x_0^2x_2dx_0 - x_0^3dx_2 \\
& \vline & x_0^2x_3 & \leadsto & \omega = x_0^2x_3dx_0 - x_0^3dx_3   \\
& \vline & x_0x_1x_2 & \leadsto & \omega = 2x_0x_1x_2dx_0 - x_0^2x_2dx_1 - x_0^2x_1dx_2 \\
&        &     &        &   \\
& \vline & x_0^3    & \leadsto & \textcolor{red}{\textrm{not defined}} \\
& \vline & x_0^2x_1 & \leadsto & \omega = x_0^2x_1dx_0 - x_0^3dx_1 \\
x_0^2, x_1^2  & \vline & x_0^2x_2 & \leadsto & \omega =           x_0^2x_2dx_0 - x_0^3dx_2 \\
& \vline & x_0^2x_3 & \leadsto & \omega = x_0^2x_3dx_0 - x_0^3dx_3   \\
&\vline & x_1^3 & \leadsto & \omega = x_1^3dx_0 - x_0x_1^2dx_1               
\end{array}$
\end{table}

\newpage

There are 26 %well resolved 
fixed points in $\mathbb{Y}$ off the indeterminacy
locus. In the following examples we will compute the
contribution of the fixed points for Bott's formula 
(\ref{BOTT}).
 
\begin{exem}\normalfont \label{fx1Y}
	Let us compute the contribution of the point $(g,f)=(x_0x_1,x_0^2x_2)$. We have 
	$$\omega=x_0x_1x_2dx_0 - 3x_0^2x_2dx_1 + 2x_0^2x_1dx_2,
	$$
	hence the induced action reads 
	$$\ba{cl}
	t \ \circ \ \omega &= t^{w_0}x_0t^{w_1}x_1t^{w_2}x_2d(t^{w_0}x_0) - 3t^{2w_0}x_0^2t^{w_2}x_2d(t^{w_1}x_1) + 2t^{2w_0}x_0^2t^{w_1}x_1d(t^{w_2}x_2)\\\na8
	&=t^{2w_0+w_1+w_2} \cdot (x_0x_1x_2dx_0 - 3x_0^2x_2dx_1 +
	2x_0^2x_1dx_2),
	\ea
	$$
	and so the equivariant Chern class is
	\begin{equation} \label{c1ww}
	c_1^T(\mathcal{W})_{|_{(\mb{x}_1' ,\mb{y}_3')}} = 2w_0+w_1+w_2.
	\end{equation}
	
	Since this component of fixed points consists of a single
	point, its normal space \linebreak coincides with the tangent space of $\mathbb{Y}$ at this point. It decomposes as
	\begin{equation} \label{tangband}
	T_{(x_0x_1,x_0^2x_2)}\mathbb{Y}=T_{(p,\ell,v)}{\bb F} \ \oplus \ T_{(x_0x_1)}\mathbb{X} \ \oplus \ T_{(x_0^2x_2)}\mathbb{Y}_{x_0x_1}=
	\end{equation}
	$$=\left( \dfrac{x_1}{x_0} + \dfrac{x_2}{x_0} + \dfrac{x_3}{x_0} + \dfrac{x_2}{x_1} + \dfrac{x_3}{x_1} + \dfrac{x_3}{x_2} \right) + \left( \dfrac{x_0^2}{x_0x_1} + \dfrac{x_0x_2}{x_0x_1} + \dfrac{x_1^2}{x_0x_1} \right) + $$ $$ + \left( \dfrac{x_0^3}{x_0^2x_2} + \dfrac{x_0^2x_1}{x_0^2x_2} + \dfrac{x_0^2x_3}{x_0^2x_2} + \dfrac{x_0x_1^2}{x_0^2x_2} \right),$$
	where this sum represents the decomposition into
	eigenspaces; the weight of $\dfrac{x_i}{x_j}$ is
	$w_i-w_j$. From this, it follows that the contribution of
	the point $(x_0x_1,x_0^2x_2)$ in  Bott's formula is the fraction 
	\begin{equation} \label{contrx}
	\dfrac{-(2w_0+w_1+w_2)^{13}}{(w_0-w_1)^3(w_2-w_0)^2(w_3-w_0)(w_2-w_1)^3(w_3-w_1)(w_3-w_2)^2(2w_1-w_0-w_2)}.
	\end{equation} $\square$
\end{exem} 

Let us illustrate the general process used to compute the tangent spaces in (\ref{tangband}). For the fiber $\mathbb{Y}_{x_0x_1}$, since $b_2=b_3=0$ and $b_1=1$, the divisibility condition (\ref{divcond1}) reads $a_5=a_6=0$, so this fiber consists of special cubics 
$$f=a_0x_0^3+a_1x_0^2x_1+a_2x_0^2x_2+a_3x_0^2x_3+a_4x_0x_1^2.$$
To compute the tangent space of this
$\mathbb{P}^4$--fiber over the fixed point $x_0^2x_2$
(\ref{QxS}), we  write
\begin{equation}\label{TY01002}
T_{x_0^2x_2}\mathbb{Y}_{x_0x_1} =  \underbrace{\langle x_0^3,\ x_0^2x_1,\ x_0^2x_3,\ x_0x_1^2  \rangle}_\text{other generators of the fiber} \otimes \underbrace{\langle x_0^2x_2 \rangle ^\vee}_\text{point}=
\end{equation}
$$= [x_0^3\otimes(x_0^2x_2)^\vee] \oplus [(x_0^2x_1)\otimes (x_0^2x_2)^\vee] \oplus [(x_0^2x_3)\otimes (x_0^2x_2)^\vee] \oplus [(x_0x_1^2)\otimes (x_0^2x_2)^\vee] = 
$$
$$\stackrel{\text{notation}}{=} \dfrac{x_0^3}{x_0^2x_2} + \dfrac{x_0^2x_1}{x_0^2x_2} + \dfrac{x_0^2x_3}{x_0^2x_2} + \dfrac{x_0x_1^2}{x_0^2x_2}. $$\\

\begin{exem}\normalfont \label{fx2Y}
	Next, take a fixed point over the exceptional divisor of $\mathbb{X}'$. Say, take the fixed point $(x_0^2,x_0x_2,x_0^2x_3)$. The tangent space of $\mathbb{Y}$ at this point decomposes as
	$$T_{(x_0^2,x_0x_2,x_0^2x_3)}\mathbb{Y}=T_{(p,\ell,v)}{\bb F} \ \oplus \ T_{(x_0^2,x_0x_2)}\mathbb{X}' \ \oplus \ T_{(x_0^2x_3)}\mathbb{Y}_{x_0^2,x_0x_2}$$
	and the tangent space of $\mathbb{X}'$ decomposes as
	\begin{equation}\label{T0002} T_{x_0^2,x_0x_2}\mathbb{X}' = \mathcal{L}_{x_0^2,x_0x_2} \ \oplus \ T_{[\mathcal{L}_{x_0^2,x_0x_2}]}\mathbb{P}(\mathcal{N}_{x_0^2|\mathbb{X}}),\end{equation}
	where $\mathcal{L}_{x_0^2,x_0x_2}$ is the line represented by the point $(x_0^2,x_0x_2)$ in $\mathbb{P}(\mathcal{N}_{x_0^2|\mathbb{X}})$.
	
	Since $x_0^2$ is a point in $\mathbb{X}$, we have 
	$$\mathcal{N}_{x_0^2|\mathbb{X}} = T_{x_0^2}\mathbb{X} = 
	\dfrac{x_0x_1}{x_0^2} + 
	\underbracket{\dfrac{x_0x_2}{x_0^2}}_{\mathcal{L}_{x_0^2,x_0x_2}} + \dfrac{x_1^2}{x_0^2}$$
	and 
	$$T_{[\mathcal{L}_{x_0^2,x_0x_2}]}\mathbb{P}(\mathcal{N}_{x_0^2|\mathbb{X}}) = \left( \dfrac{x_0x_1}{x_0^2} \cdot \dfrac{x_0^2}{x_0x_2} \right) + \left( \dfrac{x_1^2}{x_0^2} \cdot \dfrac{x_0^2}{x_0x_2} \right) = \dfrac{x_0x_1}{x_0x_2} + \dfrac{x_1^2}{x_0x_2} .$$
	
	The three eigenspaces in (\ref{T0002}) are the summands
	as displayed below:
	\begin{equation}\label{TX'02}
	T_{x_0^2,x_0x_2}\mathbb{X}' = \left( \dfrac{x_0x_2}{x_0^2} \right) + \left( \dfrac{x_0x_1}{x_0x_2} + \dfrac{x_1^2}{x_0x_2} \right).
	\end{equation}
	$\square$
\end{exem}

In a similar way of examples \ref{fx1Y} and \ref{fx2Y} we can compute the contributions for the 26 ``well resolved'' points of Table \ref{Fixosem}.

\section{Degree of a fiber}

At any point $q  \in \mathbb{Y}_4$
lying over $(p,\ell,v)\in\bb F$, 
the tangent space of $\mathbb{Y}_4$  decomposes  as
$$T_{(p,\ell,v,q)}\mathbb{Y}_4 = T_{(p,\ell,v)}{\bb F} \oplus
T_q{\mathbb{Y}_4}_{(p,\ell,v)},
$$ 
where $T_q{\mathbb{Y}_4}_{(p,\ell,v)}$ means the tangent
space of the fiber of $\mathbb{Y}_4$ over the flag
$(p,\ell,v) \in {\bb F}$. Now, since we have fixed the flag
$(p_{0},\ell_{0},v_0)$ we will omit the 
summand
$T_{(p,\ell,v)}{\bb F}$ in the following sections, and focus
on the rank 7 space $T_q{\mathbb{Y}_4}_{(p_{0},\ell_{0},v_0)} = T_q\mathbb{Y}_4$.

Note that over a fixed flag, the fiber
${\mathbb{Y}_4}_{(p,\ell,v)}$ is a smooth projective variety of
dimension 7. This gives a subvariety of exceptional
foliations which are obtained from a fixed flag. We can
ask about the degree of this 7-dimensional variety, and
the answer can be expressed as the integral (in a same way as Theorem \ref{integral})
$$\int\limits_{{\mathbb{Y}_4}_{(p,\ell,v)}}-c_1^{7}(\mathcal{W}) \cap [{\mathbb{Y}_4}_{(p,\ell,v)}].$$

From now on we fix the flag $\varphi _0$ in (\ref{pdv}) to compute the
degree of the image of that fiber. This procedure reduces the affine coordinates to 7 instead of 13, simplifying the computational aspects.

%Note that, in order to find the degree of a fiber, 
%at first we could
%try to use Aluffi's algorithm, described in
%\cite{AluffiSegre}, to compute the Segre class of the
%indeterminacy locus, as done in some cases in
%\cite{israel1}, or yet the algorithms in \cite{Helmer} to
%compute the Segre class. This could be done, in
%principle, but we have not tried yet.
%replacing $\mathbb{P}^3 \times \mathbb{P}^6$ by 
%But the basis is a subscheme of  which is described by the divisibility condition, a complete intersection only locally. This turns the actual computations 
%inviable with this procedure.
%
\section{Contributions over the blowup centers}

In $\mathbb{Y}$ there exist 30 fixed points (Table \ref{Fixosem});  for 26 of
them we may compute immediately the \linebreak contributions
as shown in (\ref{contrx}). There
are still four remaining points to consider.   

\subsection {$(g,f) = (x_1^2,x_1^3).$} \label{primpntfx}

%Let's start taking the point $g=x_1^2, \ f=x_1^3$. 
Notice
that this point lives in the  neighborhood $a_6=b_3=1$. See \ref{a6b3=1txt}.

First, we do the blowup of $\mathbb{Y}$ along the curve $C$ (\ref{expC})
$$ f=\left(\dfrac{b_1}{2}x_0 + x_1 \right)^3, g=\left(\dfrac{b_1}{2}x_0 + x_1 \right)^2.$$

It's  a local complete intersection. We will look at the local equations, and use a trick as explained below to give adequate weights to the local coefficients of a general quadric and cubic. 

From (\ref{acaoemg}), the $\bb{P}^3=\{(b_0:b_1:b_2:b_3)\}$ of special quadrics has action 
\begin{equation}\label{acaonosb}
[t,(b_0:b_1:b_2:b_3)] \mapsto (t^{2w_0}b_0:t^{w_0+w_1}b_1:t^{w_0+w_2}b_2:t^{2w_1}b_3)
\end{equation}

The trick is to look at the affine chart $b_3=1$,
$$(\mb{b}_0,\mb{b}_1,\mb{b}_2)=\left(\dfrac{b_0}{b_3},\dfrac{b_1}{b_3},\dfrac{b_2}{b_3} \right)$$ 
and set
$$[t,(\mb{b}_0,\mb{b}_1,\mb{b}_2)]\mapsto \left(\dfrac{t^{2w_0}b_0}{t^{2w_1}b_3},\dfrac{t^{w_0+w_1}b_1}{t^{2w_1}b_3},\dfrac{t^{w_0+w_2}b_2}{t^{2w_1}b_3} \right)$$

Thus, we gain naturally local weights for coordinates around the fixed point $x_1^2$:
\begin{equation}\label{pesosbx12}
w_{\mb{b}_0} = 2w_0-2w_1, \ \ \ w_{\mb{b}_1} = w_0-w_1, \ \ \ w_{\mb{b}_2} = w_0+w_2-2w_1
\end{equation}

These are weights for the \textit{points} of $\mathbb{A}^3$. For the respective linear functions the local weights have opposite sign, {\em i.e.,} the linear function $b_1$ has local weight $-w_{\mb{b}_1}=w_1-w_0$, and the linear function $b_2$ has local weight $-w_{\mb{b}_2}=2w_1-w_0-w_2$. Then,
$$t \circ (b_1x_0x_1^2) = t^{w_1-w_0}b_1t^{w_0+2w_1}x_0x_1^2 = t^{3w_1}b_1x_0x_1^2 \ ,$$
$$t \circ (b_2x_0x_1x_2) = t^{2w_1-w_0-w_2}b_2t^{w_0+w_1+w_2}x_0x_1x_2 = t^{3w_1}b_2x_0x_1x_2 \ .$$

The special cubics are (\ref{fga6b3=1})
$$f=a_0x_0^3 + a_1x_0^2x_1 + a_2x_0^2x_2+a_3x_0^2x_3 +
a_6\left(\frac{3}{2}b_1x_0x_1^2+\frac{3}{2}
b_2x_0x_1x_2 +   x_1^3\right).$$

The $\bb{P}^4=\{(a_0,a_1:a_2:a_3:a_6)\}$ of special
cubics  inherits the induced 
action 
\begin{equation}\label{acaonosa}
[t,(a_0:a_1:a_2:a_3:a_6)] \mapsto (t^{3w_0}a_0:t^{2w_0+w_1}a_1:t^{2w_0+w_2}a_2:t^{2w_0+w_3}a_3:t^{3w_1}a_6)
\end{equation}

In the affine chart $a_6=1$ we write
$$(\mb{a}_0,\mb{a}_1,\mb{a}_2,\mb{a}_3)=\left(\dfrac{a_0}{a_6},\dfrac{a_1}{a_6},\dfrac{a_2}{a_6},\dfrac{a_3}{a_6} \right)$$ 
and set
$$[t,(\mb{a}_0,\mb{a}_1,\mb{a}_2,\mb{a}_3)]\mapsto \left(\dfrac{t^{3w_0}a_0}{t^{3w_1}a_6},\dfrac{t^{2w_0+w_1}a_1}{t^{3w_1}a_6},\dfrac{t^{2w_0+w_2}a_2}{t^{3w_1}a_6},\dfrac{t^{2w_0+w_3}a_3}{t^{3w_1}a_6} \right)$$
\medskip

In this way we can set local weights around the fixed point $(x_1^2,x_1^3)$

$$\begin{array}{rclcrcl}
w_{\mb{a}_0} & = & 3w_0-3w_1 ,    & & w_{\mb{b}_0}& = & 2w_0-2w_1,\\
w_{\mb{a}_1} & = & 2w_0-2w_1 ,    & & w_{\mb{b}_1}& = & w_0-w_1,\\
w_{\mb{a}_2} & = & 2w_0-3w_1+w_2 ,& & w_{\mb{b}_2}& = & w_0-2w_1+w_2,\\
w_{\mb{a}_3} & = & 2w_0-3w_1+w_3. & &        &   & 
\end{array}$$

%In the affine chart $a_6=b_3=1$ we have
%$$
%\left\{\ba l
%f=a_0x_0^3 + a_1x_0^2x_1 + a_2x_0^2x_2+a_3x_0^2x_3 + %\frac{3}{2}b_1x_0x_1^2+\frac{3}{2}b_2x_0x_1x_2 + x_1^3,\\\na6
%g=b_0x_0^2 + b_1x_0x_1 + b_2x_0x_2 + x_1^2.
%\ea\right.$$
%We assign weights $w_{a_i},w_{b_j}$ to the coefficients
%so that $g,f$ are homogeneous. 
%Observe that $f$ must have the weighted degree
%$3w_1$, due  to the presence of the term $x_1^3$.
%Hence, the term $a_0x_0^3$ will receive this weight $3w_1$. But $x_0^3$ has weight $3w_0$, and so
%$$w_{a_0} + 3w_0 = 3w_1 \ \Rightarrow \ w_{a_0}=3w_1-3w_0.$$
%
%In the same way, a comparison with the term 
%$a_1x_0^2x_1$ reveals weights
%$$w_{a_1} + 2w_0 + w_1 \ = \ 3w_1 \Rightarrow w_{a_1}=2w_1-2w_0.$$
%
%Repeating this argument we gain naturally local weights
%for coordinates around the  fixed point $(x_1^2,x_1^3)$
%in the  indeterminacy locus: 
%$$\begin{array}{rclcrcl}
%w_{a_0} & = & 3w_1-3w_0 ,    & & w_{b_0}& = & 2w_1-2w_0,\\
%w_{a_1} & = & 2w_1-2w_0 ,    & & w_{b_1}& = & w_1-w_0,\\
%w_{a_2} & = & 3w_1-2w_0-w_2 ,& & w_{b_2}& = & 2w_1-w_0-w_2,\\
%w_{a_3} & = & 3w_1-2w_0-w_3. & &        &   & 
%\end{array}$$

The local equations of the first blowup center are (see ideal $J_{red}$ in (\ref{eqC})): 
$$\left\{\begin{array}{rcl}
8a_0-b_1^3 & & (\textrm{weight }\ 3w_1-3w_0) \\
a_2 & & (\textrm{weight }\ 3w_1-2w_0-w_2) \\
b_2 & & (\textrm{weight }\ 2w_1-w_0-w_2) \\
a_3 & & (\textrm{weight }\ 3w_1-2w_0-w_3) \\
4b_0-b_1^2 & & (\textrm{weight }\ 2w_1-2w_0) \\
4a_1-3b_1^2 & & (\textrm{weight }\ 2w_1-2w_0)
\end{array}\right.
$$

Notice that the above are all weighted homogeneous
equations. The first blowup produces as exceptional
divisor a $\mathbb{P}^5$-bundle over the curve
$C$, to wit,  the projectivized normal bundle
$\mathbb{P}(\mathcal{N}_{C|\mathbb{Y}})$. Since
$J_{red}/J_{red}^2$ is the conormal sheaf of $C$,
the weights of the normal bundle have opposite
sign by duality. Then, over a fixed point in $C$, we have a (\ref{Y1a6b3})
\begin{equation*} \label{signchg}
\mathbb{P}^5=\{(s_0:s_1:s_2:s_3:s_4:s_5) \}
\end{equation*}
 with a natural induced action and weights
$$\begin{array}{rclcrcl}
w_{s_0} & = & 3w_0-3w_1     & & w_{s_3}& = & 2w_0+w_3-3w_1\\
w_{s_1} & = & 2w_0+w_2-3w_1 & & w_{s_4}& = & 2w_0-2w_1\\
w_{s_2} & = & w_0+w_2-2w_1  & & w_{s_5}& = & 2w_0-2w_1.
\end{array}$$ 
\begin{equation} \label{acaoind}
t \circ (s_0:s_1:s_2:s_3:s_4:s_5) = (t^{w_{s_0}}s_0 : t^{w_{s_1}}s_1 : t^{w_{s_2}}s_2 : t^{w_{s_3}}s_3 : t^{w_{s_4}}s_4 : t^{w_{s_5}}s_5)
\end{equation} 

Over the fixed point $q:=(x_1^2,x_1^3) \in C$ we have a $\mathbb{P}^5$ with four isolated fixed points, and since $w_{s_4}=w_{s_5}$, a fixed one dimensional component 
\begin{equation}\label{retal1}
\ell_1=\mathbb{P}^1=\{(0:0:0:0:s_4:s_5)\}\ .
\end{equation}

Let us describe the fiber of the normal bundle over $q$. First of all, we have
\begin{equation}\label{Tx12x13}
T_q\mathbb{Y} = T_{(x_1^2)}\bb{X} \oplus T_{(x_1^3)}\bb{Y}_{x_1^2}
\end{equation}

At the point $g=x_1^2$ we have $b_1=b_2=0$, so the fiber $\bb{Y}_{x_1^2}$ is the $\bb{P}^4$ of special cubics  
$$f=a_0x_0^3 + a_1x_0^2x_1 + a_2x_0^2x_2+a_3x_0^2x_3 +
a_6x_1^3$$
and (\ref{Tx12x13}) reads
$$T_q\mathbb{Y} = \left(\dfrac{x_0^2}{x_1^2} + \dfrac{x_0}{x_1} + \dfrac{x_0x_2}{x_1^2}\right) + \left(\dfrac{x_0^3}{x_1^3} + \dfrac{x_0^2}{x_1^2} + \dfrac{x_0^2x_2}{x_1^3} + \dfrac{x_0^2x_3}{x_1^3}\right) \ . $$\medskip

The blowup center is the curve $C$ parametrized by $sx_0+x_1$, so (as in (\ref{TY01002}))
$$T_qC = \dfrac{x_0}{x_1}.$$ 

From the short exact sequence
$$T_qC \hookrightarrow T_q\mathbb{Y} \twoheadrightarrow \mathcal{N}_{qC|\mathbb{Y}},$$

one can find
\begin{equation} \label{NqCY}
\mathcal{N}_{qC|\mathbb{Y}} = \underbrace{\dfrac{x_0^2}{x_1^2}}_\text{{$s_5'$}} + \underbrace{\dfrac{x_0x_2}{x_1^2}}_\text{{$s_2'$}} + \underbrace{\dfrac{x_0^3}{x_1^3}}_\text{{$s_0'$}} + 
\underbrace{\dfrac{x_0^2}{x_1^2}}_\text{{$s_4'$}} + \underbrace{\dfrac{x_0^2x_2}{x_1^3}}_\text{{$s_1'$}} + \underbrace{\dfrac{x_0^2x_3}{x_1^3}}_\text{{$s_3'$}}
\end{equation}

Note the presence of the two summands $\dfrac{x_0^2}{x_1^2}$, representing two (independent) \linebreak eigenvectors associated to the same weight. A simple comparison of the weights \linebreak reveals the associations
with the homogeneous coordinates $s_i$; for example, the
fixed point $s_0'=(1:0:0:0:0:0)$, with weight $w_{s_0} =
3w_0-3w_1$ (cf (\ref{acaoind})) is associated to the eigenvector $\dfrac{x_0^3}{x_1^3}$, with weight $3w_0-3w_1$.\\

For short, we use the following notation for points of the $\mathbb{P}^5$-fiber:
$$s_0'=(1:0:0:0:0:0), \ \ \ s_1'=(0:1:0:0:0:0)$$
$$s_2'=(0:0:1:0:0:0), \ \ \ s_3'=(0:0:0:1:0:0)$$
$$s_4'=(0:0:0:0:1:0), \ \ \ s_5'=(0:0:0:0:0:1)$$

\begin{exem}\normalfont \label{pontox13x12s0}
	To compute the contribution of the fixed point 
	$$(q,s_0') = ((x_1^2,x_1^3),(1:0:0:0:0:0)) \in \mathbb{Y}_1,$$ write
	$$T_{(q,s_0')}\mathbb{Y}_1 = T_qC \oplus \mathcal{L}_{x_0^3/x_1^3} \oplus T_{[\mathcal{L}_{x_0^3/x_1^3}]}\mathbb{P}(\mathcal{N}_{qC|\mathbb{Y}}).$$
		This leads us to the decomposition \\
		$$T_{(q,s_0')}\mathbb{Y}_1 = \left(\dfrac{x_0}{x_1}\right) + \left(\dfrac{x_0^3}{x_1^3}\right) + \left[\left(\dfrac{x_0^2}{x_1^2} + \dfrac{x_0^2x_2}{x_1^3} + \dfrac{x_0^2x_3}{x_1^3} + \dfrac{x_0^2}{x_1^2} + \dfrac{x_0x_2}{x_1^2}\right) \cdot \left(\dfrac{x_0^3}{x_1^3}\right)^* \ \right] = $$ 
	$$= \left(\dfrac{x_0}{x_1}\right) + \left(\dfrac{x_0^3}{x_1^3}\right) + \left(\dfrac{x_1}{x_0} + \dfrac{x_2}{x_0} + \dfrac{x_3}{x_0} + \dfrac{x_1}{x_0} + \dfrac{x_1x_2}{x_0^2}\right) \ . $$
	
	Thus, 
	$$c_7^T(T_{(q,s_0')}\bb{Y}_1) = 3(w_0-w_1)^4(w_2-w_0)(w_3-w_0)(w_1+w_2-2w_0) \ .$$
	
	Now, we look for the equivariant Chern class
	$c_1^T(\mathcal{W})$. To obtain 
	the fiber of $\mathcal{W}$ at the point $(q,s_0')$ perform the blowup,
	choosing the first equation (exc$:=8a_0-b_1^3$ (\ref{eqC})) to be the
	equation of the exceptional divisor. Step by step: (cf. \ref{othviz} and \ref{fibrasx13x12})
	\begin{enumerate}
		\item Write 
		$$\left\{\ba l
		f=a_0x_0^3 + a_1x_0^2x_1 + a_2x_0^2x_2+a_3x_0^2x_3 + \frac{3}{2}b_1x_0x_1^2+\frac{3}{2}b_2x_0x_1x_2 + x_1^3\\\na9
		g=b_0x_0^2 + b_1x_0x_1 + b_2x_0x_2 + x_1^2;
		\ea\right.
		$$
		\item Compute $$\omega = \dfrac{3fdg-2gdf}{x_0};$$ 
		\item Do the substitutions into $\omega$,
		$$\left\{\begin{array}{rcl}
		a_2        & = & s_1(8a_0-b_1^3) \\
		b_2        & = & s_2(8a_0-b_1^3) \\
		a_3        & = & s_3(8a_0-b_1^3) \\
		4b_0-b_1^2 & = & s_4(8a_0-b_1^3) \\
		4a_1-b_1^2 & = & s_5(8a_0-b_1^3)
		\end{array}\right.
		$$
		\item Divide by the equation exc=$8a_0-b_1^3$;
		\item Evaluate the new obtained expression at the
		fixed point, simply by taking the \linebreak coefficients 
		$$a_0=b_1=0 \ (\textrm{the point } (f=x_1^3,g=x_1^2)),$$ 
		$$s_1=s_2=s_3=s_4=s_5=0 \ (\textrm{the point } (1:0:0:0:0:0)).$$
	\end{enumerate}
	
	After the blowup, we can see that (cf.\,Fwn$\underline{\ \ }$0 in \ref{fibrasx13x12})
	$$\omega=x_0x_1^2dx_0-x_0^2x_1dx_1,$$
	so $c_1^T(\mathcal{W})|_{(q,s_0')} = 2w_0 + 2w_1$ (see (\ref{c1ww})).
	
	$_\square$
\end{exem}

However, we will not always find a well defined $\omega$ because there are indeterminacies remaining -- we don't solve all of them with this single blowup. But in example \ref{pontox13x12s0} we were lucky at the chosen point.

Next, we present a table with the four fixed points and the fixed $\mathbb{P}^1$: \ref{fibrasx13x12}

\newpage
\begin{table}[h!]
	\centering
	\caption{Fixed points over $f=x_1^3, g=x_1^2$}
	\label{x13x12}
\begin{tabular}{ccc} \hline \hline \\
\textrm{fixed point/ }$\mathbb{P}^1$ & associated & \textrm{generator of fiber} \\
$(s_0:s_1:s_2:s_3:s_4:s_5)$ & eigenvector & \textrm{in } $\mathcal{W}$ \\
\hline
& & \\
$(1:0:0:0:0:0)$ & $x_0^3 / x_1^3$ & $\omega=x_0x_1^2dx_0-x_0^2x_1dx_1$\\
$(0:1:0:0:0:0)$ & $x_0^2x_2 / x_1^3$ & $\omega=2x_1^2x_2dx_0-3x_0x_1x_2dx_1 + x_0x_1^2dx_2$\\
$\red{(0:0:1:0:0:0)}$ & $x_0x_2 / x_1^2$ & \textcolor{red}{\textrm{not defined}}\\
$(0:0:0:1:0:0)$ & $x_0^2x_3 / x_1^3$ & $\omega=2x_1^2x_3dx_0-3x_0x_1x_3dx_1 + x_0x_1^2dx_3$\\
\blu{$(0:0:0:0:s_4:s_5)$} & $x_0^2 / x_1^2$ & $\omega=\textcolor{blue}{(2s_5-3s_4)}x_1^3dx_0-\textcolor{blue}{(2s_5-3s_4)}x_0x_1^2dx_1$
\end{tabular}
\end{table}

There are three well resolved points over $q=(x_1^2,x_1^3)$: the points $s_0', \ s_1' $ and $s_3'$. Now we need to study the indeterminacy point 
\red{(0:0:1:0:0:0)} and the fixed line \blu{$\ell_1$}. 

\subsubsection{Resolution at the point $s_2'=(0:0:1:0:0:0).$}

The new indeterminacy locus is the ruled surface $R$, defined by the ideal (see (\ref{JR}))
\begin{equation}\label{JRs2}
J''=\langle 3s_4-2s_5,s_3,s_0-b_1s_5,4s_1-3b_1,b_2 \rangle .
\end{equation}

Over the point $s_2'=(0:0:1:0:0:0)$ we have a (\ref{RRa6b3}) $$\mathbb{P}^4=\{(v_0:v_1:v_2:v_3:v_4)\}=\mathbb{P}(\mathcal{N}_{(q,s_2')R|\mathbb{Y}_1}).$$ 
The action passes to this $\mathbb{P}^4$ by looking to the normal space, see (\ref{normalR1}). But we can use our trick  to see the induced action, just like it was done in \S\ref{primpntfx}.

Indeed, we are in the neighborhood $s_2=1$, and in this affine chart,
$$(\mb{s}_0,\mb{s}_1,\mb{s}_3,\mb{s}_4,\mb{s}_5)=\left(\dfrac{s_0}{s_2},\dfrac{s_1}{s_2},\dfrac{s_3}{s_2},\dfrac{s_4}{s_2},\dfrac{s_5}{s_2}\right).$$

Thus, we have local weights for affine points around $s_2=1$, given by (see (\ref{acaoind}))
$$\begin{array}{rccclcrcccl}
w_{s_0}^L & = & w_{s_0}-w_{s_2} & = & 2w_0-w_1-w_2  & & w_{s_3}^L & = & w_{s_3}-w_{s_2} & = & w_0+w_3-w_1-w_2\\
w_{s_1}^L & = & w_{s_1}-w_{s_2} & = & w_0-w_1       & & w_{s_4}^L & = & w_{s_4}-w_{s_2} & = & w_0-w_2\\
          &   &                 &   &               & & w_{s_5}^L & = & w_{s_5}-w_{s_2} & = & w_0-w_2
\end{array}$$ 

Remember that the equations have opposite sign weigths by duality. The local equations of the blowup center $R$ are
$$\left\{\begin{array}{rcl}
3\mb{s}_4-2\mb{s}_5 & & (\textrm{weight } w_2-w_0) \\
\mb{s}_3 & & (\textrm{weight } w_1+w_2-w_0-w_3) \\
\mb{s}_0-b_1\mb{s}_5 & & (\textrm{weight } w_1+w_2-2w_0) \\
4\mb{s}_1-3b_1 & & (\textrm{weight } w_1-w_0) \\
b_2 & & (\textrm{weight } 2w_1-w_0-w_2) 
\end{array}\right.
$$
and so relations (\ref{eqbup7}) induce

$$\begin{array}{rclcrcl}
w_{v_0} & = & w_0-w_2     & & w_{v_3}& = & w_0-w_1\\
w_{v_1} & = & w_0+w_3-w_1-w_2 & & w_{v_4}& = & w_0+w_2-2w_1\\
w_{v_2} & = & 2w_0-w_1-w_2  & & 
\end{array}$$ 
\begin{equation} \label{acaoind2}
t \circ (v_0:v_1:v_2:v_3:v_4) = (t^{w_{v_0}}v_0 : t^{w_{v_1}}v_1 : t^{w_{v_2}}v_2 : t^{w_{v_3}}v_3 : t^{w_{v_4}}v_4)
\end{equation} 

This $\mathbb{P}^4$ has five isolated fixed points
since the weights are distinct pairwise.
(cf.\,\ref{fibrasx13x12s2}).

\begin{table}[h!]
	\centering
	\caption{Fixed points over $(x_1^3,x_1^2), (0:0:1:0:0:0)$}
	\label{x13x12s2}
\begin{tabular}{ccc}  \hline \hline \\
\textrm{fixed point} & associated & \textrm{generator of fiber} \\
$(v_0:v_1:v_2:v_3:v_4)$ & eigenvector & \textrm{in } $\mathcal{W}$ \\
\hline
& & \\
$(1:0:0:0:0)$ & $x_0/x_2$ & $\omega=x_1^3dx_0-x_0x_1^2dx_1$\\
$(0:1:0:0:0)$ & $x_0x_3/x_1x_2$ & $\omega=2x_1^2x_3dx_0-3x_0x_1x_3dx_1 + x_0x_1^2dx_3$\\
$(0:0:1:0:0)$ & $x_0^2/x_1x_2$ & $\omega=x_0x_1^2dx_0 - x_0^2x_1dx_1$ \\
$(0:0:0:1:0)$ & $x_0/x_1$ & $\omega=2x_1^2x_2dx_0-3x_0x_1x_2dx_1 + x_0x_1^2dx_2$\\
$(0:0:0:0:1)$ & $x_0x_2/x_1^2$ & $\omega=x_1x_2^2dx_0-2x_0x_2^2dx_1+x_0x_1x_2dx_0$
\end{tabular}
\end{table}

Now we look at the tangent spaces at these points. We are able to write (see (\ref{NqCY}))
$$\mathcal{N}_{qC|\mathbb{Y}} = \dfrac{x_0^2}{x_1^2} + \dfrac{x_0x_2}{x_1^2} + \dfrac{x_0^3}{x_1^3} + \dfrac{x_0^2}{x_1^2} + \dfrac{x_0^2x_2}{x_1^3} + \dfrac{x_0^2x_3}{x_1^3} \ .$$ 

Remember that the indeterminacy point $s_2'$ corresponds to $\dfrac{x_0x_2}{x_1^2}$ (Table \ref{x13x12}). The blowup center is the ruled surface $R$, a $\mathbb{P}^1$-bundle over $C$. The tangent space of this surface $R$ decomposes  as
\begin{equation} \label{TRR}
T_{(q,s_2')}R = T_qC \oplus T_{s_2'}\mathbb{P}^1.
\end{equation}

The equations of this $\mathbb{P}^1$ are (\ref{JRs2})
$$\left\{
\ba c
3s_4-2s_5 = 0\\\noalign{\vskip5pt}
  s_3 = 0\\\noalign{\vskip5pt}
  s_0-b_1s_5 = 0\\\noalign{\vskip5pt}
  4s_1 - 3b_1 = 0,
\ea\right.
$$
hence this is the $\mathbb{P}^1$-fiber of points
$$(b_1s_5: \dfrac{3}{4}b_1:s_2:0:\dfrac{2}{3}s_5:s_5).$$  

At the fixed point $(x_1^3,x_1^2,s_2')$ we have $b_1=0$, so this is the \begin{equation} \label{P1exc}
\mathbb{P}^1=\{(0:0:s_2:0:\dfrac{2}{3}s_5:s_5) \},
\end{equation}
and 
\begin{equation} \label{Ts2}
T_{s_2'}\mathbb{P}^1 = \dfrac{s_5}{s_2} = \dfrac{x_0^2}{x_1^2} \cdot \dfrac{x_1^2}{x_0x_2} = \dfrac{x_0}{x_2}\ .
\end{equation} \medskip

A note for (\ref{Ts2}): at the $\mathbb{P}^1 = \{
(s_2:s_5) \}$ (\ref{P1exc}) we look at the point
$s_2'=(1:0)$. The equation of this point in
$\mathbb{P}^1$ is $s_5^{\vee}=0$, the linear form
vanishing at the point. From this, the tangent space of
$\mathbb{P}^1$ at this point is
$\dfrac{s_2^{\vee}}{s_5^{\vee}} = \dfrac{s_5}{s_2}$. At
this moment we ``overload'' the notation to distinguish %tell apart 
between the equation ($s_5^\vee$) and coordinate of the point $(s_5)$. \medskip

Using the expression (\ref{Ts2}) in (\ref{TRR}) gives
$$T_{(q,s_2')}R = \dfrac{x_0}{x_1} + \dfrac{x_0}{x_2}.
$$

 The short exact sequence
$$
T_{(q,s_2')}R \hookrightarrow T_{(q,s_2')}\mathbb{Y}_1 
\twoheadrightarrow \mathcal{N}_{(q,s_2')R|\mathbb{Y}_1}
$$
allows us to write
\begin{equation}\label{normalR1}
\mathcal{N}_{(q,s_2')R|\mathbb{Y}_1} = \dfrac{x_0x_2}{x_1^2} + \dfrac{x_0}{x_2} + \dfrac{x_0^2}{x_1x_2} + \dfrac{x_0}{x_1} + \dfrac{x_0x_3}{x_1x_2}.
\end{equation}
We may conclude from this the existence of the five isolated points associated to these weights, as we have already done in  Table \ref{x13x12s2}.

\begin{exem}\normalfont
For $v_3'=(0:0:0:1:0)$ in Table \ref{x13x12s2} one has
$$c_1^T(\mathcal{W})|_{(q,s_2',v_3')} = w_0 + 2w_1 + w_2 \ ,$$
since $\mathcal{W}|_{(q,s_2',v_3')} = \langle 2x_1^2x_2dx_0-3x_0x_1x_2dx_1 + x_0x_1^2dx_2 \rangle $. \medskip

The tangent space of $\mathbb{Y}_3$ at this fixed point is
$$T_{(q,s_2',v_3')}\mathbb{Y}_3 = T_{(q,s_2')}R \oplus \mathcal{L}_{x_0/x_1} \oplus T_{[\mathcal{L}_{x_0/x_1}]}\mathbb{P}(\mathcal{N}_{(q,s_2')R|\mathbb{Y}_1})$$
$$=\left( \dfrac{x_0}{x_1} + \dfrac{x_0}{x_2} \right) + \left( \dfrac{x_0}{x_1} \right) + \left( \dfrac{x_2}{x_1} + \dfrac{x_1}{x_2} + \dfrac{x_0}{x_2} + \dfrac{x_3}{x_2} \right).$$ \medskip
This yields the equivariant class
$$c_7^T(T_{(q,s_2',v_3')}\mathbb{Y}_3) = (w_0-w_1)^2(w_0-w_2)^2(w_2-w_1)^2(w_2-w_3)\ . $$ $\square$	
\end{exem}

\subsubsection{Resolution on the fixed $\mathbb{P}^1$ over $(x_1^2,x_1^3)$.}

Inside the line $\ell_1 =\{(0:0:0:0:s_4:s_5)\}$ there's
another indeterminacy point \linebreak
$p_{\ell}=(0:0:0:0:2:3)
$: recall the fiber of $\mathcal{W}$ at this line $\ell_1 =\{(0:0:0:0:s_4:s_5)\}$ is 
\begin{equation}\label{fibral1}
\omega=\textcolor{blue}{(2s_5-3s_4)}x_1^3dx_0-\textcolor{blue}{(2s_5-3s_4)}x_0x_1^2dx_1,
\end{equation}
see Table \ref{x13x12}. We may be tempted to divide (\ref{fibral1}) by the commom factor $\blu{2s_5-3s_4}$, but we need to pay attention to some things:\\

$\bullet$ There is a single point in $\ell_1$, the
point $p_{\ell_1}=\{(0:0:0:0:2:3)\},$ in the
indeterminacy locus (note that $\omega$
(\ref{fibral1}) vanishes at this point). In fact,
this point satisfies the equations of our next blowup center $R$ (\ref{JR}). 

$\bullet$ After the blowup along $R$, over the point $p_{\ell_1}$ there is a $\mathbb{P}^4=\mathbb{P}(\mathcal{N}_{(q,p_{\ell_1})R|\bb{Y}_1})$, with five fixed points, for which the contributions to Bott's formula can not be forgotten. The fibers of $\cl W$ over these points are all distinct, see Table \ref{x13x12pl}.\\ 

Then, the contribution of the fixed component $\ell_1$ can not be computed yet. We just obtained in (\ref{NqCY})
$$\mathcal{N}_{qC|\mathbb{Y}} = \dfrac{x_0^3}{x_1^3} + \dfrac{x_0^2}{x_1^2} + \dfrac{x_0^2x_2}{x_1^3} + \dfrac{x_0^2x_3}{x_1^3} + \dfrac{x_0^2}{x_1^2} + \dfrac{x_0x_2}{x_1^2}.$$

Then 
$$T_{(q,p_{\ell_1})}\mathbb{Y}_1 = T_qC \oplus \mathcal{L}_{x_0^2/x_1^2} \oplus T_{[\mathcal{L}_{x_0^2/x_1^2}]}\mathbb{P}(\mathcal{N}_{qC|\mathbb{Y}})=$$
\begin{equation}\label{TPl}
= \dfrac{x_0}{x_1} + \dfrac{x_0^2}{x_1^2} +
\dfrac{x_2}{x_0} + \dfrac{x_0}{x_1}\blu{+1+} \dfrac{x_2}{x_1} + \dfrac{x_3}{x_1}.
\end{equation}

Note the summand \,\blu1\, corresponding to an
eigenvector with weight 0. 

The ruled surface $R$, with base $C$ has fiber over the point $q=(x_1^2,x_1^3)$ 
$$\mathbb{P}^1 = \{(0:0:s_2:0:2s_5:3s_5)\}$$ (compare with
(\ref{P1exc}), and remember that $p_{\ell_1}$ is the point $s_2=0$). It follows from this 
\begin{equation}\label{Tl1R}
T_{(q,p_{\ell_1})}R = T_qC \oplus T_{p_{\ell_1}}\mathbb{P}^1 = \dfrac{x_0}{x_1}+\dfrac{x_2}{x_0}.
\end{equation}
The exceptional $\mathbb{P}^4 = \mathbb{P}(\mathcal{N}_{(q,p_{\ell_1})R|\mathbb{Y}_1})$, where
\begin{equation} \label{Npl1}
 \mathcal{N}_{(q,p_{\ell_1})R|\mathbb{Y}_1} = \dfrac{x_0^2}{x_1^2} +\dfrac{x_0}{x_1} + 1 +\dfrac{x_2}{x_1} +\dfrac{x_3}{x_1}\ .
\end{equation}
\newpage

\begin{table}[h!]
	\centering
	\caption{Fixed points over $(x_1^3,x_1^2), p_{\ell_1}=(0:0:0:0:1:3/2)$}
	\label{x13x12pl}
	\begin{tabular}{ccc}  \hline \hline \\
\textrm{fixed point} & associated & \textrm{generator of fiber} \\
$(v_0:v_1:v_2:v_3:v_4)$ & eigenvector & \textrm{in } $\mathcal{W}$ \\
\hline
& & \\
$(1:0:0:0:0)$ & 1 & $\omega=x_1^3dx_0-x_0x_1^2dx_1$\\
$(0:1:0:0:0)$ & $x_3/x_1$ & $\omega=2x_1^2x_3dx_0-3x_0x_1x_3dx_1 + x_0x_1^2dx_3$\\
$(0:0:1:0:0)$ & $x_0/x_1$ & $\omega=x_0x_1^2dx_0 - x_0^2x_1dx_1$ \\
$(0:0:0:1:0)$ & $x_2/x_1$ & $\omega=2x_1^2x_2dx_0-3x_0x_1x_2dx_1 + x_0x_1^2dx_2$\\
$(0:0:0:0:1)$ & $x_0^2/x_1^2$ & $\omega=x_0^2x_1dx_0-x_0^3dx_1$
\end{tabular}
\end{table}

Notice that we have five fixed points, but the strict
transform $\tilde{\ell_1}$ of $\ell_1$ passes through one of them (the point $v_0':=(1:0:0:0:0)$, associated to the weight 0). Hence, we
gain the contribution of four isolated fixed  points. For
the point $v_0'$, we will compute the contribution of the
strict transform $\tilde{\ell_1}$, since this point lies
inside a higher dimensional component of the fixed point locus. \ref{fibrasx13x12pl}

\begin{exem}\normalfont
	For $v_2'=(0:0:1:0:0)$ we can compute, using (\ref{Npl1}) and (\ref{Tl1R}),
	$$T_{(q,p_{\ell_1},v_2')}\mathbb{Y}_3 = T_{(q,p_{\ell_1})}R \oplus \mathcal{L}_{x_0/x_1} \oplus T_{[\mathcal{L}_{x_0/x_1}]}\mathbb{P}(\mathcal{N}_{(q,p_{\ell_1})R|\mathbb{Y}_1})=$$
	$$= \left(\dfrac{x_0}{x_1} + \dfrac{x_2}{x_0}\right) + \left(\dfrac{x_0}{x_1}\right) + \left(\dfrac{x_0}{x_1}+\dfrac{x_1}{x_0}+ \dfrac{x_2}{x_0} + \dfrac{x_3}{x_0}\right)$$
	and $$c_7^T(T_{(q,p_{\ell},t_2')}\mathbb{Y}_3)=(w_0-w_1)^4(w_2-w_0)^2(w_0-w_3).$$
Moreover, Table \ref{x13x12pl} shows $\mathcal{W}|_{(q,p_{\ell_1},v_2')} = \langle x_0x_1^2dx_0 - x_0^2x_1dx_1 \rangle$, so
$$c_1^T(\mathcal{W})|_{(q,p_{\ell_1},v_2')} = 2w_0+2w_1 \ .$$
	$\square$
	
\end{exem}

The fixed point $v_0'=(1:0:0:0:0)$ lies in the strict
transform $\tilde \ell_1$ of $\ell_1$. The
action is trivial on $\ell_1 = \{(0:0:0:0:s_4:s_5)\}$,
since $w_{s_4} = w_{s_5}$, and  %$T \ell = 1$ (that is,
%the eigenvector is associated to the weight 0) and 
it follows from (\ref{TPl}) the decomposition of the fiber of normal bundle 
\begin{equation}\label{Nl1}
\mathcal{N}_{(q,p_{\ell_1})\ell_1 | \mathbb{Y}_1} = \underbracket{\dfrac{x_0}{x_1}}_{T_qC=\mathcal{O}_{\ell_1}} + \underbracket{\dfrac{x_0^2}{x_1^2}}_{\mathcal{L}_{x_0^2/x_1^2}=\mathcal{O}_{\ell_1}(-1)} + \underbracket{\dfrac{x_2}{x_0} + \dfrac{x_0}{x_1} + \dfrac{x_2}{x_1} + \dfrac{x_3}{x_1}}_{\mathcal{N}_{\mathbb{P}^1 | \mathbb{P}^5}=\mathcal{O}_{\ell_1}(1)^{\oplus 4}}.
\end{equation}

The term $T_qC$ corresponds to the pullback of the
tangent space of $C$, restricted to the
$\ell_1$ that lives in the
$\mathbb{P}^5$-fiber.

To find the decomposition of $\mathcal{N}_{(q,p_{\ell_1},v_0')\tilde \ell_1 |
  \mathbb{Y}_3}$, remember (\ref{Npl1}) and (\ref{Tl1R}) to write
	$$T_{(q,p_{\ell_1},v_0')}\mathbb{Y}_3 = T_{(q,p_{\ell_1})}R \oplus \mathcal{L}_{1} \oplus T_{[\mathcal{L}_{1}]}\mathbb{P}(\mathcal{N}_{(q,p_{\ell_1})R|\mathbb{Y}_1})=$$
\begin{equation}\label{Tx12x13pl1v0}
= \left(\dfrac{x_0}{x_1} + \dfrac{x_2}{x_0}\right) + 1 + \left[\left(\dfrac{x_0^2}{x_1^2}+\dfrac{x_0}{x_1}+ \dfrac{x_2}{x_1} + \dfrac{x_3}{x_1}\right) \cdot 1^{\vee} \right]
\end{equation}
Thus,
\begin{equation}\label{Nl1tilde}
\mathcal{N}_{(q,p_{\ell_1},v_0')\tilde \ell_1 |\mathbb{Y}_3} = \dfrac{x_0}{x_1} + \dfrac{x_2}{x_0} + \dfrac{x_0^2}{x_1^2}+\dfrac{x_0}{x_1}+ \dfrac{x_2}{x_1} + \dfrac{x_3}{x_1}\ .
\end{equation}

The weights in the sum (\ref{Nl1tilde}) are the same as in (\ref{Nl1}), because the point $v_0'$ is associated with weight 0. This will always be repeated on the other fixed lines ahead. 

Despite this, the global description of normal bundle $\mathcal{N}_{\tilde \ell_1 | \mathbb{Y}_3}$ do not need to be the same as in (\ref{Nl1}). There are integers $d_1, \ \ldots \ , d_6$ with
\begin{equation}\label{Nl1tildeds}
\mathcal{N}_{(q,p_{\ell_1},v_0')\tilde \ell_1 | \mathbb{Y}_3} = \underbracket{\dfrac{x_0}{x_1}}_{\mathcal{O}_{\tilde{\ell_1}}(d_1)} + \underbracket{\dfrac{x_0^2}{x_1^2}}_{\mathcal{O}_{\tilde{\ell_1}}(d_2)} + \underbracket{\dfrac{x_2}{x_0}}_{\mathcal{O}_{\tilde \ell_1}(d_3)} + \underbracket{\dfrac{x_0}{x_1}}_{\mathcal{O}_{\tilde \ell_1}(d_4)} + \underbracket{\dfrac{x_2}{x_1}}_{\mathcal{O}_{\tilde \ell_1}(d_5)} + \underbracket{\dfrac{x_3}{x_1}}_{\mathcal{O}_{\tilde \ell_1}(d_6)}.
\end{equation}

\begin{exem}\label{calcl1} \normalfont
	
\ref{l1m2} With (\ref{Nl1tildeds}), we can compute the contribution of $\tilde{\ell_1}$ for Bott. 

The equivariant Chow ring of $\tilde \ell_1 = \mathbb{P}^1$ 
is $$\mathcal{A}_*^T(\tilde \ell_1) = \mathcal{A}_*(\tilde
\ell_1) \otimes \mathbb{Z}[t],$$ where
$\mathcal{A}_*(\tilde \ell_1)=\mathbb{Z}[h]/h^2$. The
decomposition of $\mathcal{N}_{(q,p_{\ell_1},v_0')\tilde \ell_1 |
	\mathbb{Y}_3}$ gives \,$c_6^T(\mathcal{N}_{\tilde \ell_1 | \mathbb{Y}_3}) =$
$$ 
(w_0-w_1+d_1h)(2w_0-2w_1+d_2h)(w_2-w_0+d_3h)(w_0-w_1+d_4h)(w_2-w_1+d_5h)(w_3-w_1+d_6h)
$$
{\em i.e,} $c_6^T(\mathcal{N}_{\tilde \ell_1 | \mathbb{Y}_3}) = \mb{w} + \mb{d}h$,
where
$$\mb{w}=2(w_0-w_1)^3(w_2-w_1)(w_3-w_1)(w_2-w_0)$$
$$\mb{d}=(w_0-w_1)^2[2(w_2-w_0)(w_2-w_1)(w_3-w_1)d_1 + \cdots + 2(w_0-w_1)(w_2-w_1)(w_2-w_0)d_6].$$

\medskip
Along the strict transform $\tilde \ell_1$, the line bundle $\mathcal{W}|_{\tilde \ell_1}$ is trivial (the fibers are generated by $\omega=x_1^3dx_0-x_0x_1^2dx_1$ everywhere). Hence $c_1(\mathcal{W}|_{\tilde \ell_1})=0$, and then
$c_1^T(\mathcal{W}|_{\tilde \ell_1})=w_0+3w_1.$ 

The contribution of the fixed $\tilde\ell_1$ may be computed by
$$\int\limits_{\tilde \ell_1}\dfrac{-c_1^T(\mathcal{W}|_{\tilde \ell_1})^7}{c_6^T(\mathcal{N}_{\tilde \ell_1 | \mathbb{Y}_3})}
=-\int\limits_{\tilde \ell_1}\dfrac{(w_0+3w_1)^7}{\mb{w}+\mb{d}h} =
-\int\limits_{\tilde \ell_1}\dfrac{(w_0+3w_1)^7[\mb{w}-\mb{d}h]}{\mb{w}^2} \ ,$$
where in the last equality we multiplied by $\mb{w}-\mb{d}h$  both
the numerator and the denominator, and used that $h^2=0$ (since $\mathcal{A}_*(\tilde \ell_1)=\mathbb{Z}[h]/h^2$).

Taking the coefficient of $h $ we find
$$\int\limits_{\tilde \ell_1}\dfrac{-c_1^T(\mathcal{W}|_{\tilde \ell_1})^7}{c_6^T(\mathcal{N}_{\tilde \ell_1 | \mathbb{Y}_3})}
=\dfrac{d(w_0+3w_1)^7}{\mb{w}^2} \ .$$
Notice that this contribution of $\tilde \ell_1$ is a linear expression on the unknown integers $d_1, \ldots, d_6$. Further ahead, we will find some useful relations for $d_i's$, see Proposition \ref{relacoesd}.
\def\qed{\,\hfill\ensuremath{_{\mbox{$\square$}}}} \qed
\end{exem}

\subsection {$(g,g',f) = (x_0^2,x_0x_1,x_0^3).$} \label{secaox01}

We will adapt the same ideas of section \ref{primpntfx} to compute the contributions of all fixed components that appear over the point $q:=(g,g',f) = (x_0^2,x_0x_1,x_0^3)$ when solving the indeterminacies as in \S\ref{solvb0u1}.

The fixed point lie in the neighborhood $a_0=b_0=u_1=1$. Take
$$
\left\{\ba l
f=x_0^3 + a_1x_0^2x_1 + a_2x_0^2x_2+a_3x_0^2x_3 + a_4(x_0x_1^2+u_2x_0x_1x_2 + \frac{2}{3}u_3x_1^3)\\\na8
g=x_0^2 + b_1(x_0x_1 +    u_2x_0x_2 +    u_3x_1^2).
\ea\right.
$$

The coefficients inherit the weights
\begin{equation}\label{pesos01}
\begin{array}{rclcrcl} 
w_{\mb{a}_1} & = & w_1-w_0 ,    & & w_{\mb{b}_1}& = & w_1-w_0,\\
w_{\mb{a}_2} & = & w_2-w_0  ,   & & w_{\mb{u}_2}& = & w_2-w_1,\\
w_{\mb{a}_3} & = & w_3-w_0   ,  & & w_{\mb{u}_3}& = & w_1-w_0,\\
w_{\mb{a}_4} & = & 2w_1-2w_0 .& &        &   & 
\end{array}
\end{equation}
Presently the equations of $C$ are (\ref{JC}) 
$$\left\{\begin{array}{rcl}
b_1-4u_3 & & (\textrm{weight } w_0-w_1) \\
a_1-6u_3 & & (\textrm{weight } w_0-w_1) \\
a_2 & & (\textrm{weight } w_0-w_2) \\
a_3 & & (\textrm{weight } w_0-w_3) \\
a_4-12u_3^2 & & (\textrm{weight } 2w_0-2w_1) \\
u_2 & & (\textrm{weight } w_1-w_2)
\end{array}\right.$$

(Remember that, for equations, the weights have opposite sign to the dual points). Notice that these are all weighted homogeneous equations. Hence, over the indeterminacy point
$q=(x_0^2,x_0x_1,x_0^3)$ we have a
$$\mathbb{P}^5=\{(s_0:s_1:s_2:s_3:s_4:s_5)
\}=\mathbb{P}(\mathcal{N}_{qC|\mathbb{Y}})$$
with a natural induced action and weights (\ref{eqbup3})
\begin{equation}\label{wsx01}
\begin{array}{rclcrcl}
w_{s_0} & = & w_1-w_0 & & w_{s_3}& = & w_3-w_0\\
w_{s_1} & = & w_1-w_0 & & w_{s_4}& = & 2w_1-2w_0\\
w_{s_2} & = & w_2-w_0 & & w_{s_5}& = & w_2-w_1.
\end{array}
\end{equation}

Over this $\mathbb{P}^5$ we have four isolated fixed points and a fixed line $\ell_2$.
\newpage

\begin{table}[h!]
	\centering
	\caption{Fixed points over $(g,g',f)=(x_0^2,x_0x_1,x_0^3)$}
	\label{x01}
	\begin{tabular}{ccc}  \hline \hline \\
\textrm{fixed point/ }$\mathbb{P}^1$ & associated & \textrm{generator of fiber} \\
$(s_0:s_1:s_2:s_3:s_4:s_5)$ & eigenvector & \textrm{in } $\mathcal{W}$ \\
\hline
& & \\
$(s_0:s_1:0:0:0:0)$ & $x_1/x_0$ & $\omega=\textcolor{blue}{(3s_0-2s_1)}x_0^2x_1dx_0-\textcolor{blue}{(3s_0-2s_1)}x_0^3dx_1$\\
$(0:0:1:0:0:0)$ & $x_2/x_0$ & $\omega=x_0^2x_2dx_0-x_0^3dx_2$\\
$(0:0:0:1:0:0)$ & $x_3/x_0$ & $\omega=x_0^2x_3dx_0-x_0^3dx_3$\\
$(0:0:0:0:1:0)$ & $x_1^2/x_0^2$ & $\omega=x_0x_1^2dx_0 -  x_0^2x_1dx_1$\\
$(0:0:0:0:0:1)$ & $x_2/x_1$ & \textcolor{red}{\textrm{not defined}}$$
\end{tabular}
\end{table}

See \ref{fibrasx01} for calculations. In this case, there are still two points to solve, 
to wit ,
 $$s_5'=(0:0:0:0:0:1)$$ and the point
 \begin{equation} \label{defl2}
 p_{\ell_2}=(2:3:0:0:0:0) \in \ell_2=\{(s_0:s_1:0:0:0:0)\}.
 \end{equation}

The second step is to blowup along the transform of $E$ (\ref{summarized}).

\subsubsection {Resolution on the point $(x_0^2,x_0x_1,x_0^3,s_5')$}

Now, we take the blowup along the strict transform of $E$ and look at the fixed point $s_5'=(0:0:0:0:0:1) \in \mathbb{Y}_1$.

Remember that to find the equations of the transform, we have an ideal (see (\ref{idealK}))
$$J_{red}=rad(J)=\langle 3b_1^2-4a_4,a_3,a_2,2a_1-3b_1,b_1u_2,b_1(b_1-4u_3) \rangle,$$
which represents the two components, the curve $C$ and the surface
$E$. Take the transform of this ideal $J_{red}$ by the blowup along $C$.

The ideal of the strict transform of $E$ is (\ref{KE'})
$$J_E'=\langle s_3,s_2,3s_0-2s_1,6u_3+s_1u_2,3s_4+s_1^2u_2 \rangle.$$

Over the point $s_5'=$(0:0:0:0:0:1) we have a $\mathbb{P}^4=\{(t_0:t_1:t_2:t_3:t_4) \} = \mathbb{P}(\mathcal{N}_{(q,s_5')E|\mathbb{Y}_1})$.

Let's see how the action passes to this $\mathbb{P}^4$. We are in the neighborhood $s_5=1$, and in this affine chart,
$$(\mb{s}_0,\mb{s}_1,\mb{s}_2,\mb{s}_3,\mb{s}_4)=\left(\dfrac{s_0}{s_5},\dfrac{s_1}{s_5},\dfrac{s_2}{s_5},\dfrac{s_3}{s_5},\dfrac{s_4}{s_5}\right).$$

Thus, from (\ref{wsx01}) we have local weights for affine coordinates around $s_5=1$, given by
\begin{equation}\label{wsx01loc}
\begin{array}{rccclcrcl}
	w_{s_0}^L &=& w_{s_0}-w_{s_5} &=& 2w_1-w_0-w_2  & & w_{s_2}^L &=& w_1-w_0\\
	w_{s_1}^L &=& w_{s_1}-w_{s_5} &=& 2w_1-w_0-w_2       & & w_{s_3}^L &=& w_3+w_1-w_0-w_2\\
	&   &                 &   &               & & w_{s_4}^L & = & 3w_1-2w_0-w_2
\end{array}
\end{equation}

Remember that the equations have opposite sign weigths by duality. The local equations of the blowup center $E$ are (from (\ref{pesos01}) and (\ref{wsx01loc}))
$$\left\{\begin{array}{rcl}
	\mb{s}_3 & & (\textrm{weight \ } w_0+w_2-w_1-w_3) \\
	\mb{s}_2 & & (\textrm{weight \ } w_0-w_1) \\
	3\mb{s}_0-2\mb{s}_1 & & (\textrm{weight \ } w_0+w_2-2w_1) \\
	6u_3+\mb{s}_1u_2 & & (\textrm{weight \ } w_0-w_1) \\
	3\mb{s_4}+\mb{s}_1^2u_2 & & (\textrm{weight \ } 2w_0+w_2-3w_1)
\end{array}\right.
$$

and so relations (\ref{eqbup4}) induce

$$\begin{array}{rclcrcl}
	w_{t_0} & = & w_1+w_3-w_0-w_2     & & w_{t_3}& = & w_1-w_0\\
	w_{t_1} & = & w_1-w_0 & & w_{t_4}& = & 3w_1-2w_0-w_2\\
	w_{t_2} & = & 2w_1-w_0-w_2  & & 
\end{array}$$ 
\begin{equation} \label{acaoind3}
t \circ (t_0:t_1:t_2:t_3:t_4) = (t^{w_{t_0}}t_0 : t^{w_{t_1}}t_1 : t^{w_{t_2}}t_2 : t^{w_{t_3}}t_3 : t^{w_{t_4}}t_4)
\end{equation}
\medskip

An important remark is that the point $q=(x_0^2,x_0x_1,x_0^3)$ lives in the intersection of the two components $C$ and $E$. 

Now we will describe the tangent space of $\mathbb{Y}_2$ at the fixed points.

The first blowup center $C$ is parametrized by $x_0+tx_1$, hence 
\begin{equation}\label{TqCx01}
T_qC = \dfrac{x_1}{x_0} \ .
\end{equation}

The tangent space of $\mathbb{Y}$ at $q$ decomposes  as 
$$T_q\mathbb{Y} = T_{(x_0^2,x_0x_1)}\mathbb{X}' \oplus T_{x_0^3}\mathbb{Y}_{(x_0^2,x_0x_1)},$$
and the fiber $\mathbb{Y}_{(x_0^2,x_0x_1)}$ is
the projectivization of 
$$x_0^3 \oplus x_0^2x_1 \oplus x_0^2x_2 \oplus x_0^2x_3 \oplus x_0x_1^2 \ .$$

From this results $$T_q\mathbb{Y} = \left(\dfrac{x_1}{x_0} + \dfrac{x_2}{x_1} + \dfrac{x_1}{x_0}\right) + \left(\dfrac{x_1}{x_0} + \dfrac{x_2}{x_0} + \dfrac{x_3}{x_0} + \dfrac{x_1^2}{x_0^2}\right)$$ 
and the the normal bundle of $C$ in $\mathbb{Y}$ at $q$ is
\begin{equation} \label{NqCY01}
\mathcal{N}_{qC|\mathbb{Y}}=\dfrac{x_2}{x_1} + \dfrac{x_1}{x_0} + \dfrac{x_1}{x_0} + \dfrac{x_2}{x_0} + \dfrac{x_3}{x_0} + \dfrac{x_1^2}{x_0^2}\end{equation}
(cf. Table \ref{x01}). On the indeterminacy point
$s_5'=(0:0:0:0:0:1)$, related to $\dfrac{x_2}{x_1}$ (\ref{wsx01}), we have
$$T_{(q,s_5')}\mathbb{Y}_1 = T_qC \oplus \mathcal{L}_{\frac{x_2}{x_1}} \oplus T_{[\mathcal{L}_{\frac{x_2}{x_1}}]}\mathbb{P}(\mathcal{N}_{qC|\mathbb{Y}}).$$

This yields
$$T_{(q,s_5')}\mathbb{Y}_1 = \left(\dfrac{x_1}{x_0}\right) + \left(\dfrac{x_2}{x_1}\right) + \left(\dfrac{x_1^2}{x_0x_2} + \dfrac{x_1^2}{x_0x_2} + \dfrac{x_1}{x_0} + \dfrac{x_1x_3}{x_0x_2} + \dfrac{x_1^3}{x_0^2x_2}\right).$$ \medskip

Now, $E'$ is the transform of $E$ under the first blowup. But $C \cap E = \{q\}$, and so $T_{(q,s_5')}E'$ is the tangent space of the blowup of $\mathbb{P}^2$ at one point.

On $E$, spanned by $\langle x_0x_1,x_0x_2,x_1^2 \rangle$ (\ref{geraE}), we have 
$$T_{x_0x_1}E = \dfrac{x_2}{x_1} + \dfrac{x_1}{x_0}$$ 
and $$T_{(q,s_5')}E' = \mathcal{L}_{\frac{x_2}{x_1}} \oplus T_{[\mathcal{L}_{\frac{x_2}{x_1}}]}\mathbb{P}(\mathcal{N}_{\frac{x_2}{x_1}|E}) = \dfrac{x_2}{x_1} + \dfrac{x_1^2}{x_0x_2}.$$

The short exact sequence
$$T_{(q,s_5')}E' \hookrightarrow T_{(q,s_5')}\mathbb{Y}_1 \twoheadrightarrow \mathcal{N}_{({q,s_5'})E'|\mathbb{Y}_1}$$
shows us that (compare with (\ref{acaoind3}))

$$\mathcal{N}_{({q,s_5'})E'|\mathbb{Y}_1} = \underbrace{\dfrac{x_1}{x_0}}_{t_1'} + \underbrace{\dfrac{x_1^2}{x_0x_2}}_{t_2'} +\underbrace{\dfrac{x_1}{x_0}}_{t_3'} +\underbrace{\dfrac{x_1x_3}{x_0x_2}}_{t_0'} +\underbrace{\dfrac{x_1^3}{x_0^2x_2}}_{t_4'}.$$

There are 3 fixed points and a fixed line $\ell_3 := \{(0:t_1:0:t_3:0) \}$. \ref{fibrasx01s5}

\begin{table}[h!]
\centering
\caption{Fixed points over $(q,s_5')=(x_0^2, x_0x_1,x_0^3),(0:0:0:0:0:1)$}
\label{x01s5}
\begin{tabular}{ccc}  \hline \hline \\
 	\textrm{fixed point/ }$\mathbb{P}^1$ & associated & \textrm{generator of fiber} \\
		$(t_0:t_1:t_2:t_3:t_4)$ & eigenvector & \textrm{in } $\mathcal{W}$ \\
		\hline
		& & \\
		$(1:0:0:0:0)$ & $x_1x_3/x_0x_2$ & $\omega=x_0^2x_3dx_0-x_0^3dx_3$\\
		$(0:0:1:0:0)$ & $x_1^2/x_0x_2$ & $\omega=x_0^2x_1dx_0-x_0^3dx_1$\\
		$(0:0:0:0:1)$ & $x_1^3/x_0^2x_2$ & $\omega=x_0x_1^2dx_0-x_0^2x_1dx_1$\\
		$(0:t_1:0:t_3:0)$ & $x_1/x_0$ & 
		$\omega=\textcolor{blue}{(t_3-t_1)}x_0^2x_2dx_0-\textcolor{blue}{(t_3-t_1)}x_0^3dx_2$
	\end{tabular}
\end{table}

\begin{exem}\normalfont
For $t_0'=(1:0:0:0:0)$ we have $\omega=x_0^2x_3dx_0 - x_0^3dx_3$, so $$c_1^T(\mathcal{W})|_{(q,s_5',t_0')} = 3w_0 + w_3.$$

The tangent space at this point, for $c_7^T(T_{(q,s_5',t_0')}\mathbb{Y}_2)$, is
\begin{equation}\label{Tt0'}
T_{(q,s_5',t_0')}\mathbb{Y}_2 = T_{(q,s_5')}E' \oplus \mathcal{L}_{\frac{x_1x_3}{x_0x_2}} \oplus T_{[\mathcal{L}_{\frac{x_1x_3}{x_0x_2}}]}\mathbb{P}(\mathcal{N}_{(q,s_5')E'|\mathbb{Y}_1}) = 
\end{equation}
$$= \left(\dfrac{x_2}{x_1} + \dfrac{x_1^2}{x_0x_2}\right) + \left(\dfrac{x_1x_3}{x_0x_2}\right) + \left(\dfrac{x_2}{x_3} + \dfrac{x_1}{x_3} + \dfrac{x_2}{x_3} + \dfrac{x_1^2}{x_0x_3}\right).$$ $\square$
\end{exem}

There is still one indeterminacy point $p_{\ell_3}=(0:1:0:1:0) \in \ell_3$. The point $q':=(q,s_5',p_{\ell_3})$ lies in the neighborhood $t_1=1$. The ideal of the next blowup center $R$ is (ideal JRn$\underline{\ \ }$7 in \ref{fibrasx01s5})
\begin{equation}\label{JRpl3}
J_R= \langle u_2, t_3-1, t_2, t_0, 3s_1-2t_4 \rangle .
\end{equation}

Let $\mathbb{P}^4=\{(v_0:v_1:v_2:v_3:v_4)\} =
\mathbb{P}(\mathcal{N}_{q'R|\mathbb{Y}_2})$ the fiber of
the exceptional divisor over $q'$. Here, the fixed points
locus consists of  three well solved fixed points, the strict transform $\tilde \ell_3$ of $\ell_3$ and a new indeterminacy point $(q',v_0')$ \ref{fibrasx01s5pl3}.

\begin{table}[h!]
\centering
\caption{Fixed points over $q'=[(x_0^2,x_0x_1,x_0^3),s_5',p_{\ell_3}]$}
\label{x01s5pl3}
\begin{tabular}{ccc}  \hline \hline \\
		\textrm{fixed point} & associated & \textrm{generator of fiber} \\
		$(v_0:v_1:v_2:v_3:v_4)$ & eigenvector & \textrm{in } $\mathcal{W}$ \\
		\hline
		& & \\
		$(1:0:0:0:0)$ & $x_2/x_1$ & \textcolor{red}{\textrm{not defined}} \\
		$(0:1:0:0:0)$ & 1 & $\omega=x_0^2x_2dx_0-x_0^3dx_2$\\
		$(0:0:1:0:0)$ & $x_1/x_2$ & $\omega=x_0^2x_1dx_0-x_0^3dx_1$ \\
		$(0:0:0:1:0)$ & $x_3/x_2$ & $\omega=x_0^2x_3dx_0-x_0^3dx_3$ \\
		$(0:0:0:0:1)$ & $x_1^2/x_0x_2$ & $\omega=x_0x_1^2dx_0-x_0^2x_1dx_1$
	\end{tabular}
\end{table}

\textit{Weights to $\bb{P}^4=\bb{P}(\cl{N}_{q'R|\bb{Y}_2})$}. In the affine chart $t_1=1$,
$$(\mb{t}_0,\mb{t}_2,\mb{t}_3,\mb{t}_4)=\left(\dfrac{t_0}{t_1},\dfrac{t_2}{t_1},\dfrac{t_3}{t_1},\dfrac{t_4}{t_1}\right).$$

Thus, from (\ref{acaoind3}) we have local weights for affine points around $t_3=1$, given by
\begin{equation}\label{wtx01loc}
\begin{array}{rclcrcl}
w_{t_0}^L &=& w_3-w_2  & & w_{t_3}^L &=& 0\\
w_{t_2}^L &=& w_1-w_2  & & w_{t_4}^L &=& 2w_1-w_0-w_2
\end{array}
\end{equation}

The local equations of the blowup center $R$ (\ref{JRpl3}) are, from (\ref{pesos01}), (\ref{wsx01loc}) and (\ref{wtx01loc})
$$\left\{\begin{array}{rcl}
	u_2 & & (\textrm{weight \ } w_1-w_2) \\
	\mb{t}_3-1 & & (\textrm{weight \ } 0) \\
	\mb{t}_2 & & (\textrm{weight \ } w_2-w_1) \\
	\mb{t}_0 & & (\textrm{weight \ } w_2-w_3) \\
	3\mb{s}_1-2\mb{t}_4 & & (\textrm{weight \ } w_0+w_2-2w_1)
\end{array}\right.
$$

and so relations (\ref{eqbup4}) induce
$$\begin{array}{rclcrcl}
	w_{v_0} & = & w_2-w_1  & & w_{v_3}& = & w_3-w_2\\
	w_{v_1} & = & 0        & & w_{v_4}& = & 2w_1-w_0-w_2\\
	w_{v_2} & = & w_1-w_2  & & 
\end{array}$$ 
\begin{equation} \label{acaoind4}
t \circ (v_0:v_1:v_2:v_3:v_4) = (t^{w_{v_0}}v_0 : t^{w_{v_1}}v_1 : t^{w_{v_2}}v_2 : t^{w_{v_3}}v_3 : t^{w_{v_4}}v_4)
\end{equation}
\medskip

Note that the fixed point $v_1'=(0:1:0:0:0)$ is
associated to the weight $w_{v_1}=0$, and we write  1 as in Table \ref{x01s5pl3} for this eigenvalue. It is quite useful for calculations, like in (\ref{Tx12x13pl1v0}).

Let's see the decomposition of $\mathcal{N}_{q'R|\mathbb{Y}_2}$ . As done for $t_0'$ (\ref{Tt0'}), for $p_{\ell_3}$ we have
\begin{equation}\label{Tx01s5t1}
T_{(q,s_5',p_{\ell_3})}\mathbb{Y}_2 = \left(\dfrac{x_2}{x_1} + \dfrac{x_1^2}{x_0x_2}\right) + \left(\dfrac{x_1}{x_0}\right) + \left(\dfrac{x_3}{x_2} + \dfrac{x_1}{x_2} + \dfrac{x_1^2}{x_0x_2} + 1\right).
\end{equation}

Among the equations of $R$ (\ref{JRpl3}) we have exc$_1=u_2$ (the basis $C$) and the $\mathbb{P}^1$ $$\{(t_0:t_1:t_2:t_3:t_4)=(0:1:0:1:3s_1/2)\}.$$

%From the equations of $R$ (\ref{JRpl3}) follows 
%$$R=\{(0:s_1:s_2:0:0:1),(0:1:0:1:3s_1/2)\},$$
%and then 

The point $p_{\ell_3}=(0:1:0:1:0)$ is the point of equation $s_1=0$, and $T_{p_{\ell_3}}\bb{P}^1 = \dfrac{t_3^\vee}{s_1^\vee} = \dfrac{\mb{s}_1}{\mb{t}_3}$.
$$T_{q'}R = T_qC \oplus T_{p_{\ell_3}}\bb{P}^1 = \dfrac{x_1}{x_0} (\ref{TqCx01}) + \dfrac{x_1^2}{x_0x_2} ((\ref{wsx01loc})+(\ref{wtx01loc})),$$
$$\mathcal{N}_{q'R|\mathbb{Y}_2} = \dfrac{x_2}{x_1} + \dfrac{x_3}{x_2} + \dfrac{x_1}{x_2} + \dfrac{x_1^2}{x_0x_2} + 1 \ .$$

\begin{exem}\normalfont \label{exl2}
\ref{l2m2} The fixed point $v_1'$ lies in the strict transform $\tilde \ell_3$ of $\ell_3$. From (\ref{Tx01s5t1}) we have the decomposition of the normal bundle 
$$\mathcal{N}_{\ell_3 | \mathbb{Y}_2} = \underbracket{\dfrac{x_2}{x_1} + \dfrac{x_1^2}{x_0x_2}}_{T_{(q,s_5')}E=\mathcal{O}_{\ell_3}^{\oplus 2}} + \underbracket{\dfrac{x_1}{x_0}}_{\mathcal{L}_{x_1/x_0}=\mathcal{O}_{\ell_3}(-1)} + \underbracket{\dfrac{x_3}{x_2} + \dfrac{x_1}{x_2} +  \dfrac{x_1^2}{x_0x_2}}_{\mathcal{N}_{\mathbb{P}^1 | \mathbb{P}^4}=\mathcal{O}_{\ell_3}(1)^{\oplus 3}}.$$

The decomposition for the normal bundle of $\tilde \ell_3$ is, for some integers $d_7, \ldots, d_{12}$,

$$\mathcal{N}_{\tilde \ell_3 | \mathbb{Y}_3} = \underbracket{\dfrac{x_2}{x_1}}_{\mathcal{O}_{\tilde \ell_3}(d_7)}+\underbracket{\dfrac{x_1^2}{x_0x_2}}_{\mathcal{O}_{\tilde \ell_3}(d_8)} + \underbracket{\dfrac{x_1}{x_0}}_{\mathcal{O}_{\tilde \ell_3}(d_9)} + \underbracket{\dfrac{x_1}{x_2}}_{\mathcal{O}_{\tilde \ell_3}(d_{10})} + \underbracket{\dfrac{x_3}{x_2}}_{\mathcal{O}_{\tilde \ell_3}(d_{11})} + \underbracket{\dfrac{x_1^2}{x_0x_2}}_{\mathcal{O}_{\tilde \ell_3}(d_{12})}.$$

The contribution of the fixed line $\tilde \ell_3$ may be computed as
$$\int\limits_{\tilde \ell_3}\dfrac{-c_1^T(\mathcal{W}|_{\tilde \ell_3})^7}{c_6^T(\mathcal{N}_{\tilde \ell_3 | \mathbb{Y}_3})}
=-\int\limits_{\tilde \ell_3}\dfrac{(3w_0+w_2)^7}{\mb{w}+\mb{d}h} = -\int\limits_{\tilde \ell_3}\dfrac{(3w_0+w_2)^7[\mb{w}-\mb{d}h]}{\mb{w}^2} \ ,$$
where
$$\mb{w}=-(w_2-w_1)^2(2w_1-w_0-w_2)^2(w_1-w_0)(w_3-w_2)$$
$$\mb{d}=(2w_1-w_0-w_2)^2(w_1-w_0)(w_1-w_2)(w_3-w_2)d_7+\cdots+(w_2-w_1)^2(w_0-w_1)(w_3-w_2)(2w_1-w_0-w_2)d_{12}  $$
and we multiplied by $\mb{w}-\mb{d}h$ in both the numerator and the denominator, and used that $h^2=0$.

Taking the coefficient of $h $ we find a linear expression on $d_7,\ldots ,d_{12},$
$$\int\limits_{\tilde \ell _3}\dfrac{-c_1^T(\mathcal{W}|_{\tilde \ell_3})^7}{c_6^T(\mathcal{N}_{\tilde \ell_3 | \mathbb{Y}_3})}
=\dfrac{\mb{d}(3w_0+w_2)^7}{\mb{w}^2} \ .$$
\qed%$\square$
\end{exem}

We still have to solve the indeterminacy over the
point $(q',v_0')$, see Table \ref{x01s5pl3}. We will take a look at the neighborhood $v_0=1$.
% and in the blowup along $R$, with the equation $u_2$ for the exceptional divisor.

%This blowup is done by the changes (mimicking Macaulay2)
%$$\begin{array}{rclcrcl}
%	t_3-1 &=>& v_1u_2 & \ \ \ \ \ & t_0  &=>& v_3u_2 \\
%	t_2   &=>& v_2u_2 & \ \ \ \ \ & 3s_1 &=>& v_4u_2 + 2t_4
%\end{array}$$
%followed by the division by the divisorial component $u_2$. 

The affine variables now are $(v_1,v_2,v_3,v_4,u_2,t_4,s_2)$. The final blowup center is the one dimensional component
$L$ (see page \pageref{summarized}) defined by the
equations $$s_2=v_1=v_2=v_3=v_4=t_4=0.$$ The
component $L$ is parametrized by $u_2$, and $T_{(q',v_0')}L = w_{u_2} = \dfrac{x_2}{x_1}$ (\ref{pesos01}).

Over the point $(q',v_0')$ is defined a $\mathbb{P}^5 =
\mathbb{P}(\mathcal{N}_{(q',v_0')L|\mathbb{Y}_3}) =
\{(z_0:z_1:z_2:z_3:z_4:z_5)\}$. Since we are in the
affine chart $v_0=1$, we have local weights assigned to
coordinates $\underline{v}$ from (\ref{acaoind4}), and
with (\ref{wsx01loc}), (\ref{wtx01loc}) we gain Table
\ref{x01s5pl3v0} (cf.\,\ref{fibrasx01s5pl3v0}).

\begin{table}[h!]
\centering
\caption{Fixed points over $[(x_0^2,x_0x_1,x_0^3),s_5',(0:1:0:1:0),v_0']$ }
\label{x01s5pl3v0}
\begin{tabular}{ccc}  \hline \hline \\
		\textrm{fixed point} & associated & \textrm{generator of fiber} \\
		$(z_0:z_1:z_2:z_3:z_4:z_5)$ & eigenvector & \textrm{in } $\mathcal{W}$ \\
		\hline
		& & \\
		$(1:0:0:0:0:0)$ & $x_1/x_0$ & $\omega=x_0x_2^2dx_0-x_0^2x_2dx_2$ \\
		$(0:1:0:0:0:0)$ & $x_1/x_2$ & $\omega=x_0^2x_2dx_0-x_0^3dx_2$\\
		$(0:0:1:0:0:0)$ & $x_1^2/x_2^2$ & $\omega=x_0^2x_1dx_0-x_0^3dx_1$ \\
		$(0:0:0:1:0:0)$ & $x_1x_3/x_2^2$ & $\omega=x_0^2x_3dx_0-x_0^3dx_3$ \\
		$(0:0:0:0:1:0)$ & $x_1^3/x_0x_2^2$ & $\omega=x_0x_1^2dx_0-x_0^2x_1dx_1$ \\
		$(0:0:0:0:0:1)$ & $x_1^2/x_0x_2$ & $\omega=2x_0x_1x_2dx_0-x_0^2x_2dx_1-x_0^2x_1dx_2$
	\end{tabular}
\end{table}

To get the contributions of these points, write
$$T_{(q',v_0')}\mathbb{Y}_3 = T_{q'}R \oplus \mathcal{L}_{\frac{x_2}{x_1}} \oplus T_{[\mathcal{L}_{\frac{x_2}{x_1}}]}\mathbb{P}(\mathcal{N}_{(q')R|\mathbb{Y}_2}) = $$
$$= \left(\dfrac{x_1^2}{x_0x_2} + \dfrac{x_1}{x_0}\right) + \left(\dfrac{x_2}{x_1}\right) + \left(\dfrac{x_1}{x_2} + \dfrac{x_1^2}{x_2^2} + \dfrac{x_1x_3}{x_2^2} + \dfrac{x_1^3}{x_0x_2^2}\right)\ .$$
The short exact sequence
$$T_{(q',v_0')}L \hookrightarrow T_{(q,v_0')}\mathbb{Y}_3 \twoheadrightarrow \mathcal{N}_{({q',v_0'})L|\mathbb{Y}_3}$$
shows

\begin{equation}
\label{Nqv0LY3}
\mathcal{N}_{({q',v_0'})}L|\mathbb{Y}_3 = \underbrace{\dfrac{x_1}{x_0}}_\text{{$z_0'$}} + 
\underbrace{\dfrac{x_1}{x_2}}_\text{{$z_1'$}} + \underbrace{\dfrac{x_1^2}{x_2^2}}_\text{{$z_2'$}} + \underbrace{\dfrac{x_1x_3}{x_2^2}}_\text{{$z_3'$}} +
\underbrace{\dfrac{x_1^3}{x_0x_2^2}}_\text{{$z_4'$}} + \underbrace{\dfrac{x_1^2}{x_0x_2}}_\text{{$z_5'$}} \ .
\end{equation}

\begin{exem}\normalfont
For $z_4'=(0:0:0:0:1:0)$ we have $\omega=x_0x_1^2dx_0 - x_0^2x_1dx_1$, so $$c_1^T(\mathcal{W})|_{(q',v_0',z_4')} = 2w_0 + 2w_1.$$

The tangent space at this point is
\begin{equation}
T_{(q',v_0',z_4')}\mathbb{Y}_4 = T_{(q',v_0')}L
\oplus \mathcal{L}_{\frac{x_1^3}{x_0x_2^2}} \oplus
T_{[\mathcal{L}_
    {\frac{x_1^3}{x_0x_2^2}}]}\mathbb{P}(\mathcal{N}_{(q',v_0')L|\mathbb{Y}_3}) = 
\end{equation}
$$\ba l= \left(\dfrac{x_2}{x_1}\right) + \left(\dfrac{x_1^3}{x_0x_2^2}\right) + \left(\dfrac{x_1}{x_0} + \dfrac{x_1}{x_2} + \dfrac{x_1^2}{x_2^2} + \dfrac{x_1x_3}{x_2^2} + \dfrac{x_1^2}{x_0x_2}\right) \cdot \underbrace{\dfrac{x_0x_2^2}{x_1^3}}_{w_{z_4}^\vee}\\
= \left(\dfrac{x_2}{x_1}\right) +
\left(\dfrac{x_1^3}{x_0x_2^2}\right) +
\left(\dfrac{x_2^2}{x_1^2} + \dfrac{x_0x_2}{x_1^2} +
  \dfrac{x_0}{x_1} + \dfrac{x_0x_3}{x_1^2} +
  \dfrac{x_2}{x_1}\right)\Rightarrow
\ea$$\medskip
$c_7^T(T_{(q',v_0',z_4')}\mathbb{Y}_4) = 2(w_2-w_1)^3(3w_1-w_0-2w_2)(w_0+w_2-2w_1)(w_0-w_1)(w_0+w_3-2w_1).$\\ 
 $\square$
\end{exem}

\subsubsection{Resolution on the fixed line over $(x_0^2,x_0x_1,x_0^3)$}

After the blowup of $\mathbb{Y}$ along $C$ we obtain, over the fixed point $q=(x_0^2,x_0x_1,x_0^3)$, the \linebreak indeterminacy point $s_5'$ which was fully solved in the previous section. But there is another indeterminacy point at the fixed line $\ell_2$ (\ref{defl2}, see Table \ref{x01}).

Let $\mathbb{P}^4=\{(t_0:t_1:t_2:t_3:t_4)\}=\mathbb{P}(\mathcal{N}_{(q,p_{\ell _2})E'|\mathbb{Y}_2})$ the fiber of the projectivized normal bundle over  $p_{\ell_2}=(1:3/2:0:0:0:0)$. The ideal of $E'$ is, in affine chart $s_0=1,$
$$J_{E'}= \langle s_4-3u_3, s_3, s_2, 2s_1-3, b_1 \rangle .$$
Note that $s_0=1$ represents the choice of
exc$_1=b_1-4u_3$ for the local equation of the first exceptional divisor, see (\ref{eqbup3}). 
The local weights induced by (\ref{wsx01}) are
\begin{equation}\label{wpl2x01loc}
\begin{array}{rccclcrcl}
w_{s_1}^L &=& w_{s_1}-w_{s_0} &=& 0  & & w_{s_3}^L &=& w_3-w_1\\
w_{s_2}^L &=& w_{s_2}-w_{s_0} &=& w_2-w_1       & & w_{s_4}^L &=& w_1-w_0\\
&   &                 &   &               & & w_{s_5}^L & = & w_0+w_2-2w_1
\end{array}
\end{equation}

From equations of $J_{E'}$ and (\ref{wpl2x01loc})+(\ref{pesos01}) we obtain the weights on this $\bb{P}^4$. \ref{fibrasx01pl2} 

\begin{table}[h!]
\centering
\caption{Fixed points over $(x_0^2,x_0x_1,x_0^3), p_{\ell_2}=(1:3/2:0:0:0:0)$}
\label{x01pl2}
\begin{tabular}{ccc}  \hline \hline \\
		\textrm{fixed point} / $\mathbb{P}^1$ & associated & \textrm{generator of fiber} \\
		$(t_0:t_1:t_2:t_3:t_4)$ & eigenvector & \textrm{in } $\mathcal{W}$ \\
		\hline
		& & \\
		$(t_0:0:0:0:t_4)$ & $x_1/x_0$ & $\omega=\textcolor{blue}{(3t_4-8t_0)}x_0x_1^2dx_0-\textcolor{blue}{(3t_4-8t_0)}x_0^2x_1dx_1$\\
		$(0:1:0:0:0)$ & $x_3/x_1$ & $\omega=x_0^2x_3dx_0 - x_0^3dx_3$\\
		$(0:0:1:0:0)$ & $x_2/x_1$ & $\omega=x_0^2x_2dx_0 - x_0^3dx_2$ \\
		$(0:0:0:1:0)$ & 1 & $\omega=x_0^2x_1dx_0 - x_0^3dx_1$
	\end{tabular}
\end{table}

Oh no! There is another fixed line $\ell _4=\{(t_0:0:0:0:t_4)\}$.  We also have two fixed points and the strict transform $\tilde \ell_2 $ of $\ell_2$. 

Remember that (\ref{NqCY01})
$$\mathcal{N}_{qC|\mathbb{Y}}=\dfrac{x_2}{x_1} + \dfrac{x_1}{x_0} + \dfrac{x_1}{x_0} + \dfrac{x_2}{x_0} + \dfrac{x_3}{x_0} + \dfrac{x_1^2}{x_0^2}\ .$$

Since $p_{\ell_2}$ has weight $\dfrac{x_1}{x_0}$
(Table \ref{x01}), 
\be
\ba r
T_{(q,p_{\ell_2})}\mathbb{Y}_1 = T_qC \oplus
\mathcal{L}_{\dfrac{x_1}{x_0}} \oplus
T_{[\mathcal{L}_
    {\dfrac{x_1}{x_0}}]}\mathbb{P}(\mathcal{N}_{qC|\mathbb{Y}})\\\na9
\label{Tpl2Y1}
= \dfrac{x_1}{x_0} + \dfrac{x_1}{x_0} + \dfrac{x_0x_2}{x_1^2} + \dfrac{x_2}{x_1}+1+ \dfrac{x_3}{x_1} + \dfrac{x_1}{x_0}\ .
\ea
\end{equation}

Since $E$ is spanned by $\langle x_0x_1,x_0x_2,x_1^2 \rangle$ (\ref{geraE}) we have 
$$T_{x_0x_1}E = \dfrac{x_2}{x_1} + \dfrac{x_1}{x_0}$$ 
and 
\begin{equation}\label{TE'pl2}
T_{(q,p_{\ell_2})}E' = \mathcal{L}_{\frac{x_1}{x_0}} \oplus T_{[\mathcal{L}_{\frac{x_1}{x_0}}]}\mathbb{P}(\mathcal{N}_{\frac{x_1}{x_0}|E}) = \dfrac{x_1}{x_0} + \dfrac{x_0x_2}{x_1^2}\ .
\end{equation}

Hence the fiber $\mathbb{P}^4 =
\mathbb{P}(\mathcal{N}_{(q,p_{\ell})E'|\mathbb{Y}_1})$ of
the exceptional divisor, where (cf. Table \ref{x01pl2})
\begin{equation}\label{NE'pl2}
\mathcal{N}_{(q,p_{\ell_2})E'|\mathbb{Y}_1} = \dfrac{x_1}{x_0} +\dfrac{x_2}{x_1} + 1 +\dfrac{x_3}{x_1} +\dfrac{x_1}{x_0}\ .
\end{equation}

This will enable us to find the contributions of the fixed points.

\begin{exem}\normalfont \label{exl3}
\ref{l3m2} For the fixed line $\ell_2$, we have from (\ref{Tpl2Y1})
\begin{equation}\label{Nl2Y1}
\mathcal{N}_{\ell_2 | \mathbb{Y}_1} = \underbracket{\dfrac{x_1}{x_0}}_{T_qC=\mathcal{O}_{\ell_2}} + \underbracket{\dfrac{x_1}{x_0}}_{\mathcal{L}_{x_1/x_0}=\mathcal{O}_{\ell_2}(-1)} + \underbracket{\dfrac{x_0x_2}{x_1^2} + \dfrac{x_2}{x_1} + \dfrac{x_3}{x_1} + \dfrac{x_1}{x_0}}_{\mathcal{N}_{\mathbb{P}^1 | \mathbb{P}^5}=\mathcal{O}_{\ell_2}(1)^{\oplus 4}}\ .
\end{equation}

The normal bundle of $\tilde \ell_2$ is, for some integers $d_{13}, \ldots, d_{18}$,

$$\mathcal{N}_{\tilde \ell_2 | \mathbb{Y}_2} = \underbracket{\dfrac{x_1}{x_0}}_{\mathcal{O}_{\tilde\ell_2}(d_{13})} + \underbracket{\dfrac{x_1}{x_0}}_{\mathcal{O}_{\tilde \ell_2}(d_{14})} + \underbracket{\dfrac{x_0x_2}{x_1^2}}_{\mathcal{O}_{\tilde \ell_2}(d_{15})} + \underbracket{\dfrac{x_2}{x_1}}_{\mathcal{O}_{\tilde\ell_2}(d_{16})} + \underbracket{\dfrac{x_3}{x_1}}_{\mathcal{O}_{\tilde\ell_2}(d_{17})} + \underbracket{\dfrac{x_1}{x_0}}_{\mathcal{O}_{\tilde\ell_2}(d_{18})}$$
so
$$c_6^T(\mathcal{N}_{\tilde \ell_2 | \mathbb{Y}_2}) = [\prod_{i=13,14,18}(w_1-w_0+d_{i}h)](w_0+w_2-2w_1+d_{15}h)(w_2-w_1+d_{16}h)(w_3-w_1+d_{17}h)$$
(mod\,$h^2$).
From Table \ref{x01pl2}, $\cl{W}|_{\tilde p \in \tilde \ell_2} = \langle x_0^2x_1dx_0 - x_0^3dx_1 \rangle$ for all $\tilde p \in \tilde \ell_2$ (the restriction $\cl{W}|_{\tilde \ell_2}$ is trivial), and so
$c_1^T(\cl{W}|_{\tilde \ell_2}) = 3w_0+w_1$. Like Examples \ref{calcl1} and \ref{exl2}, we gain a linear expression on $d_{13},\ldots ,d_{18}$

$$\int\limits_{\tilde \ell_2}\dfrac{-c_1^T(\mathcal{W}|_{\tilde \ell_2})^7}{c_6^T(\mathcal{N}_{\tilde \ell_2 | \mathbb{Y}_2})}
=\dfrac{\mb{d}(3w_0+w_1)^7}{\mb{w}^2} \ .$$
\qed%$\square$
\end{exem}

In view of Table \ref{x01pl2}, the fiber of $\mathcal{W}$ over a point $(t_0:t_4) \in \ell_4$ is spanned by 
$$\omega=\textcolor{blue}{(3t_4-8t_0)}x_0x_1^2dx_0-\textcolor{blue}{(3t_4-8t_0)}x_0^2x_1dx_1.
$$
Thus, the fixed line $\ell_4$ meets the next blowup center $R$ at the fixed point $$p_{\ell_4}=(3:0:0:0:8).$$ 

The ideal of $R$ is locally given by
\begin{equation}\label{JRx01}
J_R=\langle 3t_4-8, t_3, t_1, t_2-4s_5, 2s_4-9u_3 \rangle .
\end{equation}

Let $\mathbb{P}^4=\{(v_0:v_1:v_2:v_3:v_4)\}=\mathbb{P}(\mathcal{N}_{(q,p_{\ell_2},p_{\ell_4})R|\mathbb{Y}_2})$ over $p_{\ell_4}=(1:0:0:0:8/3)$. After taking the convenient weights we are able to find the Table \ref{x01pl2pl4} below (cf.\,\ref{fibrasx01pl2pl4}).

\begin{table}[h!]
\centering
\caption{Fixed points over $(x_0^2,x_0x_1,x_0^3), p_{\ell_2}, (1:0:0:0:8/3)$}
\label{x01pl2pl4}
\begin{tabular}{ccc}  \hline \hline \\
		\textrm{fixed point} & associated & \textrm{generator of fiber} \\
		$(v_0:v_1:v_2:v_3:v_4)$ & eigenvector & \textrm{in } $\mathcal{W}$ \\
		\hline
		& & \\
		$(1:0:0:0:0)$ & 1 & $\omega=x_0x_1^2dx_0-x_0^2x_1dx_1$\\
		$(0:1:0:0:0)$ & $x_0/x_1$ & $\omega=x_0^2x_1dx_0 - x_0^3dx_1$\\
		$(0:0:1:0:0)$ & $x_0x_3/x_1^2$ & $\omega=x_0^2x_3dx_0 - x_0^3dx_3$ \\
		$(0:0:0:1:0)$ & $x_0x_2/x_1^2$ & $\omega=x_0^2x_2dx_0 - x_0^3dx_2$ \\
		$(0:0:0:0:1)$ & $x_1/x_0$ & $\omega=x_1^3dx_0 - x_0x_1^2dx_1$
	\end{tabular}
\end{table}

For short, write here $q'' = [(x_0^2,x_0x_1,x_0^3),p_{\ell_2},p_{\ell_4}]$. From (\ref{TE'pl2})+(\ref{NE'pl2}),
\begin{equation}\label{Tq''Y2}
\ba r
T_{q''}\mathbb{Y}_2 = T_{(q,p_{\ell_2})}E' \oplus \mathcal{L}_{x_1/x_0} \oplus T_{[\mathcal{L}_{x_1/x_0}]}\mathbb{P}(\mathcal{N}_{(q,p_{\ell_2})E'|\mathbb{Y}_1})\\\na9
= \left(\dfrac{x_1}{x_0} + \dfrac{x_0x_2}{x_1^2}\right) +
\left(\dfrac{x_1}{x_0}\right) +
\left(\dfrac{x_0x_2}{x_1^2} + \dfrac{x_0}{x_1} +
  \dfrac{x_0x_3}{x_1^2} + 1\right).
\ea
\end{equation}

The ruled surface $R$ has equations $2s_4-9u_3$, representing the (transform of) basis $C$ and the other four linear equations of (\ref{JRx01}) represent the $\mathbb{P}^1$
$$\{(t_0:t_1:t_2:t_3:t_4)=(1:0:4s_5:0:8/3)\}.$$  
The point $p_{\ell_4}=(1:0:0:0:8/3)$ is the point of equation $s_5=0$, thus (\ref{TqCx01})+(\ref{wpl2x01loc})
$$T_{q''}R = T_qC \oplus T_{p_{\ell_4}}\bb{P}^1 = \dfrac{x_1}{x_0} + \dfrac{t_0^\vee}{{\mb{s}_5^L}^\vee} = \dfrac{x_1}{x_0} + \dfrac{\mb{s}_5^L}{1} = \dfrac{x_1}{x_0} + \dfrac{x_0x_2}{x_1^2}\ .$$

%From the equations of $R$ (\ref{JRx01}) follows  
%$$R=\{(1:3/2:0:0:(9/2)u_3:s_5),(1:0:4s_5:0:8/3)\},$$
%and then
%$$T_{q''}R = w_{u_3^{\vee}} + \dfrac{s_5}{s_0} = \dfrac{x_1}{x_0} + %\dfrac{x_0x_2}{x_1^2},$$
From (\ref{Tq''Y2}) remains 
$$\mathcal{N}_{q''R|\mathbb{Y}_2} = \dfrac{x_1}{x_0} +\dfrac{x_0x_2}{x_1^2} +\dfrac{x_0}{x_1} +\dfrac{x_0x_3}{x_1^2}+1.$$

\begin{exem}\normalfont
For $v_1'=(0:1:0:0:0)$, associated to $x_0/x_1$ (Table \ref{x01pl2pl4}),
$$\ba c
T_{(q'',v_1')}\mathbb{Y}_3 = T_{q''}R \oplus \mathcal{L}_{\frac{x_0}{x_1}} \oplus T_{[\mathcal{L}_{\frac{x_0}{x_1}}]}\mathbb{P}(\mathcal{N}_{q''R|\mathbb{Y}_2}) = 
\dfrac{x_1}{x_0} + \dfrac{x_0x_2}{x_1^2} + \dfrac{x_0}{x_1} + \dfrac{x_1^2}{x_0^2} + \dfrac{x_2}{x_1} +  \dfrac{x_3}{x_1} + \dfrac{x_1}{x_0}\Rightarrow\\\na8
c_7^T(T_{(q'',v_1')}\bb{Y}_3) =
-2(w_0-w_1)^4(w_0+w_2-2w_1)(w_2-w_1)(w_3-w_1).
\ea$$

And $\cl{W}|_{(q'',v_1')} = \langle x_0^2x_1dx_0 - x_0^3dx_1 \rangle$, so $c_1^T(\cl{W}|_{(q'',v_1')}) = 3w_0+w_1.$\\ \qed
\end{exem}

\begin{exem}\normalfont \label{exl4}
\ref{l4m2} For the fixed line $\ell_4$, we have
$$\mathcal{N}_{\ell_4 | \mathbb{Y}_2} = \underbracket{\dfrac{x_1}{x_0}+\dfrac{x_0x_2}{x_1^2}}_{T_{(q,p_{\ell_2})}E=\mathcal{O}_{\ell_4}^{\oplus 2}} + \underbracket{\dfrac{x_1}{x_0}}_{\mathcal{L}_{x_1/x_0}=\mathcal{O}_{\ell_4}(-1)} + \underbracket{\dfrac{x_0x_2}{x_1^2} + \dfrac{x_0x_3}{x_1^2} + \dfrac{x_0}{x_1}}_{\mathcal{N}_{\mathbb{P}^1 | \mathbb{P}^4}=\mathcal{O}_{\ell_4}(1)^{\oplus 3}}\ ,$$
and for unknown integers $d_{19}, \cdots, d_{24}$,

$$\mathcal{N}_{\tilde \ell_4 | \mathbb{Y}_3} = \underbracket{\dfrac{x_1}{x_0}}_{\mathcal{O}_{\tilde \ell_4}(d_{19})} + \underbracket{\dfrac{x_0x_2}{x_1^2}}_{\mathcal{O}_{\tilde \ell_4}(d_{20})} + \underbracket{\dfrac{x_1}{x_0}}_{\mathcal{O}_{\tilde \ell_4}(d_{21})} + \underbracket{\dfrac{x_0x_2}{x_1^2}}_{\mathcal{O}_{\tilde \ell_4}(d_{22})} + \underbracket{\dfrac{x_0x_3}{x_1^2}}_{\mathcal{O}_{\tilde \ell_4}(d_{23})} + \underbracket{\dfrac{x_0}{x_1}}_{\mathcal{O}_{\tilde \ell_4}(d_{24})}.$$

$$\int\limits_{\tilde \ell_4}\dfrac{-c_1^T(\mathcal{W}|_{\tilde \ell_4})^7}{c_6^T(\mathcal{N}_{\tilde \ell_4 | \mathbb{Y}_3})}
=\dfrac{\mb{d}\cdot (2w_0+2w_1)^7}{\mb{w}^2} \ ,
$$
where
$$\mb{w}=-(w_1-w_0)^3(w_0+w_2-2w_1)^2(w_0+w_3-2w_1),
$$ 
$$\mb{d}=[\prod_{i=19,21}(w_1-w_0+d_ih)][\prod_{j=20,22}(w_0+w_2-2w_1+d_jh)](w_0+w_3-2w_1+d_{23}h)(w_0-w_1+d_{24}h)$$
(mod\,$h^2$). \hfill
$\square$
\end{exem}

Finally we complete the contributions over the fixed point $(x_0^2,x_0x_1,x_0^3)\ !$

\subsection {$(g,g',f)=(x_0^2,x_0x_2,x_0^3).$}

Next, we study the indeterminacy point $q:=(x_0^2,x_0x_2,x_0^3).$

It lies in the neighborhood $a_0=b_0=u_2=1$, so take
$$\left\{
\ba l
f=x_0^3 + a_1x_0^2x_1 + a_2x_0^2x_2+a_3x_0^2x_3 + a_5(u_1x_0x_1^2+x_0x_1x_2 + \frac{2}{3}u_3x_1^3),
\\
g=x_0^2 + b_2u_1x_0x_1 + b_2x_0x_2 + b_2u_3x_1^2.
\ea\right.$$

The local coefficients are equipped with the weights
\begin{equation}\label{pesosx02}
\begin{array}{rclcrcl}
w_{\mb{a}_1} & = & w_1-w_0      & & w_{\mb{b}_2}& = & w_2-w_0\\
w_{\mb{a}_2} & = & w_2-w_0      & & w_{\mb{u}_1}& = & w_1-w_2\\
w_{\mb{a}_3} & = & w_3-w_0      & & w_{\mb{u}_3}& = & 2w_1-w_0-w_2\\
w_{\mb{a}_5} & = & w_1+w_2-2w_0\ . & &        &   & 
\end{array}
\end{equation}

Notice that in this neighborhood there are no points in the first blowup center $C$.

The local equations of $E$ are 
$$\left\{\begin{array}{rcl}
	a_1 & & (\textrm{degree } w_0-w_1) \\
	a_2 & & (\textrm{degree } w_0-w_2) \\
	a_3 & & (\textrm{degree } w_0-w_3) \\
	a_5 & & (\textrm{degree } 2w_0-w_1-w_2) \\
	b_2 & & (\textrm{degree } w_0-w_2)
\end{array}\right.
$$

Notice that these are all homogeneous equations. Hence, over the indeterminacy point $(x_0^2,x_0x_2,x_0^3)$ we have a $\mathbb{P}^4=\{(t_0:t_1:t_2:t_3:t_4) \} = \mathbb{P}(\mathcal{N}_{qE|\mathbb{Y}_1})$ with a natural induced action and weights
\begin{equation}\label{acaoind6}
\begin{array}{rclcrcl}
w_{t_0} & = & w_1-w_0 & & w_{t_3}& = & -2w_0+w_1+w_2\\
w_{t_1} & = & w_2-w_0 & & w_{t_4}& = & w_2-w_0\\
w_{t_2} & = & w_3-w_0 & &        &   & 
\end{array}
\end{equation}

The methods to compute the contributions of $c_1^T(\mathcal{W})$ and the tangent space at the fixed points are similar to those explained in the previous section \ref{secaox01}. Now, we have 
\begin{equation}\label{TEx02}
T_qE'=T_qE = \dfrac{x_1}{x_2} + \dfrac{x_1^2}{x_0x_2}
\end{equation} (without any changes by the blowup in $C$).

Over this $\mathbb{P}^4$ there are three fixed points and a fixed  $\ell_5=\{(0:t_1:0:0:t_4)\}$. \ref{fibrasx02}

\begin{table}[h!]
\centering
\caption{Fixed points over $(g,g',f)=(x_0^2, x_0x_2,x_0^3)$}
\label{x02}
\begin{tabular}{ccc}  \hline \hline \\
		\textrm{fixed point/ }$\mathbb{P}^1$ & associated & \textrm{generator of fiber} \\
		$(t_0:t_1:t_2:t_3:t_4)$ & eigenvector & \textrm{in } $\mathcal{W}$ \\
		\hline
		& & \\
		$(1:0:0:0:0)$ & $x_1/x_0$ & $\omega=x_0^2x_1dx_0-x_0^3dx_1$\\
		$(0:0:1:0:0)$ & $x_3/x_0$ & $\omega=x_0^2x_3dx_0-x_0^3dx_3$\\
		$(0:0:0:1:0)$ & $x_1x_2/x_0^2$ & $\omega=2x_0x_1x_2dx_0 - x_0^2x_2dx_1 - x_0^2x_1dx_2$\\
		$(0:t_1:0:0:t_4)$ & $x_2/x_0$ & $\omega=\textcolor{blue}{(3t_4-2t_1)}x_0^2x_2dx_0-\textcolor{blue}{(3t_4-2t_1)}x_0^3dx_2$
	\end{tabular}
\end{table}

There's one indeterminacy point $p_{\ell_5}:=(0:3:0:0:2)$ inside the fixed line $\ell_5$. 
The remaining indeterminacy locus is 
\begin{equation}\label{eqLx02}
L : u_3=t_3=t_2=2t_1-3=2t_0-3u_1=b_2=0.
\end{equation}

The affine variables are $b_2,u_1,u_3,t_0,t_1,t_2,t_3$, and $L$ is parametrized by $u_1$, so 
\begin{equation}\label{TLx02}
T_{(q,p_{\ell_5})}L=w_{u_1}=\dfrac{x_1}{x_2}\ .
\end{equation}

The blowup of $\bb{Y}_2$ along $L$ gives over $p_{\ell_5}=(0:\frac{3}{2}:0:0:1)$, as fiber of exceptional divisor, $$\bb{P}^5 = \{(z_0:z_1:z_2:z_3:z_4:z_5)\}   = \bb{P}(\cl{N}_{(q,p_{\ell_5})L|\bb{Y}_2}).$$

\textit{Weights to $\bb{P}^5=\bb{P}(\cl{N}_{(q,p_{\ell_5})L|\bb{Y}_2})$}. In the affine chart $t_4=1$,
$$(\mb{t}_0,\mb{t}_1,\mb{t}_2,\mb{t}_3)=\left(\dfrac{t_0}{t_4},\dfrac{t_1}{t_4},\dfrac{t_2}{t_4},\dfrac{t_3}{t_4}\right).$$

Thus, from (\ref{acaoind6}) we have local weights for affine points around $t_4=1$, given by
\begin{equation}\label{wtx02loc}
\begin{array}{rclcrcl}
w_{t_0}^L &=& w_1-w_2  & & w_{t_2}^L &=& w_3-w_2\\
w_{t_1}^L &=& 0        & & w_{t_3}^L &=& w_1-w_0
\end{array}
\end{equation}

The local equations of the blowup center $L$ (\ref{eqLx02}) are, from (\ref{pesosx02}) and (\ref{wtx02loc})
$$\left\{\begin{array}{rcl}
	u_3 & & (\textrm{weight \ } w_0-2w_1+w_2) \\
	\mb{t}_3 & & (\textrm{weight \ } w_0-w_1) \\
	\mb{t}_2 & & (\textrm{weight \ } w_2-w_3) \\
	2\mb{t}_1-3 & & (\textrm{weight \ } 0) \\
	2\mb{t}_0-3u_1 & & (\textrm{weight \ } w_2-w_1)\\
	b_2 & & (\textrm{weight \ } w_0-w_2)
\end{array}\right.
$$
thus, 
$$\begin{array}{rclcrcl}
	w_{z_0} & = & 2w_1-w_0-w_2  & & w_{z_3}& = & 0\\
	w_{z_1} & = & w_1-w_0       & & w_{z_4}& = & w_1-w_2\\
	w_{z_2} & = & w_3-w_2       & & w_{z_5}& = & w_2-w_0
\end{array}$$ 
\begin{equation} \label{acaoind7}
t \circ (z_0:z_1:z_2:z_3:z_4:z_5) = (t^{w_{z_0}}z_0 : t^{w_{z_1}}z_1 : t^{w_{z_2}}z_2 : t^{w_{z_3}}z_3 : t^{w_{z_4}}z_4: t^{w_{z_5}}z_5)
\end{equation}
From (\ref{acaoind7}) we have the following table \ref{fibrasx02pl5}

\begin{table}[h!]
\centering
\caption{Fixed points over $(x_0^2,x_0x_2,x_0^3),p_{\ell_5}$}
\label{x02pl5}
\begin{tabular}{ccc} \hline \hline \\
		\textrm{fixed point } & associated & \textrm{generator of fiber} \\
		$(z_0:z_1:z_2:z_3:z_4:z_5)$ & eigenvector & \textrm{in } $\mathcal{W}$ \\
		\hline
		& & \\
		$(1:0:0:0:0:0)$ & $x_1^2/x_0x_2$ & $\omega=x_0x_1^2dx_0-x_0^2x_1dx_1$\\
		$(0:1:0:0:0:0)$ & $x_1/x_0$ & $\omega=2x_0x_1x_2dx_0-x_0^2x_2dx_1-x_0^2x_1dx_2$\\
		$(0:0:1:0:0:0)$ & $x_3/x_2$ & $\omega=x_0^2x_3dx_0-x_0^3dx_3$\\
		$(0:0:0:1:0:0)$ & 1 & $\omega=x_0^2x_2dx_0-x_0^3dx_2$\\
		$(0:0:0:0:1:0)$ & $x_1/x_2$ & $\omega=x_0^2x_1dx_0-x_0^3dx_1$ \\
		$(0:0:0:0:0:1)$ & $x_2/x_0$ & $\omega=x_0x_2^2dx_0-x_0^2x_2dx_2$
\end{tabular}
\end{table}

\begin{exem}\normalfont \label{exl5}
\ref{l5m2} On the point $q=(x_0^2,x_0x_2,x_0^3)$, we have (cf. (\ref{TX'02})) $$T_q\bb{Y}=\dfrac{x_2}{x_0}+\dfrac{x_1}{x_2}+\dfrac{x_1^2}{x_0x_2}+\dfrac{x_1}{x_0}+\dfrac{x_2}{x_0}+\dfrac{x_3}{x_0}+\dfrac{x_1x_2}{x_0^2}\ .$$
From (\ref{TEx02}), and the summand $x_2/x_0$ for $p_{\ell_5}$ (Table \ref{x02}), one can find
$$\mathcal{N}_{\ell_5 | \mathbb{Y}_2} = \underbracket{\dfrac{x_1}{x_2}+\dfrac{x_1^2}{x_0x_2}}_{T_qE=\mathcal{O}_{\ell_5}^{\oplus 2}} + \underbracket{\dfrac{x_2}{x_0}}_{\mathcal{L}_{x_2/x_0}=\mathcal{O}_{\ell_5}(-1)} + \underbracket{\dfrac{x_1}{x_2} + \dfrac{x_3}{x_2} + \dfrac{x_1}{x_0}}_{\mathcal{N}_{\mathbb{P}^1 | \mathbb{P}^4}=\mathcal{O}_{\ell_5}(1)^{\oplus 3}}\ .$$
For some integers $d_{25}, \ldots, d_{30}$ we have
$$\mathcal{N}_{\tilde \ell_5 | \mathbb{Y}_4} = \underbracket{\dfrac{x_1}{x_2}}_{\mathcal{O}_{\tilde \ell_5}(d_{25})} + \underbracket{\dfrac{x_1^2}{x_0x_2}}_{\mathcal{O}_{\tilde \ell_5}(d_{26})} + \underbracket{\dfrac{x_2}{x_0}}_{\mathcal{O}_{\tilde \ell_5}(d_{27})} + \underbracket{\dfrac{x_1}{x_2}}_{\mathcal{O}_{\tilde \ell_5}(d_{28})} + \underbracket{\dfrac{x_3}{x_2}}_{\mathcal{O}_{\tilde \ell_5}(d_{29})} + \underbracket{\dfrac{x_1}{x_0}}_{\mathcal{O}_{\tilde \ell_5}(d_{30})}\ .$$

Since $c_1^T(\cl{W}|_{\tilde{\ell_5}})=3w_0+w_2$ (cf. Table \ref{x02pl5}), like examples (\ref{calcl1}), (\ref{exl2}), (\ref{exl3}), (\ref{exl4}) we obtain a linear expression in $d_{25}, \ldots, d_{30}$,
$$\int\limits_{\tilde \ell_5}\dfrac{-c_1^T(\mathcal{W}|_{\tilde \ell_5})^7}{c_6^T(\mathcal{N}_{\tilde \ell_5 | \mathbb{Y}_4})}
=\dfrac{\mb{d}(3w_0+w_2)^7}{\mb{w}^2} \ .$$
\qed
\end{exem}

\subsection {$(g,g',f)=(x_0^2,x_1^2,x_0^3).$}

At last, let's see how to solve indeterminacies around the point $q:=(x_0^2,x_1^2,x_0^3).$

It lies in the neighborhood $a_0=b_0=u_3=1$, and then
$$\left\{\ba l
f=x_0^3 + a_1x_0^2x_1 + a_2x_0^2x_2+a_3x_0^2x_3 + a_6(\frac{3}{2}u_1x_0x_1^2+\frac{3}{2}u_2x_0x_1x_2 + x_1^3)\\
g=x_0^2 + b_3u_1x_0x_1 + b_3u_2x_0x_2 + b_3x_1^2.
\ea\right.$$

The weights of the local coefficients are 
\begin{equation}\label{pesosx11}
\begin{array}{rclcrcl}
w_{\mb{a}_1} & = & w_1-w_0   & & w_{\mb{b}_3}& = & 2w_1-2w_0\\
w_{\mb{a}_2} & = & w_2-w_0   & & w_{\mb{u}_1}& = & w_0-w_1\\
w_{\mb{a}_3} & = & w_3-w_0   & & w_{\mb{u}_2}& = & w_0-2w_1+w_2\\
w_{\mb{a}_6} & = & 3w_1-3w_0 & &        &   & 
\end{array}
\end{equation}

The ideal of the curve $C$ is  
$$J_C=\langle
u_2,a_3,a_2,2a_1-3b_3u_1,2a_6-b_3^2u_1,b_3u_1^2-4 \rangle.$$

Note the presence of the equation $b_3u_1^2=4$. Since on the fixed point  $q=(x_0^2,x_1^2,x_0^3)$ we have $b_3=u_1=0$, it means that the fixed point is outside of this first blowup center.

Next, we look at the second blowup center, the component $E$ with equations 

$$\left\{\begin{array}{rcl}
	b_3 & & (\textrm{degree } 2w_0-2w_1) \\
	a_1 & & (\textrm{degree } w_0-w_1) \\
	a_2 & & (\textrm{degree } w_0-w_2) \\
	a_3 & & (\textrm{degree } w_0-w_3) \\
	a_6 & & (\textrm{degree } 3w_0-3w_1).
\end{array}\right.
$$

We have a $\mathbb{P}^4=\{(t_0:t_1:t_2:t_3:t_4) \}=\mathbb{P}(\mathcal{N}_{qE|\mathbb{Y}})$ over the point $q=(x_0^2,x_1^2,x_0^3)$, with induced action and weights
\begin{equation}\label{wtx11}
\begin{array}{rclcrcl}
w_{t_0} & = & 2w_1-2w_0 & & w_{t_3}& = & w_3-w_0\\
w_{t_1} & = & w_1-w_0   & & w_{t_4}& = & 3w_1-3w_0\\
w_{t_2} & = & w_2-w_0   & & &  & 
\end{array}
\end{equation}

Here, there are five isolated fixed points. The table over this $\mathbb{P}^4$ is \ref{fibrasx12}

\begin{table}[h!]
\centering
\caption{Fixed points over $(x_0^2,x_1^2,x_0^3)$}
\label{x12}
\begin{tabular}{ccc}  \hline \hline \\
		\textrm{fixed point } & associated & \textrm{generator of fiber} \\
		$(t_0:t_1:t_2:t_3:t_4)$ & eigenvector & \textrm{in } $\mathcal{W}$ \\
		\hline
		& & \\
		$(1:0:0:0:0)$ & $x_1^2/x_0^2$ & $\omega=x_0x_1^2dx_0-x_0^2x_1dx_1$\\
		$(0:1:0:0:0)$ & $x_1/x_0$ & $\omega=x_0^2x_1dx_0-x_0^3dx_1$\\
		$(0:0:1:0:0)$ & $x_2/x_0$ & $\omega=x_0^2x_2dx_0-x_0^3dx_2$\\
		$(0:0:0:1:0)$ & $x_3/x_0$ & $\omega=x_0^2x_3dx_0-x_0^3dx_3$\\
		$(0:0:0:0:1)$ & $x_1^3/x_0^3$ & $\omega=x_1^3dx_0-x_0x_1^2dx_1$
	\end{tabular}
\end{table}

This completes the description of the fixed points.

\section{Sum of the contributions}

The methods to solve the indeterminacies of the map $(g,f) \mapsto \omega$ show  us that there are over a fixed flag: \\

$\cdot$ 26 fixed points at $\mathbb{Y}$ (see Table \ref{Fixosem}, page \pageref{Fixosem}) ; 

$\cdot$ 6 fixed points at $\mathbb{Y}_1$ (\textcolor{blue}{3} on Table \ref{x13x12} and \textcolor{blue}{3} on Table \ref{x01});

$\cdot$ 13 fixed points at $\mathbb{Y}_2$ (\textcolor{blue}{3} on Table \ref{x01s5}, \textcolor{blue}{2} on Table \ref{x01pl2}, \textcolor{blue}{5} on Table \ref{x12} and \textcolor{blue}{3} on Table \ref{x02});

$\cdot$ 16 fixed points at $\mathbb{Y}_3$ (\textcolor{blue}{5} on Table \ref{x13x12s2}, \textcolor{blue}{4} on Table \ref{x13x12pl}, \textcolor{blue}{3} on Table \ref{x01s5pl3} and \textcolor{blue}{4} on Table \ref{x01pl2pl4});

$\cdot$ 11 fixed points at $\mathbb{Y}_4$ (\textcolor{blue}{6} on Table \ref{x01s5pl3v0} and \textcolor{blue}{5} on Table \ref{x02pl5});

$\cdot$ 5 fixed $\mathbb{P}^1$ (\textcolor{blue}{1} on Tables \ref{x13x12pl}, \ref{x01s5pl3}, \ref{x01pl2}, \ref{x01pl2pl4} and \ref{x02pl5}). \\

Running over all flags, we obtain a total of $72 \cdot 24
= 1728$ fixed points and $5 \cdot 24 = 120$ fixed lines.

There are so many points! We can do the sum of all
contributions using 
the freeware  \textit{Macaulay2}. The scripts can be found in the appendix \ref{degfib}, \ref{degexc}.

By examples \ref{calcl1}, \ref{exl2}, \ref{exl3}, \ref{exl4} and \ref{exl5} this sum is a linear expression in unknowns $d_1, \ldots, d_{30}$.

\begin{lema} The integers $d_1, \ldots, d_{30}$ satisfy the relations
\begin{equation}\label{relacoesd}
\begin{array}{c|c|c|c}
2d_1+d_2+2d_4+2=0 & d_3=1 & d_5=0 & d_6=0 \\
2d_7-2d_{10}-d_{26}-2d_{30}+1=0 & d_8+d_{12}=0 & d_9+d_{30}+1=0 & d_{11}=0 \\
d_{13}+d_{14}+d_{18}+2=0 & d_{15}=1 & d_{16}=0 & d_{17}=0 \\
d_{19}+d_{21}-d_{24}+2=0 & d_{20}+d_{22}=0 & d_{23}=0 & \\
2d_{25}+d_{26}+2d_{28}+2d_{30}+1=0 & d_{27}=-2 & d_{29}=0 &
\end{array}
\end{equation}

\begin{proof}
Take the particular weights $w_0=0$, $w_1=1$, $w_2=5$, $w_3=25$ (cf. \ref{pesospart}).

The degree of fibers is the same for every base flag that we look. For example, for the flag $\varphi _0$ (\ref{pdv}) the degree is (see somad$\underline{ \ \ }$(0,0,1) in \ref{somad's}
\begin{equation}\label{grauphi0}
-\dfrac{729}{320}d_1 - \dfrac{729}{640}d_2 + \dfrac{729}{1600}d_3 - \dfrac{729}{320}d_4 + \cdots + \dfrac{3125}{192}d_{30} + \dfrac{49642909}{3974400}
\end{equation} 
and for the flag $\varphi _1$,
\begin{equation*}
\varphi_1:\ 
p_1=\{x_0=x_1=x_3=0\} \in \ell_0 = \{x_0=x_1=0\} 
\subset v_1 =
\{x_1=0\} ,
\end{equation*}
the degree is (see somad$\underline{ \ \ }$(1,0,0) in \ref{somad's})
\begin{equation}\label{grauphi1}
-\dfrac{1}{6000}d_1 - \dfrac{1}{12000}d_2 - \dfrac{1}{144000}d_3 - \dfrac{1}{6000}d_4 + \cdots - \dfrac{210827008}{121875}d_{30} + \dfrac{7491488544437}{5070000000}\ .
\end{equation} 

Since the degrees are the same, we gain the equation (\ref{grauphi0})=(\ref{grauphi1}). And we have 24 fixed flags, so actually we get 23 linear equations on the 30 variables $d_i$. This system of equations has the 18 independent linear equations (\ref{relacoesd}), cf. ideal II in \ref{degfibercalc}.
\end{proof}
\end{lema}

\begin{remk}\normalfont
The relations \ref{relacoesd} are enough to compute the Chern classes of the normal bundles of fixed lines, and also for compute the degree of a fiber. To example, for the line $\tilde \ell_1$ we have (\ref{Nl1tildeds})
$$
\mathcal{N}_{(q,p_{\ell_1},v_0')\tilde \ell_1 | \mathbb{Y}_3} = \underbracket{\dfrac{x_0}{x_1}}_{\mathcal{O}_{\tilde{\ell_1}}(d_1)} + \underbracket{\dfrac{x_0^2}{x_1^2}}_{\mathcal{O}_{\tilde{\ell_1}}(d_2)} + \underbracket{\dfrac{x_2}{x_0}}_{\mathcal{O}_{\tilde \ell_1}(d_3)} + \underbracket{\dfrac{x_0}{x_1}}_{\mathcal{O}_{\tilde \ell_1}(d_4)} + \underbracket{\dfrac{x_2}{x_1}}_{\mathcal{O}_{\tilde \ell_1}(d_5)} + \underbracket{\dfrac{x_3}{x_1}}_{\mathcal{O}_{\tilde \ell_1}(d_6)},
$$
Since $d_3=1, \ d_5=0$ and $d_6=0$, we obtain
$$
c_6^T(\mathcal{N}_{\tilde \ell_1 | \mathbb{Y}_3}) =
(w_0-w_1+d_1h)(2w_0-2w_1+d_2h)(w_0-w_1+d_4h)(w_2-w_0+h)(w_2-w_1)(w_3-w_1)
$$
 $c_6^T(\mathcal{N}_{\tilde \ell_1 | \mathbb{Y}_3}) = \mb{w} + \mb{d}h$,
where
$$\mb{w}=2(w_0-w_1)^3(w_2-w_1)(w_3-w_1)(w_2-w_0)$$
$$\mb{d}=[(w_0-w_1)^2(w_2-w_1)(w_3-w_1)]\cdot [\underbrace{(2d_1+d_2+2d_4)}_{=-2}(w_2-w_0)+2(w_0-w_1)],$$
$$\mb{d}=2(w_0-w_1)^2(w_2-w_1)(w_3-w_1)(2w_0-w_1-w_2).$$

\end{remk}

\begin{prop} (Degree of a Fiber) \label{degreefibre}
The degree of the 7-dimensional variety of exceptional foliations  over a fixed flag is \mb{21}.

\begin{proof}
Simply replace the relations (\ref{relacoesd}) in (\ref{grauphi0}), cf. \ref{degfibercalc}
\end{proof}
\end{prop}

When running over all  the complete flags of
$\mathbb{P}^3$ we obtain the full exceptional component of foliations. The contributions on Bott's formula
are computed replacing the power 7 on the numerators by
13, and multiplying the denominators by
$c_6^T(T_{(p,\ell,v)}{\bb F})$, where $(p,\ell,v)$ runs
over the fixpoints of the variety ${\bb F}$ of complete 
flags in $\mathbb{P}^3$ (see (\ref{BOTT}) and (\ref{tangband})). \\

If $\langle x_i,x_j,x_k,x_m \rangle $ is a basis of $(\mathbb{C}^{4})^{\vee}$ and the fixed complete flag is 
$$p_{ijk}=\{x_i=x_j=x_k=0\} \in \ell_{ij}=\{x_i=x_j=0\} \subset v_i=\{x_i=0\},$$
then 
$$T_{(p_{ijk},\ell_{ij},v_i)}{\bb F} = \dfrac{x_j}{x_i} + \dfrac{x_k}{x_i} + \dfrac{x_m}{x_i} + \dfrac{x_k}{x_j} + \dfrac{x_m}{x_j} + \dfrac{x_m}{x_k},$$
and so
$$c_6^T(T_{(p_{ijk},\ell_{ij},v_i)}{\bb F})=(w_j-w_i)(w_k-w_i)(w_m-w_i)(w_k-w_j)(w_m-w_j)(w_m-w_k).$$

\begin{teo} \label{GRAUE3}
The degree of the exceptional component of foliations of codimension one and degree two in $\mathbb{P}^3$ is \mb{168208}.
\end{teo}

\section{A geometric interpretation of the degree}

For a codimension one foliation in $\mathbb{P}^n$, given by the differential form
$$\omega = \sum_{i=0}^{n} A_idx_i,$$
the hyperplane defined by the distribution at a point $p \in \mathbb{P}^n$ is 
$$H=\sum_{i=0}^{n} A_i(p)x_i = 0.$$

Thus, a tangent direction $v=(v_0: \ldots : v_n) \in \mathbb{P}(T_p\mathbb{P}^3)$ lies on this hyperplane $H$ if 
$$\sum_{i=0}^{n} A_i(p) \cdot v_i = 0.$$
This can be thought of as a linear equation on the
coefficients of the $A_i$.
Hence, the point $(p,v) \in \mathbb{P}(T\mathbb{P}^3) $ 
defines a hyperplane in the projective space of distributions.
\begin{equation}\label{condlinear}
\omega(p) \cdot v = 0.
\end{equation}

The equation (\ref{condlinear}) shows that the degree of
a $m$--dimensional component of the space of codimension
one foliations in $\mathbb{P}^n$ can be interpreted as
the number of such foliations that are tangent  to $m$ 
general
directions in $\mathbb{P}^n$.  

In particular, Theorem \ref{GRAUE3} means that there are 168208 exceptional foliations tangent to 13 general directions in $\mathbb{P}^3$.  

%\begin{figure}
%\centering\ifpdf
%\includegraphics{tangente3}\fi
%\caption{There are 168208 exceptional foliations tangent to 13 general directions.}
%\end{figure}

%\include{apendice/apendice}
\chapter{Appendix A}\label{apA}

Here, we present scripts for the software \textit{Macaulay2} to perform some calculations.
You may open a session in
\url{http://habanero.math.cornell.edu:3690/} and
perform the computations just by cut\&paste.
%In order
%to see an explicit output, remove the ``;'' at the
%end of line.

\section{The Standard Exceptional Foliation}\label{primdec}

\begin{verbatim}
Rx=QQ[x_0..x_3]; R = Rx[dx_0..dx_3,SkewCommutative=>true];
X = vars Rx; dX = vars R;
f=x_0^2*x_3-x_0*x_1*x_2+(1/3)*x_1^3; 
g=x_0*x_2-(1/2)*x_1^2;

-- computing df and dg
matrix{{f}}; diff(X,oo); mdf=sub(oo,R); df = determinant(oo * (transpose dX));
matrix{{g}}; diff(X,oo); mdg=sub(oo,R); dg = determinant(oo * (transpose dX));

-- the explicity 1-form
w = 2*g*df - 3*f*dg; w=w//x_0

-- the singular locus. Taking the polinomials A_i
(m,c)=coefficients w;
c=sub(c,Rx); I=ideal c

-- the components of the singular locus
primaryDecomposition I
\end{verbatim}

\section{Resolution of the map}

We start by the resolution of the indeterminacies of the map $(g,f) \mapsto \omega$. This was done by looking over standard neighborhoods that cover $\mathbb{X}'$. We present here two neighborhoods, leaving the others to the reader.

\subsection{Neighborhood $b_3=a_6=1$} 

\subsubsection{Computing $\omega$}\label{compw1}
\begin{verbatim}
Ra=QQ[a_0 .. a_6,b_0..b_3,s_0..v_5];
Rx = Ra[x_0 .. x_3]; R = Rx[dx_0 .. dx_3,SkewCommutative=>true];
X = vars Rx; dX = vars R; X3 = basis(3,Rx); X4 = basis(4,Rx);
f1 = a_0*x_0^3 + a_1*x_0^2*x_1 + a_2*x_0^2*x_2 + a_3*x_0^2*x_3;
f2 = a_4*x_0*x_1^2 + a_5*x_0*x_1*x_2 + a_6*x_1^3;
f=f1+f2;
g = b_0*x_0^2 + b_1*x_0*x_1 + b_2*x_0*x_2 + b_3*x_1^2;

-- the divisibility condition on this neighborhood
f=sub(f,{a_4=>(3/2)*b_1,a_5=>(3/2)*b_2,a_6=>1});
g=sub(g,b_3=>1);

-- codes to compute df and dg
matrix{{f}}; diff(X,oo); mdf=sub(oo,R); df = determinant(oo * (transpose dX));
matrix{{g}}; diff(X,oo); mdg=sub(oo,R); dg = determinant(oo * (transpose dX));
w = 2*g*df - 3*f*dg; 

-- at this point, one can verify if w is divisible by x_0.
-- for this, calculate w % x_0 to verify that this remainder is zero.

w=w // x_0;

\end{verbatim}

\subsubsection{Indeterminacy Locus}\label{obtemJ1}

\begin{verbatim}
-- the following codes will extract the ideal J of the indeterminacy locus
-- execute after loading the codes in previous section
matrix{{w}}; (p,r) = quotientRemainder(oo,dX);
(p,r) = quotientRemainder(sub(transpose p,Rx), X3);
J = ideal p; J = trim J; J=sub(J,Ra);

-- the ideal J_red is the radical of J
-- We mean K=J_red. This is the local ideal of C 
K=radical J;
\end{verbatim}

\subsubsection{Blowup along $C$}\label{blowC1}

\begin{verbatim}
-- the local equations of C are the equations of K
e_0=8*a_0-b_1^3;
e_1=a_2;
e_2=b_2;
e_3=a_3;
e_4=4*b_0-b_1^2;
e_5=4*a_1-3*b_1^2;

-- the substitution to local equation a_2 are done
-- by list m_1 below
m_11={a_0=>(1/8)*b_1^3+(1/8)*a_2*s_0, b_2=>a_2*s_2, a_3=>a_2*s_3};
m_12={b_0=>(1/4)*b_1^2+(1/4)*a_2*s_4, a_1=>(3/4)*b_1^2+(1/4)*a_2*s_5};
m_1=flatten(append(m_11,m_12));

wnt_1=sub(w,m_1);
-- to verify that wn_1t is divisible by a_2, compute 
-- wnt_1 % a_2, and see if this remainder is zero.

wn_1=wnt_1 // a_2;

-- collect the new indeterminacy locus ideal Jn_1
matrix{{wn_1}}; (p,r) = quotientRemainder(oo,dX);
(p,r) = quotientRemainder(sub(transpose p,Rx), X3);
Jn_1 = ideal p; Jn_1 = trim Jn_1; Jn_1=sub(Jn_1,Ra);
\end{verbatim}

\subsubsection{Running other neighboorhods}\label{othviz}

\begin{verbatim}
-- now, we will run over all 6 local equations

for i from 0 to 5 do (
exc=e_i;
ll1={a_0=>(1/8)*(s_0*exc + b_1^3), a_2=>s_1*exc, b_2=>s_2*exc};
ll2={a_3=>s_3*exc, b_0=>(1/4)*(s_4*exc + b_1^2), a_1=>(1/4)*(s_5*exc + 3*b_1^2)};
l=flatten(append(ll1,ll2));
m_i=delete(l_i,l);
Jn_i=sub(J,m_i); 
Jn_i=Jn_i : exc;
wn_i=sub(w,m_i);
wn_i=wn_i // exc; 
)

-- the Jn_i are the new indeterminacy locus ideals, seen locally. 
-- all Jn_i, 0 <= i <= 5, are irreducible and reduced
-- and Jn_3=1 means that the map was solved in this neighborhood

for i from 0 to 5 do (
print(# primaryDecomposition Jn_i, radical Jn_i == Jn_i);     
)

for i from 0 to 5 do (
exc=e_i;
JR_i = ideal(3*s_4-2*s_5, s_3, exc, 3*b_1*s_2-4*s_1, s_0-b_1*s_5);
print(sub(JR_i,s_i=>1)==Jn_i);   
)
\end{verbatim}

\subsubsection{Fibers to Table \ref{x13x12}}\label{fibrasx13x12}

\begin{verbatim}
-- in the previous section, we've computed the expressions 
-- wn_i of w after the blowup along C. Now, we look 
-- at the fixed points for the referred Table2.
M1={a_0=>0,a_1=>0,a_2=>0,a_3=>0,b_0=>0,b_1=>0,b_2=>0};
M2={s_0=>0,s_1=>0,s_2=>0,s_3=>0};
M3={s_4=>0,s_5=>0};
M12=flatten(append(M1,M2));

-- a little caution now: there's a fixed P1 (s_4:s_5).
-- for this reason we've set  separated M3 above.
 
for i from 0 to 5 do (
Fwn_i=sub(wn_i,M12);
if i<4 then Fwn_i=sub(Fwn_i,M3);
print(s'_i,Fwn_i);  
)

-- look closely at the fibers over the fixed P1!
\end{verbatim}

\subsubsection{Fibers to Table \ref{x13x12s2}}\label{fibrasx13x12s2}

\begin{verbatim}
-- Fwn_2=0 means indeterminacy. We continue with equations of R
e_6=3*s_4-2*s_5;
e_7=s_3;
e_8=s_0-b_1*s_5;
e_9=4*s_1-3*b_1;
e_10=b_2;

for i from 6 to 10 do (
exc=e_i;
ll1={s_4=>(1/3)*(v_0*exc+2*s_5), s_3=>v_1*exc, s_0=>v_2*exc+b_1*s_5};
ll2={s_1=>(1/4)*(v_3*exc+3*b_1), b_2=>v_4*exc};
l=flatten(append(ll1,ll2));
m_i=delete(l_(i-6),l);
Jn2_i=sub(Jn_2,m_i); 
Jn2_i=Jn2_i : exc;
wn2_i=sub(wn_2,m_i);
wn2_i=wn2_i // exc; 
)

M4={v_0=>0,v_1=>0,v_2=>0,v_3=>0,v_4=>0};

for i from 6 to 10 do (
Fwn2_i=sub(wn2_i,M12);
Fwn2_i=sub(Fwn2_i,M3);
Fwn2_i=sub(Fwn2_i,M4);
print(v'_(i-6),Fwn2_i);  
)
\end{verbatim}

\subsubsection{Fibers to Table \ref{x13x12pl}}\label{fibrasx13x12pl}

\begin{verbatim}
-- the point pl=(0:0:0:0:1:3/2) has s_4=1, s_5=3/2
-- ideal of next blowup center R is Jn_4 
Jn_4==ideal(2*s_5-3,s_3,2*s_0-3*b_1,4*s_1-3*b_1*s_2,4*b_0-b_1^2)

e_6=2*s_5-3;
e_7=s_3;
e_8=2*s_0-3*b_1;
e_9=4*s_1-3*b_1*s_2;
e_10=4*b_0-b_1^2;

for i from 6 to 10 do (
exc=e_i;
ll1={s_5=>(1/2)*(v_0*exc + 3), s_3=>v_1*exc, s_0=>(1/2)*(v_2*exc + 3*b_1)};
ll2={s_1=>(1/4)*(v_3*exc+3*b_1*s_2), b_0=>(1/4)*(v_4*exc+b_1^2)};
l=flatten(append(ll1,ll2));
m_i=delete(l_(i-6),l);
Jn4_i=sub(Jn_4,m_i); 
Jn4_i=Jn4_i : exc;
wn4_i=sub(wn_4,m_i);
wn4_i=wn4_i // exc; )

-- all Jn4_i=1, 6<=i<=10.

for i from 6 to 10 do (
Fwn4_i=sub(wn4_i,M12);
Fwn4_i=sub(Fwn4_i,s_5=>3/2);
Fwn4_i=sub(Fwn4_i,M4);
print(v'_(i-6),Fwn4_i);  
)

-- Fwn4_i, i=6..10 is the generator of W over (x_1^2,x_1^3,pl,v_(6-i)) 
-- this completes the resolution over (x_1^2,x_1^3).

\end{verbatim}

\subsection{Neighborhood $a_0=b_0=u_1=1$}

\subsubsection{Computing $\omega$ }\label{compw2}
\begin{verbatim}

restart
Ra=QQ[a_0..b_6,s_0..v_5,z_0..z_5];
Rx = Ra[x_0 .. x_3]; R = Rx[dx_0 .. dx_3,SkewCommutative=>true];
X = vars Rx; dX = vars R; X3 = basis(3,Rx); X4 = basis(4,Rx);
f1 = a_0*x_0^3 + a_1*x_0^2*x_1 + a_2*x_0^2*x_2 + a_3*x_0^2*x_3;
f2 = a_4*x_0*x_1^2 + a_5*x_0*x_1*x_2 + a_6*x_1^3;
f = f1+f2;
g = b_0*x_0^2 + b_1*x_0*x_1 + b_2*x_0*x_2 + b_3*x_1^2;

-- first blowup X' --> X
g=sub(g,{b_0=>1, b_2=>u_2*b_1, b_3=>u_3*b_1});
f=sub(f,{a_5=>a_4*u_2, a_6=>(2/3)*a_4*u_3});

matrix{{f}}; diff(X,oo); mdf=sub(oo,R); df = determinant(oo * (transpose dX));
matrix{{g}}; diff(X,oo); mdg=sub(oo,R); dg = determinant(oo * (transpose dX));

-- the 1-form w
w = 2*g*df - 3*f*dg; w=w // x_0;

-- indeterminacy locus of (g,f) --> w
matrix{{w}}; (p,r) = quotientRemainder(oo,dX);
(p,r) = quotientRemainder(sub(transpose p,Rx), X3);
J = ideal p; J = trim J; J=sub(J,Ra);

-- if a_0=0 there's no indeterminacy
trim(sub(J,a_0=>0))
-- this means that a_0=...=a_4=0, impossible.


\end{verbatim}

\subsubsection{Components of reduced indeterminacy locus}\label{compCE}

\begin{verbatim}
-- next, we can set a_0=1 to play with indeterminacy locus

w=sub(w,a_0=>1);
J=sub(J,a_0=>1);
K=trim(radical J);
primaryDecomposition K

-- The primary decomposition of K reveals two ideals.
-- Namely, the codimension one component is C and the
-- codimension two component is E.

-- These two components meet at one point
trim(oo_0+oo_1)

\end{verbatim}

\subsubsection{First blowup center: C}\label{primC}
 
\begin{verbatim}
-- we start with the local ideal of C 
-- and choose u_2 for local equation

-- the list of substitutions is m_5 below:

m_51={b_1=>4*u_3+u_2*s_0, a_1=>6*u_3 + u_2*s_1, a_2=>u_2*s_2};
m_52={a_3=>u_2*s_3, a_4=>12*u_3^2+u_2*s_4};
m_5=flatten(append(m_51,m_52));

-- now, plug it into w
wn_5=sub(w,m_5);

-- verify if the new expression is divisible by u_2
wn_5 % u_2 == 0

-- since the answer is true, may perform the division
wn_5=wn_5 // u_2;

-- the new indeterminacy locus
Jn_5=sub(J,m_5);
Jn_5=trim(Jn_5 : u_2);

-- and the strict transform of K=J_red
Kn_5=sub(K,m_5);
Kn_5=trim(Kn_5 : u_2);
 
-- The strict transform of E
JE=ideal(a_1..a_4,b_1);
SJE=sub(JE,m_5):u_2;
SJE==Kn_5
-- this means that the strict transform E' of E
-- coincides with the strict transform of J_red
\end{verbatim}

\subsubsection{Second blowup center: E}\label{secE}
\begin{verbatim}
-- After the blowup along C, the ruled surface
-- R appears in the reduction of the
-- indeterminacy locus, rad(Jn_5).
-- This reduction consists of E' and R

primaryDecomposition(radical Jn_5);
trim(oo_0+oo_1)
-- this means that this intersection is a curve

-- Now, we'll perform the blowup along E
-- this is done by changes m_7 below:

m_71={s_3=>t_0*s_2, s_0=>(1/3)*(t_2*s_2+2*s_1)};
m_72={u_3=>(1/6)*(t_3*s_2 - s_1*u_2), s_4=>(1/3)*(t_4*s_2-s_1^2*u_2)};
m_7=flatten(append(m_71,m_72));

-- now, put it into wn_5
wn_7=sub(wn_5,m_7);

-- verify if the new expression is divisible by s_2
wn_7 % s_2 == 0

-- since the answer is true, perform the division
wn_7=wn_7 // s_2;

-- the new indeterminacy locus
Jn_7=sub(Jn_5,m_7);
Jn_7=trim(Jn_7 : s_2);

-- this indeterminacy locus Jn_7 is reduced
radical Jn_7 == Jn_7

-- but not irreducible
primaryDecomposition Jn_7

\end{verbatim}

\subsubsection{Third blowup center: R}\label{tercR}

\begin{verbatim}

-- The blowup on v_0=1 is done by changes m_11 below

m_111={t_3=>v_1*u_2+1, t_2=> v_2*u_2, t_0=> v_3*u_2};
m_112={s_1=>(1/3)*(v_4*u_2+2*t_4)};
m_11=flatten(append(m_111,m_112));

-- now, put it into wn_7
wn_11=sub(wn_7,m_11);

-- verify if the new expression is divisible by u_2
wn_11 % u_2 == 0

-- since the answer is true, perform the division
wn_11=wn_11 // u_2;

-- the new indeterminacy locus
Jn_11=sub(Jn_7,m_11);
Jn_11=trim(Jn_11 : u_2);

-- this indeterminacy locus Jn_11 is reduced and irreducible !!!!
radical Jn_11 == Jn_11
primaryDecomposition Jn_11 
Jn_11

\end{verbatim}

\subsubsection{Fourth blowup center: L}\label{quartL}

\begin{verbatim}

-- we'll run all 6 standard neighborhoods z_i=1, 0<=i<=5

e_16=s_2;
e_17=v_1;
e_18=v_2;
e_19=v_3;
e_20=v_4;
e_21=t_4;

-- in the following codes, the list m_(16+i) does the
-- necessary substitutions to the blowup on z_i=1 

for i from 16 to 21 do (
exc=e_i;
ll1={s_2=>z_0*exc, v_1=>z_1*exc, v_2=>z_2*exc, v_3=>z_3*exc};
ll2={v_4=>z_4*exc, t_4=>z_5*exc};
l=flatten(append(ll1,ll2));
m_i=delete(l_(i-16),l);
Jn_i=sub(Jn_11,m_i);  
Jn_i=Jn_i : exc;    
print Jn_i;      
)

-- the indeterminacy locus on z_i=1 is given by Jn_(i+16)
-- as they are all equal to 1, the map was solved.
\end{verbatim}

\subsubsection{Running other neighboorhods}\label{othviz2}

\begin{verbatim}

-- equations of C

e_0=b_1-4*u_3;
e_1=a_1-6*u_3;
e_2=a_2;
e_3=a_3;
e_4=a_4-12*u_3^2;
e_5=u_2;

-- we'll run all six possible choices of local
-- equation for exceptional divisor.
-- Jn_i is the new indeterminacy locus for s_i=1.

for i from 0 to 5 do (
exc=e_i;
ll1={b_1=>s_0*exc + 4*u_3, a_1=>s_1*exc + 6*u_3, a_2=>s_2*exc};
ll2={a_3=>s_3*exc, a_4=>s_4*exc + 12*u_3^2, u_2=>s_5*exc};
l=flatten(append(ll1,ll2));
m_i=delete(l_i,l);
Jn_i=sub(J,m_i);
Jn_i=trim(Jn_i : exc);
Kn_i=sub(K,m_i);
Kn_i=trim(Kn_i : exc);
wn_i=sub(w,m_i); 
wn_i=wn_i // exc; )

-- wn_i is the expression of w in the new affine coordinates
-- notice that Jn_3=1 means that there's 
-- no indeterminacy when s_3=1.

-- the second blowup center is Kn_i. Note Kn_2 = 1.
-- from now on, there are many neighborhoods to consider.
-- we continue by neighborhood s_0=1, leaving the others 
-- to the reader since they are similar

-- radical of Jn_0 is composed by two components:
primaryDecomposition Jn_0


-- the components E (represented by Kn_0) and
-- the component R (the other ideal in p.D. Jn_0) 
-- we take the ideal of R

JR_0 = (radical Jn_0) : Kn_0;

-- equations of Kn_0

e_6=s_4-3*u_3;
e_7=s_3;
e_8=s_2;
e_9=2*s_1-3;
e_10=b_1;


for i from 6 to 10 do (
exc=e_i;
ll1={s_4=>t_0*exc + 3*u_3, s_3=>t_1*exc, s_2=>t_2*exc};
ll2={s_1=>(1/2)*(t_3*exc + 3), b_1=>t_4*exc};
l=flatten(append(ll1,ll2));
m_i=delete(l_(i-6),l);
Jn0_i=sub(Jn_0,m_i);
Jn0_i=trim(Jn0_i : exc);
Kn0_i=sub(Kn_0,m_i);
Kn0_i=trim(Kn0_i : exc);
JRn_i=sub(JR_0,m_i);
JRn_i=trim(JRn_i : exc);
wn0_i=sub(wn_0,m_i); 
wn0_i=wn0_i // exc;    
)

-- the ideal JRn_i represents the transform of R
-- 6 <= i <= 10. The new indeterminacy locus is Jn0_i.
-- The magic now is the equality in all neighborhoods

for i from 6 to 10 do print (Jn0_i == JRn_i); 

-- Henceforth, the next blowup along R will solve
-- the map in any neighboorhood, without another blowup!
\end{verbatim}

\subsubsection{Fibers to Table \ref{x01}}\label{fibrasx01}

\begin{verbatim}
M1={a_1=>0,a_2=>0,a_3=>0,a_4=>0,b_1=>0,u_2=>0,u_3=>0};
M2={s_2=>0,s_3=>0,s_4=>0,s_5=>0};
M3={s_0=>0,s_1=>0};
M12=flatten(append(M1,M2));

-- a little caution now: there's a fixed P1 (s_0:s_1).
-- it is for this reason the separated M3 above.

for i from 0 to 5 do (
Fwn_i=sub(wn_i,M12);
if i>1 then Fwn_i=sub(Fwn_i,M3);
print(s'_i,Fwn_i);  
)

\end{verbatim}

\subsubsection{Fibers to Table \ref{x01s5}}\label{fibrasx01s5}

\begin{verbatim}
-- equations of E' on s_5=1

e_6=s_3;
e_7=s_2;
e_8=3*s_0-2*s_1;
e_9=6*u_3+s_1*u_2;
e_10=3*s_4+s_1^2*u_2;

-- equations of R
JR_5 = (radical Jn_5) : Kn_5;

for i from 6 to 10 do (
exc=e_i;
ll1={s_3=>t_0*exc, s_2=>t_1*exc, s_0=>(1/3)*(t_2*exc+2*s_1)};
ll2={u_3=>(1/6)*(t_3*exc-s_1*u_2), s_4=>(1/3)*(t_4*exc-s_1^2*u_2)};
l=flatten(append(ll1,ll2));
m_i=delete(l_(i-6),l);
Jn5_i=sub(Jn_5,m_i);
Jn5_i=trim(Jn5_i : exc);
Kn5_i=sub(Kn_5,m_i);
Kn5_i=trim(Kn5_i : exc);
JRn_i=sub(JR_5,m_i);
JRn_i=trim(JRn_i : exc);
wn5_i=sub(wn_5,m_i); 
wn5_i=wn5_i // exc;    
)

-- the JRn_i are the ideals of next blowup center R 

M41={t_0=>0,t_2=>0,t_4=>0};
M42={t_1=>0,t_3=>0};

for i from 6 to 10 do (
Fwn5_i=sub(wn5_i,M41);
Fwn5_i=sub(Fwn5_i,M12);
Fwn5_i=sub(Fwn5_i,M3);
teste_i=(i-6)*(i-8)*(i-10);
if teste_i==0 then Fwn5_i=sub(Fwn5_i,M42);
print(t'_(i-6),Fwn5_i);  
)

\end{verbatim}

\subsubsection{Fibers to Table \ref{x01s5pl3}}\label{fibrasx01s5pl3}

\begin{verbatim}
-- Equations of R
e_11=u_2;
e_12=t_3-1;
e_13=t_2;
e_14=t_0;
e_15=3*s_1-2*t_4;

for i from 11 to 15 do (
exc=e_i;
ll1={u_2=>v_0*exc, t_3=>v_1*exc+1, t_2=>v_2*exc};
ll2={t_0=>v_3*exc, s_1=>(1/3)*(v_4*exc+2*t_4)};
l=flatten(append(ll1,ll2));
m_i=delete(l_(i-11),l);
Jn57_i=sub(Jn5_7,m_i);
Jn57_i=trim(Jn57_i : exc);
wn57_i=sub(wn5_7,m_i); 
wn57_i=wn57_i // exc;    
)

M5=apply(5,i->v_i=>0);
-- build list {v_i=>0}
-- remember that t_3=1 on pl_3

for i from 11 to 15 do (
Fwn57_i=sub(wn57_i,M41);
Fwn57_i=sub(Fwn57_i,M12);
Fwn57_i=sub(Fwn57_i,M3);
Fwn57_i=sub(Fwn57_i,M5);
Fwn57_i=sub(Fwn57_i,t_3=>1);
print(v'_(i-11),Fwn57_i);  
)

\end{verbatim}

\subsubsection{Fibers to Table \ref{x01s5pl3v0}}\label{fibrasx01s5pl3v0}

\begin{verbatim}
-- The blowup on neighborhood s_0=1 was done in 
-- "running other neighborhoods". 
-- fibers on neighborhood s_0=1 are wn0_i, 6<=i<=10

-- there's a fixed P^1={(t_0:t_4)}
M71={t_1=>0,t_2=>0,t_3=>0};
M72={t_0=>0,t_4=>0};

for i from 6 to 10 do (
Fwn0_i=sub(wn0_i,M12);
Fwn0_i=sub(Fwn0_i,M3);
Fwn0_i=sub(Fwn0_i,M71);
teste_i=(i-7)*(i-8)*(i-9);
if teste_i==0 then Fwn0_i=sub(Fwn0_i,M72);
print(t'_(i-6),Fwn0_i);  
)

\end{verbatim}

\begin{verbatim}

-- indeterminacy locus Jn57_11
-- Equations of L
e_16=s_2;
e_17=v_1;
e_18=v_2;
e_19=v_3;
e_20=v_4;
e_21=t_4;

for i from 16 to 21 do (
exc=e_i;
ll1={s_2=>z_0*exc, v_1=>z_1*exc, v_2=>z_2*exc};
ll2={v_3=>z_3*exc, v_4=>z_4*exc, t_4=>z_5*exc};
l=flatten(append(ll1,ll2));
m_i=delete(l_(i-16),l);
Jnz_i=sub(Jn57_11,m_i);
Jnz_i=trim(Jnz_i : exc);
wnz_i=sub(wn57_11,m_i); 
wnz_i=wnz_i // exc;    
)

M6=apply(6,i->z_i=>0);

for i from 16 to 21 do (
Fwnz_i=sub(wnz_i,M41);
Fwnz_i=sub(Fwnz_i,M12);
Fwnz_i=sub(Fwnz_i,M3);
Fwnz_i=sub(Fwnz_i,M5);
Fwnz_i=sub(Fwnz_i,M6);
print(z'_(i-16),Fwnz_i);  
)

-- finally done the fixed points over 
-- (x_0^2,x_0x_1,x_0^3,s_5')

\end{verbatim}

\subsubsection{Fibers to Table \ref{x01pl2}}\label{fibrasx01pl2}

\begin{verbatim}

-- The blowup on neighborhood s_0=1 was done in 
-- "running other neighborhoods". 
-- fibers on neighborhood s_0=1 are wn0_i, 6<=i<=10

-- there's a fixed P^1={(t_0:t_4)}
M71={t_1=>0,t_2=>0,t_3=>0};
M72={t_0=>0,t_4=>0};

for i from 6 to 10 do (
Fwn0_i=sub(wn0_i,M12);
Fwn0_i=sub(Fwn0_i,M3);
Fwn0_i=sub(Fwn0_i,M71);
teste_i=(i-7)*(i-8)*(i-9);
if teste_i==0 then Fwn0_i=sub(Fwn0_i,M72);
print(t'_(i-6),Fwn0_i);  
)

\end{verbatim}

\subsubsection{Fibers to Table \ref{x01pl2pl4}}\label{fibrasx01pl2pl4}

\begin{verbatim}

-- Equations of R
e_11=3*t_4-8;
e_12=t_3;
e_13=t_1;
e_14=t_2-4*s_5;
e_15=2*s_4-9*u_3;

for i from 11 to 15 do (
exc=e_i;
ll1={t_4=>(1/3)*(v_0*exc+8), t_3=>v_1*exc, t_1=>v_2*exc};
ll2={t_2=>v_3*exc+4*s_5, s_4=>(1/2)*(v_4*exc+9*u_3)};
l=flatten(append(ll1,ll2));
m_i=delete(l_(i-11),l);
Jn06_i=sub(Jn0_6,m_i);
Jn06_i=trim(Jn06_i : exc);
wn06_i=sub(wn0_6,m_i); 
wn06_i=wn06_i // exc;    
)

for i from 11 to 15 do (
Fwn06_i=sub(wn06_i,M12);
Fwn06_i=sub(Fwn06_i,M5);
Fwn06_i=sub(Fwn06_i,t_4=>(8/3));
print(v'_(i-11),Fwn06_i);  
)

-- Finally solved over (x_0^2,x_0x_1,x_0^3) !!!!

\end{verbatim}

\subsection{Neighborhood $a_0=b_0=u_2=1$.}\label{a0b0u2=1}

\begin{verbatim}
restart;
Ra=QQ[a_0 .. a_6,b_0..b_3,s_0..v_5,z_0..z_5];
Rx = Ra[x_0 .. x_3]; R = Rx[dx_0 .. dx_3,SkewCommutative=>true];
X = vars Rx; dX = vars R; X3 = basis(3,Rx); X4 = basis(4,Rx);
f1 = a_0*x_0^3 + a_1*x_0^2*x_1 + a_2*x_0^2*x_2 + a_3*x_0^2*x_3;
f2= a_4*x_0*x_1^2 + a_5*x_0*x_1*x_2 + a_6*x_1^3;
f=f1+f2;
g = b_0*x_0^2 + b_1*x_0*x_1 + b_2*x_0*x_2 + b_3*x_1^2;

g=sub(g,{b_0=>1, b_1=>u_1*b_2, b_3=>u_3*b_2});
f=sub(f,{a_4=>a_5*u_1, a_6=>(2/3)*a_5*u_3});

matrix{{f}}; diff(X,oo); mdf=sub(oo,R); df = determinant(oo * (transpose dX));
matrix{{g}}; diff(X,oo); mdg=sub(oo,R); dg = determinant(oo * (transpose dX));

w = 2*g*df - 3*f*dg; w=w // x_0;

matrix{{w}}; (p,r) = quotientRemainder(oo,dX);
(p,r) = quotientRemainder(sub(transpose p,Rx), X3);
J = ideal p; J = trim J; J=sub(J,Ra);

-- no indeterminacies when a_0=0
trim(sub(J,a_0=>0))==ideal(a_1,a_2,a_3,a_5)

-- this means (a_0:a_1:a_2:a_3:a_5)=(0:0:0:0:0), impossible.
-- Hence, to deal with indeterminacy locus we can take a_0=1.

w=sub(w,a_0=>1); J=trim(sub(J,a_0=>1));

-- The indeterminacy locus J is non-reduced. But the reduction presents a 
-- single component E (this neighborhood is outside 
-- of first blowup center C)

K=radical J;

-- equations of blowup center E

e_0=a_1;
e_1=a_2;
e_2=a_3;
e_3=a_5;
e_4=b_2;

for i from 0 to 4 do (
exc=e_i;
ll1={a_1=>t_0*exc, a_2=>t_1*exc, a_3=>t_2*exc};
ll2={a_5=>t_3*exc, b_2=>t_4*exc};
l=flatten(append(ll1,ll2));
m_i=delete(l_i,l);
Jn_i=sub(J,m_i); 
Jn_i=Jn_i : exc;
wn_i=sub(w,m_i);
wn_i=wn_i // exc; )

-- solved for i=2,3. Continue over t_4=1
-- equations of next blowup center Jn_4
-- the line L

e_5=u_3;
e_6=t_3;
e_7=t_2;
e_8=2*t_1-3;
e_9=2*t_0-3*u_1;
e_10=b_2; 

for i from 5 to 10 do (
exc=e_i;
ll1={u_3=>z_0*exc, t_3=>z_1*exc, t_2=>z_2*exc};
ll2={t_1=>(1/2)*(z_3*exc+3), t_0=>(1/2)*(z_4*exc+3*u_1), b_2=>z_5*exc};
l=flatten(append(ll1,ll2));
m_i=delete(l_(i-5),l);
Jn4_i=sub(Jn_4,m_i); 
Jn4_i=Jn4_i : exc;
wn4_i=sub(wn_4,m_i);
wn4_i=wn4_i // exc; )

-- solved everywhere, since Jn4_i=1, for 5 <= i <= 10.

\end{verbatim}

\subsubsection{Fibers to Table \ref{x02}}\label{fibrasx02}

\begin{verbatim}

-- the fixed points. Note the P^1={(t_1:t_4)}

M1={a_1=>0,a_2=>0,a_3=>0,a_5=>0,b_2=>0,u_1=>0,u_3=>0};
M2={t_0=>0,t_2=>0,t_3=>0};
M3={t_1=>0,t_4=>0};
M12=flatten(append(M1,M2));

for i from 0 to 4 do (
Fwn_i=sub(wn_i,M12);
teste_i=i*(i-2)*(i-3);
if teste_i==0 then Fwn_i=sub(Fwn_i,M3);
print(t'_i,Fwn_i);  
)

\end{verbatim}

\subsubsection{Fibers to Table \ref{x02pl5}}\label{fibrasx02pl5}

\begin{verbatim}

-- neighborhood t_4=1. In pl_5 have t_1=3/2.
M4=apply(6,i->z_i=>0);

for i from 5 to 10 do (
Fwn4_i=sub(wn4_i,M12);
Fwn4_i=sub(Fwn4_i,t_1=>(3/2));
Fwn4_i=sub(Fwn4_i,M4);
print(z'_(i-5),Fwn4_i);  
)

\end{verbatim}

\subsection{Neighborhood $a_0=b_0=u_3=1$.}\label{a0b0u3=1}

\begin{verbatim}

Ra=QQ[a_0 .. a_6,b_0..b_3,s_0..v_5,z_0..z_5];
Rx = Ra[x_0 .. x_3]; R = Rx[dx_0 .. dx_3,SkewCommutative=>true];
X = vars Rx; dX = vars R; X3 = basis(3,Rx); X4 = basis(4,Rx);
f1 = a_0*x_0^3 + a_1*x_0^2*x_1 + a_2*x_0^2*x_2 + a_3*x_0^2*x_3;
f2= a_4*x_0*x_1^2 + a_5*x_0*x_1*x_2 + a_6*x_1^3;
f=f1+f2;
g = b_0*x_0^2 + b_1*x_0*x_1 + b_2*x_0*x_2 + b_3*x_1^2;

g=sub(g,{b_0=>1, b_1=>u_1*b_3, b_2=>u_2*b_3});
f=sub(f,{a_4=>(3/2)*a_6*u_1, a_5=>(3/2)*a_6*u_2});

matrix{{f}}; diff(X,oo); mdf=sub(oo,R); df = determinant(oo * (transpose dX));
matrix{{g}}; diff(X,oo); mdg=sub(oo,R); dg = determinant(oo * (transpose dX));
w = 2*g*df - 3*f*dg; w=w // x_0;

matrix{{w}}; (p,r) = quotientRemainder(oo,dX);
(p,r) = quotientRemainder(sub(transpose p,Rx), X3);
J = ideal p; J = trim J; J=sub(J,Ra);

-- no indeterminacies when a_0=0
trim(sub(J,a_0=>0))==ideal(a_1,a_2,a_3,a_6)

-- this means (a_0:a_1:a_2:a_3:a_6)=(0:0:0:0:0), impossible.
-- Hence, to deal with indeterminacy locus we can take a_0=1.

w=sub(w,a_0=>1); J=trim(sub(J,a_0=>1));

-- the first blowup center is C, of ideal
MMM=ideal(u_2,a_3,a_2,9*b_3-a_1^2,a_1*u_1-6,a_1^3-27*a_6)

-- the second blowup center is (the transform) JE, the ideal of E
JE= (radical J) : MMM;

-- the p.D. of radical J is composed by these two ideals MMM and JE
primaryDecomposition radical J

-- choose a_1u_1-6 for the exceptional equation
exc = a_1*u_1 - 6;

m_41={u_2=>s_0*exc, a_3=>s_1*exc, a_2=>s_2*exc};
m_42={b_3=>(1/9)*(s_3*exc + a_1^2), a_6=>(1/27)*(s_5*exc + a_1^3)};
m_4=flatten(append(m_41,m_42));

-- new indeterminacy locus 
Jn_4=sub(J,m_4);
Jn_4=Jn_4 : exc;

-- the second blowup center E
JEn_4=sub(JE,m_4);
JEn_4=JEn_4 : exc;

wn_4=sub(w,m_4);
wn_4 = wn_4 // exc;


-- equations of E'

e_6=s_3;
e_7=a_1;
e_8=s_2;
e_9=s_1;
e_10=s_5;

for i from 6 to 10 do (
exc=e_i;
ll1={s_3=>t_0*exc, a_1=>t_1*exc, s_2=>t_2*exc};
ll2={s_1=>t_3*exc, s_5=>t_4*exc};
l=flatten(append(ll1,ll2));
m_i=delete(l_(i-6),l);
Jn4_i=sub(Jn_4,m_i); 
Jn4_i=Jn4_i : exc;
wn4_i=sub(wn_4,m_i);
wn4_i=wn4_i // exc; )

-- After this two blowups, the indeterminacy locus becomes reduced
for i from 6 to 10 do print(radical Jn4_i == Jn4_i);

-- moreover, Jn4_i is irreducible, 6 <= i <= 10
-- note that Jn4_9=1 (the map was solved here), thus 
-- # primaryDecomposition Jn4_i will return 0. 
-- The ideals Jn4_i, i \not= 9, represents
-- the ruled surface R, and so codim = 5.

for i from 6 to 10 do print(# primaryDecomposition Jn4_i, codim Jn4_i);

-- Thus, the next blowup along R will solve the map.

-- we leave to the reader as an exercise the verifications over
-- the neighborhood [b_1=1]. Over [b_2=1] there's no indeterminacies.
\end{verbatim}

\subsubsection{Fibers to Table \ref{x12}}\label{fibrasx12}

\begin{verbatim}

-- since the indeterminacy point q is outside of first 
-- blowup center C, we will go directly to equations of E
e_0=b_3;
e_1=a_1;
e_2=a_2;
e_3=a_3;
e_4=a_6;

M1={a_1=>0,a_2=>0,a_3=>0,a_6=>0,b_3=>0,u_1=>0,u_2=>0};
M2=apply(5,i->t_i=>0);
M12=flatten(append(M1,M2));

for i from 0 to 4 do (
exc=e_i;
ll1={b_3=>t_0*exc, a_1=>t_1*exc, a_2=>t_2*exc};
ll2={a_3=>t_3*exc, a_6=>t_4*exc};
l=flatten(append(ll1,ll2));
m_i=delete(l_i,l);
wn_i=sub(w,m_i);
wn_i=wn_i // exc; 
Fwn_i=sub(wn_i,M12);
print(t_i, Fwn_i);
)

-- solved in all fixed points.
\end{verbatim}

\subsection{Neighborhood $b_2=1$}\label{b2=1}

\begin{verbatim}
Ra=QQ[a_0 .. a_6,b_0..b_3,s_0..v_5,z_0..z_5];
Rx = Ra[x_0 .. x_3]; R = Rx[dx_0 .. dx_3,SkewCommutative=>true];
X = vars Rx; dX = vars R; X3 = basis(3,Rx); X4 = basis(4,Rx);
f1 = a_0*x_0^3 + a_1*x_0^2*x_1 + a_2*x_0^2*x_2 + a_3*x_0^2*x_3;
f2= a_4*x_0*x_1^2 + a_5*x_0*x_1*x_2 + a_6*x_1^3;
f=f1+f2;
g = b_0*x_0^2 + b_1*x_0*x_1 + b_2*x_0*x_2 + b_3*x_1^2;

-- divisibiliy condition
g=sub(g,b_2=>1);
f=sub(f,{a_4=>a_5*b_1, a_6=>(2/3)*a_5*b_3});

matrix{{f}}; diff(X,oo); mdf=sub(oo,R); df = determinant(oo * (transpose dX));
matrix{{g}}; diff(X,oo); mdg=sub(oo,R); dg = determinant(oo * (transpose dX));
w = 2*g*df - 3*f*dg; w=w // x_0;

matrix{{w}}; (p,r) = quotientRemainder(oo,dX);
(p,r) = quotientRemainder(sub(transpose p,Rx), X3);
J = ideal p; J = trim J; J=sub(J,Ra);

-- (a_0:a_1:a_2:a_3:a_5)=(0:0:0:0:0) is impossible. This
-- means that there's no ideterminacies when b_2=1.

\end{verbatim}

\subsection{Fibers to Table \ref{Fixosem}}\label{fibrasY}

\begin{verbatim}
restart;
Ra=QQ[a_0 .. a_6,b_0..b_3,s_0..v_5,z_0..z_5];
Rx = Ra[x_0 .. x_3]; R = Rx[dx_0 .. dx_3,SkewCommutative=>true];
X = vars Rx; dX = vars R; X3 = basis(3,Rx); X4 = basis(4,Rx);
f1 = a_0*x_0^3 + a_1*x_0^2*x_1 + a_2*x_0^2*x_2 + a_3*x_0^2*x_3;
f2= a_4*x_0*x_1^2 + a_5*x_0*x_1*x_2 + a_6*x_1^3;
f=f1+f2;
g = b_0*x_0^2 + b_1*x_0*x_1 + b_2*x_0*x_2 + b_3*x_1^2;

Za=apply(7,i->a_i=>0);

-- g=x_0x_1

ggg=sub(g,{b_0=>0, b_1=>1, b_2=>0, b_3=>0});
fff=sub(f,{a_5=>0, a_6=>0});
matrix{{fff}}; diff(X,oo); mdf=sub(oo,R); df = determinant(oo * (transpose dX));
matrix{{ggg}}; diff(X,oo); mdg=sub(oo,R); dg = determinant(oo * (transpose dX));
w = 2*ggg*df - 3*fff*dg; w=w // x_0;

for i from 0 to 4 do (
ff_i=sub(fff,a_i=>1);
ff_i=sub(ff_i,Za);
ww_i=sub(w,a_i=>1);
ww_i=sub(ww_i,Za);
print(ggg,ff_i,ww_i); )

-- the print corresponds to the three columns of the table

-- g=x_0x_2

ggg=sub(g,{b_0=>0, b_1=>0, b_2=>1, b_3=>0});
fff=sub(f,{a_4=>0, a_6=>0});
matrix{{fff}}; diff(X,oo); mdf=sub(oo,R); df = determinant(oo * (transpose dX));
matrix{{ggg}}; diff(X,oo); mdg=sub(oo,R); dg = determinant(oo * (transpose dX));
w = 2*ggg*df - 3*fff*dg; w=w // x_0;

for i from 0 to 5 do (
ff_i=sub(fff,a_i=>1);
ff_i=sub(ff_i,Za);
ww_i=sub(w,a_i=>1);
ww_i=sub(ww_i,Za);
if i!=4 then print(ggg,ff_i,ww_i); )

-- g=x_1^2

ggg=sub(g,{b_0=>0, b_1=>0, b_2=>0, b_3=>1});
fff=sub(f,{a_4=>0, a_5=>0});
matrix{{fff}}; diff(X,oo); mdf=sub(oo,R); df = determinant(oo * (transpose dX));
matrix{{ggg}}; diff(X,oo); mdg=sub(oo,R); dg = determinant(oo * (transpose dX));
w = 2*ggg*df - 3*fff*dg; w=w // x_0;

for i from 0 to 6 do (
ff_i=sub(fff,a_i=>1);
ff_i=sub(ff_i,Za);
ww_i=sub(w,a_i=>1);
ww_i=sub(ww_i,Za);
if i^2-9*i!=-20 then print(ggg,ff_i,ww_i); )

-- g=x_0^2, g'=x_0x_1

ggg=sub(g,{b_0=>1, b_1=>0, b_2=>0, b_3=>0});
fff=sub(f,{a_5=>0, a_6=>0});
matrix{{fff}}; diff(X,oo); mdf=sub(oo,R); df = determinant(oo * (transpose dX));
matrix{{ggg}}; diff(X,oo); mdg=sub(oo,R); dg = determinant(oo * (transpose dX));
w = 2*ggg*df - 3*fff*dg; w=w // x_0;

for i from 0 to 4 do (
ff_i=sub(fff,a_i=>1);
ff_i=sub(ff_i,Za);
ww_i=sub(w,a_i=>1);
ww_i=sub(ww_i,Za);
print(ggg,ff_i,ww_i); )

-- g=x_0^2, g'=x_0x_2

ggg=sub(g,{b_0=>1, b_1=>0, b_2=>0, b_3=>0});
fff=sub(f,{a_4=>0, a_6=>0});
matrix{{fff}}; diff(X,oo); mdf=sub(oo,R); df = determinant(oo * (transpose dX));
matrix{{ggg}}; diff(X,oo); mdg=sub(oo,R); dg = determinant(oo * (transpose dX));
w = 2*ggg*df - 3*fff*dg; w=w // x_0;

for i from 0 to 5 do (
ff_i=sub(fff,a_i=>1);
ff_i=sub(ff_i,Za);
ww_i=sub(w,a_i=>1);
ww_i=sub(ww_i,Za);
if i!=4 then print(ggg,ff_i,ww_i); )

-- g=x_0^2, g'=x_1^2

ggg=sub(g,{b_0=>1, b_1=>0, b_2=>0, b_3=>0});
fff=sub(f,{a_4=>0, a_5=>0});
matrix{{fff}}; diff(X,oo); mdf=sub(oo,R); df = determinant(oo * (transpose dX));
matrix{{ggg}}; diff(X,oo); mdg=sub(oo,R); dg = determinant(oo * (transpose dX));
w = 2*ggg*df - 3*fff*dg; w=w // x_0;

for i from 0 to 6 do (
ff_i=sub(fff,a_i=>1);
ff_i=sub(ff_i,Za);
ww_i=sub(w,a_i=>1);
ww_i=sub(ww_i,Za);
if i^2-9*i!=-20 then print(ggg,ff_i,ww_i); )
\end{verbatim}

\section{Degree of a fiber}\label{degfib}

Now, the script to compute the sum of all
contributions. Although similar, first we will present the degree of a fiber for the sake of clarity, and then the degree of the exceptional component.

\begin{verbatim}
Ra=QQ[d_1..d_30];
R=Ra[x_0..x_3];
B={2*x_0,x_0+x_1,x_0+x_2,2*x_1};
A={3*x_0,2*x_0+x_1,2*x_0+x_2,2*x_0+x_3};
A0_1=append(A,(x_0+2*x_1));
A0_2=append(A,(x_0+x_1+x_2));
A0_3=append(A,(3*x_1));

-- in A and B we put the weights for 
-- cubics and quadrics
\end{verbatim}

\subsection{Contributions on Table \ref{Fixosem}}

\begin{verbatim}
-- contributions in Y ---> X
-- here, B instead of X and A instead of Y

-- T_(01)B

for k from 1 to 3 do (
B0_k=delete(B_k,B);
TB0_k=1;
for i from 0 to 2 do (
TB_(k,i)=(B0_k)_i-B_k;
TB0_k=TB0_k * TB_(k,i); ); )

-- T_pA_(01)


for k from 1 to 3 do (
for i from 0 to 4 do (
TA0_(k,i)=1;
FA0_(k,i)=delete((A0_k)_i,A0_k);
for j from 0 to (#(FA0_(k,i)) - 1) do (
TA_(k,i,j)=(FA0_(k,i))_j - (A0_k)_i;
TA0_(k,i)=TA0_(k,i) * TA_(k,i,j); ); ); )

for k from 1 to 3 do (
for i from 0 to 4 do ( 
de_(k,i)=TA0_(k,i) * TB0_k; 
nu_(k,i)=B_k+(A0_k)_i-x_0;
); )

-- the quotient nu_(k,i)^7 / de_(k,i) represents the contribution of 
-- the pair (k,i), 1<=k<=3, 0<=i<=4, (k,i) \not= (3,4), in the order
-- k: (x_0x_1,x_0x_2,x_1^2) = (1,2,3)
-- i: (x_0^3,x_0^2x_1,x_0^2x_2,x_0^2x_3) = (0,1,2,3)
-- i=4 : x_0x_1^2, or x_0x_1x_2, or x_1^3, depending on  k.
-- note that the pair (k,i)=(3,4) does not represent anything
-- since it is an indeterminacy point.

-- contributions in Y -----> X'

B00r=delete(B_0,B);
for k from 1 to 3 do (
B00_k=delete(B_k,B00r);
TB00_k=(B_k-B_0)*((B00_k)_0-B_k)*((B00_k)_1-B_k);
for i from 0 to 4 do (
de0_(k,i)=(TA0_(k,i)*TB00_k); 
nu0_(k,i)=B_0+(A0_k)_i-x_0;
); )

-- the quotient nu0_(k,i)^7 / de0_(k,i) 
-- represents the contribution of (x_0^2,k,i), in the order
-- k: (x_0x_1,x_0x_2,x_1^2) = (1,2,3)
-- i: (x_0^3,x_0^2x_1,x_0^2x_2,x_0^2x_3) = (0,1,2,3)
-- i=4 : x_0x_1^2, ou x_0x_1x_2, ou x_1^3, depending on the k.
-- notice that (x_0^2,k,0) does not represent anything
-- since they are indeterminacies. 
-- 1<=k<=3, 1<=i<=4.

\end{verbatim}

\subsection{Contributions on Table \ref{x13x12}}

\begin{verbatim}
-- contributions over Y_1
-- point (3,4)=x_1^2,x_1^3. Adjusting the tangent space

pesosB={TB_(3,0),TB_(3,1),TB_(3,2)};
pesosA={TA_(3,4,0),TA_(3,4,1),TA_(3,4,2),TA_(3,4,3)};
pesos=flatten(append(pesosA,pesosB));
normal=delete(x_0-x_1,pesos);

-- 3 fixed points. The weights on the table are

curva_1=3*x_0-3*x_1;
curva_2=2*x_0+x_2-3*x_1;
curva_3=2*x_0+x_3-3*x_1;

for k from 1 to 3 do (
Cv_k=delete(curva_k,normal);
TC_(3,4,k)=(x_0-x_1)*(curva_k);
for i from 0 to (#Cv_k - 1) do (
TC1_(k,i)=(Cv_k)_i - curva_k;
TC_(3,4,k)=TC_(3,4,k) * TC1_(k,i); );
de_(3,4,k)=TC_(3,4,k);
nu_(3,4,k)=5*x_1 - x_0 + curva_k; )

-- the quotient nu_(3,4,k)^7 / de_(3,4,k) 
-- represents the contribution of the three isolated
-- fixed points solved over the pair (x_1^2, x_1^3)
-- 1<=k<=3
\end{verbatim}

\subsection{Contributions on Table \ref{x01}}

\begin{verbatim}

normalB={(B00_1)_0-B_1 , (B00_1)_1-B_1};
normalA={TA_(1,0,0),TA_(1,0,1),TA_(1,0,2),TA_(1,0,3)};
normal=flatten(append(normalA,normalB));

-- 3 fixed points. The weights on the table are
-- (curva_4 for future use in the next exceptional)

curva_1=-x_0+x_2;
curva_2=-x_0+x_3;
curva_3=-2*x_0+2*x_1;
curva_4=-x_1+x_2;

for k from 1 to 4 do (
Cv_k=delete(curva_k,normal);
TC_(1,0,k)=(x_1-x_0)*(curva_k);
for i from 0 to (#Cv_k - 1) do (
TC1_(k,i)=(Cv_k)_i - curva_k;
TC_(1,0,k)=TC_(1,0,k) * TC1_(k,i); );
de_(1,0,k)=(TC_(1,0,k));
nu_(1,0,k)=4*x_0 + curva_k; )

-- the quotient nu_(1,0,k)^7 / de_(1,0,k) 
-- represents the contribution of three isolated fixed points solved
-- over the point (x_0^2,x_0x_1,x_0^3), 1<=k<=3
\end{verbatim}

\subsection{Contributions on Table \ref{x01s5}}

\begin{verbatim}
TEt=(x_2-x_1)*(2*x_1-x_0-x_2);
normal={TC1_(4,0),TC1_(4,1),TC1_(4,2),TC1_(4,3),x_1-x_0};
P2_1=-x_0+2*x_1-x_2;
P2_2=-x_0+x_1-x_2+x_3;
P2_3=-2*x_0+3*x_1-x_2;
for k from 1 to 3 do (
P2E_k=delete(P2_k,normal);
TE_k=(TEt)*(P2_k);
for i from 0 to (#P2E_k - 1) do (
TP2_(k,i)=(P2E_k)_i - P2_k;
TE_k=TE_k * TP2_(k,i); );
de01A2_k=TE_k;
nu01A2_k=4*x_0+x_2-x_1+P2_k; )

-- the quotient nu01A2_k^7 / de01A2_k 
-- represents contribution of the three isolated fixed points over
-- (x_0^2,x_0x_1,x_0^3,(x_2/x_1)),  1<=k<=3.
\end{verbatim}

\subsection{Contributions on Table \ref{x02}}

\begin{verbatim}
pesosB={x_2-x_0,TB_(2,1),TB_(2,2)};
pesosA={TA_(2,0,0),TA_(2,0,1),TA_(2,0,2),TA_(2,0,3)};
pesos=flatten(append(pesosA,pesosB));
normal=delete(TB_(2,1),pesos);
normal=delete(TB_(2,2),normal);
TEt=(x_1-x_2)*(2*x_1-x_0-x_2);
P2_1=-x_0+x_1;
P2_2=-x_0+x_3;
P2_3=-2*x_0+x_1+x_2;
for k from 1 to 3 do (
P2E_k=delete(P2_k,normal);
TE_k=(TEt)*(P2_k);
for i from 0 to (#P2E_k - 1) do (
TP2_(k,i)=(P2E_k)_i - P2_k;
TE_k=TE_k * TP2_(k,i); );
de02A2_k=TE_k;
nu02A2_k=4*x_0+P2_k; )
-- the quotient nu02A2_k^7 / de02A2_k 
-- represents the contribution of the three isolated
-- fixed points over (x_0^2,x_0x_2,x_0^3), 1<=k<=3.
\end{verbatim}

\subsection{Contributions on Table \ref{x12}}

\begin{verbatim}
pesosB={-TB_(3,0),TB_(3,1),TB_(3,2)};
pesosA={TA_(3,0,0),TA_(3,0,1),TA_(3,0,2),TA_(3,0,3)};
pesos=flatten(append(pesosA,pesosB));
normal=delete(TB_(3,1),pesos);
normal=delete(TB_(3,2),normal);
TEt=(x_0-x_1)*(-2*x_1+x_0+x_2);
for k from 1 to 5 do (
P2_k=normal_(k-1);
P2E_k=delete(P2_k,normal);
TE_k=(TEt)*(P2_k);
for i from 0 to (#P2E_k - 1) do (
TP2_(k,i)=(P2E_k)_i - P2_k;
TE_k=TE_k * TP2_(k,i); );
de03A2_k=TE_k;
nu03A2_k=4*x_0+P2_k; )

-- the quotient nu03A2_k^7 / de03A2_k 
-- represents the contribution of the five isolated
-- fixed points over (x_0^2,x_1^2,x_0^3), 1<=k<=5.
\end{verbatim}

\subsection{Contributions on Table \ref{x13x12s2}}

\begin{verbatim}
pesosB={TB_(3,0),TB_(3,1),TB_(3,2)};
pesosA={TA_(3,4,0),TA_(3,4,1),TA_(3,4,2),TA_(3,4,3)};
pesos=flatten(append(pesosA,pesosB));
normal=delete(TB_(3,1),pesos);
pesos2=delete(TB_(3,2),normal);
tgreta={};
for i from 0 to 4 do (
nm2_i=pesos2_i-TB_(3,2);
tgreta=append(tgreta,nm2_i); )
tgreta=unique tgreta;
normal2=append(tgreta,TB_(3,2));
TEt=(x_0-x_1)*(x_0-x_2);
for k from 1 to 5 do (
P2_k=normal2_(k-1);
P2E_k=delete(P2_k,normal2);
TE_k=(TEt)*(P2_k);
for i from 0 to (#P2E_k - 1) do (
TP2_(k,i)=(P2E_k)_i - P2_k;
TE_k=TE_k * TP2_(k,i); );
deA3_k=TE_k;
nuA3_k=5*x_1+TB_(3,2)+P2_k-x_0; )

-- the quotient nuA3_k^7 / deA3_k 
-- represents the contribution of the five fixed points over
--  (x_1^2,x_1^2,x_1^3,x_0x_2/x_1^2), 1<=k<=5.
\end{verbatim}

\subsection{Contribution on Table \ref{x13x12pl}}

\begin{verbatim}
normal={2*x_0-2*x_1,x_0-x_1,x_2-x_1,x_3-x_1};
TEt=(x_0-x_1)*(x_2-x_0);



for k from 1 to 4 do (
P2_k=normal_(k-1);
P2E_k=delete(P2_k,normal);
TE_k=-1*(TEt)*(P2_k)^2;
for i from 0 to (#P2E_k - 1) do (
TP2_(k,i)=(P2E_k)_i - P2_k;
TE_k=TE_k * TP2_(k,i); );
deA3sp_k=TE_k;
nuA3sp_k=5*x_1+P2_k-x_0+2*x_0-2*x_1; )

-- the quotient nuA3sp_k^7 / deA3sp_k 
-- represents the contribution of the four fixed points over
-- (x_1^2,x_1^2,x_1^3,x_0^2/x_1^2), 1<=k<=4.
\end{verbatim}

\subsection{Contribution on Table \ref{x01s5pl3}}

\begin{verbatim}
normal={x_3-x_2,x_1-x_2,2*x_1-x_0-x_2};
TEt=(2*x_1-x_0-x_2)*(x_1-x_0);
for k from 1 to 3 do (
P2_k=normal_(k-1);
P2E_k=delete(P2_k,normal);
P2E_k=append(P2E_k,x_2-x_1);
TE_k=-1*(TEt)*(P2_k)^2;
for i from 0 to (#P2E_k - 1) do (
TP2_(k,i)=(P2E_k)_i - P2_k;
TE_k=TE_k * TP2_(k,i); );
de01A3sp_k=TE_k;
nu01A3sp_k=4*x_0+x_2-x_1+x_1-x_0+P2_k; )

-- the quotient nu01A3sp_k^7 / de01A3sp_k 
-- represents the contribution of the three fixed points over
-- (x_0^2,x_0x_1,x_0^3,(x_2/x_1),x_1/x_0), 1<=k<=3.
\end{verbatim}

\subsection{Contribution on Table \ref{x01s5pl3v0}}

\begin{verbatim}
normal11={x_1-x_0,x_1-x_2,2*x_1-2*x_2,x_1+x_3-2*x_2};
normal22={3*x_1-x_0-2*x_2,2*x_1-x_0-x_2};
normal=flatten(append(normal11,normal22));
TEt=(x_2-x_1);

for k from 1 to 6 do (
P2_k=normal_(k-1);
P2E_k=delete(P2_k,normal);
TE_k=(TEt)*(P2_k);
for i from 0 to (#P2E_k - 1) do (
TP2_(k,i)=(P2E_k)_i - P2_k;
TE_k=TE_k * TP2_(k,i); );
dev0A4_k=TE_k;
nuv0A4_k=4*x_0+x_2-x_1+x_1-x_0+x_2-x_1+P2_k; )

-- the quotient nuv0A4_k^7 / dev0A4_k 
-- represents the contribution of the six fixed points over
-- (x_0^2,x_0x_1,x_0^3,(x_2/x_1),(x_1/x_0),(x_2/x_1)), 1<=k<=6.
\end{verbatim}

\subsection{Contribution on Table \ref{x01pl2}}

\begin{verbatim}
normal={x_3-x_1,x_2-x_1};
TEt=(-2*x_1+x_0+x_2)*(x_1-x_0);
for k from 1 to 2 do (
P2_k=normal_(k-1);
P2E_k=delete(P2_k,normal);
P2E_k=append(P2E_k,x_1-x_0);
P2E_k=append(P2E_k,x_1-x_0);
TE_k=-1*(TEt)*(P2_k)^2;
for i from 0 to (#P2E_k - 1) do (
TP2_(k,i)=(P2E_k)_i - P2_k;
TE_k=TE_k * TP2_(k,i); );
de01plA2_k=TE_k;
nu01plA2_k=4*x_0+x_1-x_0+P2_k; )

-- the quotient nu01plA2_k^7 / de01plA2_k 
-- represents the contribution of the two fixed points over
-- (x_0^2,x_0x_1,x_0^3,(x_1/x_0)), 1<=k<=2.
\end{verbatim}

\subsection{Contribution on Table \ref{x01pl2pl4}}

\begin{verbatim}
normal={x_0-x_1,x_0+x_3-2*x_1,x_0+x_2-2*x_1,x_1-x_0};
TEt=(-2*x_1+x_0+x_2)*(x_1-x_0);


for k from 1 to 4 do (
P2_k=normal_(k-1);
P2E_k=delete(P2_k,normal);
TE_k=-1*(TEt)*(P2_k)^2;
for i from 0 to (#P2E_k - 1) do (
TP2_(k,i)=(P2E_k)_i - P2_k;
TE_k=TE_k * TP2_(k,i); );
de01qlA3_k=TE_k;
nu01qlA3_k=4*x_0+2*x_1-2*x_0+P2_k; )

-- the quotient nu01qlA3_k^7 / de01qlA3_k 
-- represents the contribution of the four fixed points over
-- (x_0^2,x_0x_1,x_0^3,(x_1/x_0),(x_1/x_0)), 1<=k<=4.
\end{verbatim}

\subsection{Contribution on Table \ref{x02pl5}}

\begin{verbatim}
normal={2*x_1-x_0-x_2,x_1-x_0,x_3-x_2,x_1-x_2,x_2-x_0};
TEt=(x_1-x_2);
for k from 1 to 5 do (
P2_k=normal_(k-1);
P2E_k=delete(P2_k,normal);
TE_k=-1*(TEt)*(P2_k)^2;
for i from 0 to (#P2E_k - 1) do (
TP2_(k,i)=(P2E_k)_i - P2_k;
TE_k=TE_k * TP2_(k,i); );
de02plnA4_k=TE_k;
nu02plnA4_k=4*x_0+x_2-x_0+P2_k; )

-- the quotient nu02plnA4_k^7 / de02plnA4_k 
-- represents the contribution of the five fixed points over
-- (x_0^2,x_0x_2,x_0^3,p_ln=(x_2/x_0)), 1<=k<=5.
\end{verbatim}

\subsection{Contribution on fixed lines}

\begin{verbatim}
-- function to compute the contribution of a fixed P1
Rh=R[h]/h^2;
contl = (Nl,Nll,peso) -> (
denl=product(Nl+h*Nll);
conjl=denl-2*(denl%h);
denl=denl*conjl;
denl=sub(denl,R);
numl=peso^7;
numl=numl*conjl; 
numl=numl//h; 
numl=sub(numl,R);
)
\end{verbatim}

\subsubsection{$\tilde\ell_1$: Example \ref{calcl1}}\label{l1m2}
\begin{verbatim}
Nl={x_0-x_1,2*x_0-2*x_1,x_2-x_0,x_0-x_1,x_2-x_1,x_3-x_1};
Nll=apply(6,j->d_(j+1));
peso=3*x_1+x_0;
contl(Nl,Nll,peso);
bp_1=numl;
cp_1=denl;

-- the quotient bp1 / cp1 
-- represents the contribution of the fixed line ~l_1
\end{verbatim}

\subsubsection{$\tilde\ell_3$: Example \ref{exl2}}\label{l2m2}
\begin{verbatim}
Nl={x_2-x_1,2*x_1-x_0-x_2,x_1-x_0,x_1-x_2,x_3-x_2,2*x_1-x_0-x_2};
Nll=apply(6,j->d_(j+7));
peso=3*x_0+x_2;
contl(Nl,Nll,peso);
bp_2=numl;
cp_2=denl;

-- the quotient bp2 / cp2 
-- represents the contribution of the line ~l_3
\end{verbatim}

\subsubsection{$\tilde\ell_2$: Example \ref{exl3}}\label{l3m2}
\begin{verbatim}
Nl={x_1-x_0,x_1-x_0,x_0+x_2-2*x_1,x_2-x_1,x_3-x_1,x_1-x_0};
Nll=apply(6,j->d_(j+13));
peso=3*x_0+x_1;
contl(Nl,Nll,peso);
bp_3=numl;
cp_3=denl;

-- the quotient bp3 / cp3 
-- represents the contribution of the line ~l_2.
\end{verbatim}

\subsubsection{$\tilde\ell_4$: Example \ref{exl4}}\label{l4m2}
\begin{verbatim}
Nl={x_1-x_0,x_0+x_2-2*x_1,x_1-x_0,x_0+x_2-2*x_1,x_0+x_3-2*x_1,x_0-x_1};
Nll=apply(6,j->d_(j+19));
peso=2*x_1+2*x_0;
contl(Nl,Nll,peso);
bp_4=numl;
cp_4=denl;

-- the quotient bp4 / cp4 
-- represents the contribution of the fixed line ~l_4
\end{verbatim}

\subsubsection{$\tilde\ell_5$: Example \ref{exl5}}\label{l5m2}
\begin{verbatim}
Nl={x_1-x_2,2*x_1-x_0-x_2,x_2-x_0,x_1-x_2,x_3-x_2,x_1-x_0};
Nll=apply(6,j->d_(j+25));
peso=3*x_0+x_2;
contl(Nl,Nll,peso);
bp_5=numl;
cp_5=denl;

-- the quotient bp5 / cp5 
-- represents the contribution of the fixed line ~l_5
\end{verbatim}

\subsection{Sum of contributions}\label{somad's}
\begin{verbatim}
Ry=R[y_0..y_3];

-- take the weights
-- y_0 => 0
-- y_1 => 1
-- y_2 => 5
-- y_3 => 25

pesos={y_0=>0,y_1=>1,y_2=>5,y_3=>25};
Y=ideal(y_0..y_3);

-- a trick to vary the flag

for a from 0 to 3 do (  
bandeira_a=append({},x_0=>y_a);
Ys_a=sub(Y,Ry/y_a);
Ys_a=trim(sub(Ys_a,Ry));
for b from 0 to 2 do (
bandeira_(a,b)=append(bandeira_a,x_1=>(Ys_a)_(2-b));
Ys_(a,b)=sub(Ys_a,Ry/((Ys_a)_(2-b)));
Ys_(a,b)=trim(sub(Ys_(a,b),Ry));
for c from 0 to 1 do (
bandeira_(a,b,c)=append(bandeira_(a,b),x_2=>(Ys_(a,b))_c);
Ys_(a,b,c)=sub(Ys_(a,b),Ry/(Ys_(a,b))_c);
Ys_(a,b,c)=trim(sub(Ys_(a,b,c),Ry));
bandeira_(a,b,c)=append(bandeira_(a,b,c),x_3=>(Ys_(a,b,c))_0);
); ); )

-- in the ideal II we will store all equalities
-- somad_(a,b,c) = somad_(a',b',c'),
-- where somad_(a,b,c) represents the degree of the
-- fiber over the flag (a,b,c), see bandeira_(a,b,c)
II=ideal();

for a from 0 to 3 do (
for b from 0 to 2 do (
for c from 0 to 1 do (
soma=0;
pntsfixos=0;
P1s=0;
for k from 1 to 3 do (
for i from 0 to 4 do (
nun_(k,i)=sub(nu_(k,i),bandeira_(a,b,c));
nun_(k,i)=sub(nun_(k,i),pesos);
nun_(k,i)=sub(nun_(k,i),QQ);
den_(k,i)=sub(de_(k,i),bandeira_(a,b,c));
den_(k,i)=sub(den_(k,i),pesos);
den_(k,i)=sub(den_(k,i),QQ);
ctr_(k,i)=nun_(k,i)^7 / den_(k,i);
soma=soma+ctr_(k,i); 
pntsfixos=pntsfixos + 1; ); );
soma=soma-ctr_(3,4);
pntsfixos=pntsfixos-1;
for k from 1 to 3 do (
for i from 1 to 4 do (
nu0n_(k,i)=sub(nu0_(k,i),bandeira_(a,b,c));
de0n_(k,i)=sub(de0_(k,i),bandeira_(a,b,c));
nu0n_(k,i)=sub(nu0n_(k,i),pesos);
nu0n_(k,i)=sub(nu0n_(k,i),QQ);
de0n_(k,i)=sub(de0n_(k,i),pesos);
de0n_(k,i)=sub(de0n_(k,i),QQ);
ctr0_(k,i)=nu0n_(k,i)^7 / de0n_(k,i);
soma=soma+ctr0_(k,i); 
pntsfixos=pntsfixos + 1;
); );
for k from 1 to 3 do (
nun_(3,4,k)=sub(nu_(3,4,k),bandeira_(a,b,c));  
den_(3,4,k)=sub(de_(3,4,k),bandeira_(a,b,c));
nun_(3,4,k)=sub(nun_(3,4,k),pesos);
nun_(3,4,k)=sub(nun_(3,4,k),QQ);
den_(3,4,k)=sub(den_(3,4,k),pesos);
den_(3,4,k)=sub(den_(3,4,k),QQ); 
ctr_(3,4,k)=nun_(3,4,k)^7 / den_(3,4,k);
soma=soma+ctr_(3,4,k);
pntsfixos=pntsfixos + 1;
);
for k from 1 to 3 do (
nun_(1,0,k)=sub(nu_(1,0,k),bandeira_(a,b,c));  
den_(1,0,k)=sub(de_(1,0,k),bandeira_(a,b,c));
nun_(1,0,k)=sub(nun_(1,0,k),pesos);
nun_(1,0,k)=sub(nun_(1,0,k),QQ);
den_(1,0,k)=sub(den_(1,0,k),pesos);
den_(1,0,k)=sub(den_(1,0,k),QQ);
ctr_(1,0,k)=nun_(1,0,k)^7 / den_(1,0,k);
soma=soma+ctr_(1,0,k);
pntsfixos=pntsfixos + 1;
);
for k from 1 to 3 do (
nun01A2_k=sub(nu01A2_k,bandeira_(a,b,c));  
den01A2_k=sub(de01A2_k,bandeira_(a,b,c));
nun01A2_k=sub(nun01A2_k,pesos);
nun01A2_k=sub(nun01A2_k,QQ);
den01A2_k=sub(den01A2_k,pesos);
den01A2_k=sub(den01A2_k,QQ); 
ctr01A2_k=nun01A2_k^7 / den01A2_k;
soma=soma+ctr01A2_k;  
pntsfixos=pntsfixos + 1;
);
for k from 1 to 3 do (
nun02A2_k=sub(nu02A2_k,bandeira_(a,b,c));  
den02A2_k=sub(de02A2_k,bandeira_(a,b,c));
nun02A2_k=sub(nun02A2_k,pesos);
nun02A2_k=sub(nun02A2_k,QQ);
den02A2_k=sub(den02A2_k,pesos);
den02A2_k=sub(den02A2_k,QQ); 
ctr02A2_k=nun02A2_k^7 / den02A2_k;
soma=soma+ctr02A2_k;  
pntsfixos=pntsfixos + 1;
);
for k from 1 to 5 do (
nun03A2_k=sub(nu03A2_k,bandeira_(a,b,c));  
den03A2_k=sub(de03A2_k,bandeira_(a,b,c));  
nun03A2_k=sub(nun03A2_k,pesos);
nun03A2_k=sub(nun03A2_k,QQ);
den03A2_k=sub(den03A2_k,pesos);
den03A2_k=sub(den03A2_k,QQ); 
ctr03A2_k=nun03A2_k^7 / den03A2_k;
soma=soma+ctr03A2_k;  
pntsfixos=pntsfixos + 1;
);
for k from 1 to 5 do (
nunA3_k=sub(nuA3_k,bandeira_(a,b,c));  
denA3_k=sub(deA3_k,bandeira_(a,b,c));
nunA3_k=sub(nunA3_k,pesos);
nunA3_k=sub(nunA3_k,QQ);
denA3_k=sub(denA3_k,pesos);
denA3_k=sub(denA3_k,QQ); 
ctrA3_k=nunA3_k^7 / denA3_k;
soma=soma+ctrA3_k;  
pntsfixos=pntsfixos + 1;
);
for k from 1 to 4 do (
nunA3sp_k=sub(nuA3sp_k,bandeira_(a,b,c));  
denA3sp_k=sub(deA3sp_k,bandeira_(a,b,c));
nunA3sp_k=sub(nunA3sp_k,pesos);
nunA3sp_k=sub(nunA3sp_k,QQ);
denA3sp_k=sub(denA3sp_k,pesos);
denA3sp_k=sub(denA3sp_k,QQ); 
ctrA3sp_k=nunA3sp_k^7 / denA3sp_k;
soma=soma+ctrA3sp_k;  
pntsfixos=pntsfixos + 1;
);
for k from 1 to 3 do (
nun01A3sp_k=sub(nu01A3sp_k,bandeira_(a,b,c));  
den01A3sp_k=sub(de01A3sp_k,bandeira_(a,b,c));
nun01A3sp_k=sub(nun01A3sp_k,pesos);
nun01A3sp_k=sub(nun01A3sp_k,QQ);
den01A3sp_k=sub(den01A3sp_k,pesos);
den01A3sp_k=sub(den01A3sp_k,QQ); 
ctr01A3sp_k=nun01A3sp_k^7 / den01A3sp_k;
soma=soma+ctr01A3sp_k;  
pntsfixos=pntsfixos + 1;
);
for k from 1 to 6 do (
nunv0A4_k=sub(nuv0A4_k,bandeira_(a,b,c));  
denv0A4_k=sub(dev0A4_k,bandeira_(a,b,c));
nunv0A4_k=sub(nunv0A4_k,pesos);
nunv0A4_k=sub(nunv0A4_k,QQ);
denv0A4_k=sub(denv0A4_k,pesos);
denv0A4_k=sub(denv0A4_k,QQ); 
ctrv0A4_k=nunv0A4_k^7 / denv0A4_k;
soma=soma+ctrv0A4_k;  
pntsfixos=pntsfixos + 1;
);
for k from 1 to 2 do (
nun01plA2_k=sub(nu01plA2_k,bandeira_(a,b,c));  
den01plA2_k=sub(de01plA2_k,bandeira_(a,b,c));
nun01plA2_k=sub(nun01plA2_k,pesos);
nun01plA2_k=sub(nun01plA2_k,QQ);
den01plA2_k=sub(den01plA2_k,pesos);
den01plA2_k=sub(den01plA2_k,QQ); 
ctr01plA2_k=nun01plA2_k^7 / den01plA2_k;
soma=soma+ctr01plA2_k;  
pntsfixos=pntsfixos + 1;
);
for k from 1 to 4 do (
nun01qlA3_k=sub(nu01qlA3_k,bandeira_(a,b,c));  
den01qlA3_k=sub(de01qlA3_k,bandeira_(a,b,c));
nun01qlA3_k=sub(nun01qlA3_k,pesos);
nun01qlA3_k=sub(nun01qlA3_k,QQ);
den01qlA3_k=sub(den01qlA3_k,pesos);
den01qlA3_k=sub(den01qlA3_k,QQ); 
ctr01qlA3_k=nun01qlA3_k^7 / den01qlA3_k;
soma=soma+ctr01qlA3_k;  
pntsfixos=pntsfixos + 1;
);
for k from 1 to 5 do (
nun02plnA4_k=sub(nu02plnA4_k,bandeira_(a,b,c));  
den02plnA4_k=sub(de02plnA4_k,bandeira_(a,b,c));
nun02plnA4_k=sub(nun02plnA4_k,pesos);
nun02plnA4_k=sub(nun02plnA4_k,QQ);
den02plnA4_k=sub(den02plnA4_k,pesos);
den02plnA4_k=sub(den02plnA4_k,QQ); 
ctr02plnA4_k=nun02plnA4_k^7 / den02plnA4_k;
soma=soma+ctr02plnA4_k;  
pntsfixos=pntsfixos + 1;
);
for k from 1 to 5 do ( 
cpn_k=sub(cp_k,bandeira_(a,b,c));  
bpn_k=sub(bp_k,bandeira_(a,b,c));
cpn_k=sub(cpn_k,pesos);
cpn_k=sub(cpn_k,QQ);
bpn_k=sub(bpn_k,pesos);
ctrP1_k=bpn_k / cpn_k;
soma=soma + ctrP1_k;
P1s=P1s+1;
); 
somad_(a,b,c)=soma;
II=II+(somad_(a,b,c)-somad_(0,0,0));
); ); )
\end{verbatim}

\subsection{Degree to Proposition \ref{degreefibre}}\label{degfibercalc}

\begin{verbatim}
-- The relations for d's
II=trim(sub(II,Ra))

-- and the degree of a fiber 
soma=sub(-somad_(0,0,1),Ra/II)

\end{verbatim}

The minus sign is due to the minus in
$$\int\limits_{{\mathbb{Y}_4}_{(p,\ell,v)}}-c_1^{7}(\mathcal{W})
\cap [{\mathbb{Y}_4}_{(p,\ell,v)}].
$$
The variable soma returns the result, \,\mb{21},\, of  
Proposition \ref{degreefibre}.

\section{Degree of the Exceptional Component}\label{degexc}

The scripts are similar to the presented in the previous section.

\begin{verbatim}
R=QQ[d_1..d_30];
Rx=R[x_0..x_3];
B={2*x_0,x_0+x_1,x_0+x_2,2*x_1};
A={3*x_0,2*x_0+x_1,2*x_0+x_2,2*x_0+x_3};
A0_1=append(A,(x_0+2*x_1));
A0_2=append(A,(x_0+x_1+x_2));
A0_3=append(A,(3*x_1));
vtgflag={x_1-x_0,x_2-x_0,x_3-x_0,x_2-x_1,x_3-x_1,x_3-x_2};
tgflag=product vtgflag;

-- Y ---> X

for k from 1 to 3 do (
B0_k=delete(B_k,B);
TB0_k=1;
for i from 0 to 2 do (
TB_(k,i)=(B0_k)_i-B_k;
TB0_k=TB0_k * TB_(k,i); ); )

for k from 1 to 3 do (
for i from 0 to 4 do (
TA0_(k,i)=1;
FA0_(k,i)=delete((A0_k)_i,A0_k);
for j from 0 to (#(FA0_(k,i)) - 1) do (
TA_(k,i,j)=(FA0_(k,i))_j - (A0_k)_i;
TA0_(k,i)=TA0_(k,i) * TA_(k,i,j); ); ); )

for k from 1 to 3 do (
for i from 0 to 4 do ( 
de_(k,i)=TA0_(k,i) * TB0_k * tgflag; 
nu_(k,i)=B_k+(A0_k)_i-x_0;    ); )

-- Y -----> X'
B00r=delete(B_0,B);
for k from 1 to 3 do (
B00_k=delete(B_k,B00r);
TB00_k=(B_k-B_0)*((B00_k)_0-B_k)*((B00_k)_1-B_k);
for i from 0 to 4 do (
de0_(k,i)=(TA0_(k,i)*TB00_k)*(tgflag); 
nu0_(k,i)=B_0+(A0_k)_i-x_0;     ); )

-- point (3,4)=x_1^2,x_1^3
pesosB={TB_(3,0),TB_(3,1),TB_(3,2)};
pesosA={TA_(3,4,0),TA_(3,4,1),TA_(3,4,2),TA_(3,4,3)};
pesos=flatten(append(pesosA,pesosB));
normal=delete(x_0-x_1,pesos);
curva_1=3*x_0-3*x_1;
curva_2=2*x_0+x_2-3*x_1;
curva_3=2*x_0+x_3-3*x_1;
for k from 1 to 3 do (
Cv_k=delete(curva_k,normal);
TC_(3,4,k)=(x_0-x_1)*(curva_k);
for i from 0 to (#Cv_k - 1) do (
TC1_(k,i)=(Cv_k)_i - curva_k;
TC_(3,4,k)=TC_(3,4,k) * TC1_(k,i); );
de_(3,4,k)=(TC_(3,4,k))*(tgflag);
nu_(3,4,k)=5*x_1 - x_0 + curva_k; )

-- point 0_(1,0)=x_0^2,x_0x_1,x_0^3
normalB={(B00_1)_0-B_1 , (B00_1)_1-B_1};
normalA={TA_(1,0,0),TA_(1,0,1),TA_(1,0,2),TA_(1,0,3)};
normal=flatten(append(normalA,normalB));

curva_1=-x_0+x_2;
curva_2=-x_0+x_3;
curva_3=-2*x_0+2*x_1;
curva_4=-x_1+x_2;

for k from 1 to 4 do (
Cv_k=delete(curva_k,normal);
TC_(1,0,k)=(x_1-x_0)*(curva_k);
for i from 0 to (#Cv_k - 1) do (
TC1_(k,i)=(Cv_k)_i - curva_k;
TC_(1,0,k)=TC_(1,0,k) * TC1_(k,i); );
de_(1,0,k)=(TC_(1,0,k))*(tgflag);
nu_(1,0,k)=4*x_0 + curva_k; )

-- point 0_(1,0),(x_2/x_1).
TEt=(x_2-x_1)*(2*x_1-x_0-x_2);
normal={TC1_(4,0),TC1_(4,1),TC1_(4,2),TC1_(4,3),x_1-x_0};
P2_1=-x_0+2*x_1-x_2;
P2_2=-x_0+x_1-x_2+x_3;
P2_3=-2*x_0+3*x_1-x_2;
for k from 1 to 3 do (
P2E_k=delete(P2_k,normal);
TE_k=(TEt)*(P2_k);
for i from 0 to (#P2E_k - 1) do (
TP2_(k,i)=(P2E_k)_i - P2_k;
TE_k=TE_k * TP2_(k,i); );
de01A2_k=(TE_k)*(tgflag);
nu01A2_k=4*x_0+x_2-x_1+P2_k; )

-- point 0_(2,0)=x_0^2,x_0x_2,x_0^3
pesosB={x_2-x_0,TB_(2,1),TB_(2,2)};
pesosA={TA_(2,0,0),TA_(2,0,1),TA_(2,0,2),TA_(2,0,3)};
pesos=flatten(append(pesosA,pesosB));
normal=delete(TB_(2,1),pesos);
normal=delete(TB_(2,2),normal);
TEt=(x_1-x_2)*(2*x_1-x_0-x_2);
P2_1=-x_0+x_1;
P2_2=-x_0+x_3;
P2_3=-2*x_0+x_1+x_2;
for k from 1 to 3 do (
P2E_k=delete(P2_k,normal);
TE_k=(TEt)*(P2_k);
for i from 0 to (#P2E_k - 1) do (
TP2_(k,i)=(P2E_k)_i - P2_k;
TE_k=TE_k * TP2_(k,i); );
de02A2_k=(TE_k)*(tgflag);
nu02A2_k=4*x_0+P2_k; )

-- point 0_(3,0)=x_0^2,x_1^2,x_0^3
pesosB={-TB_(3,0),TB_(3,1),TB_(3,2)};
pesosA={TA_(3,0,0),TA_(3,0,1),TA_(3,0,2),TA_(3,0,3)};
pesos=flatten(append(pesosA,pesosB));
normal=delete(TB_(3,1),pesos);
normal=delete(TB_(3,2),normal);
TEt=(x_0-x_1)*(-2*x_1+x_0+x_2);
for k from 1 to 5 do (
P2_k=normal_(k-1);
P2E_k=delete(P2_k,normal);
TE_k=(TEt)*(P2_k);
for i from 0 to (#P2E_k - 1) do (
TP2_(k,i)=(P2E_k)_i - P2_k;
TE_k=TE_k * TP2_(k,i); );
de03A2_k=(TE_k)*(tgflag);
nu03A2_k=4*x_0+P2_k; )

-- over (x_1^3,x_1^2, x_0x_2/x_1^2).
pesosB={TB_(3,0),TB_(3,1),TB_(3,2)};
pesosA={TA_(3,4,0),TA_(3,4,1),TA_(3,4,2),TA_(3,4,3)};
pesos=flatten(append(pesosA,pesosB));
normal=delete(TB_(3,1),pesos);
pesos2=delete(TB_(3,2),normal);
tgreta={};
for i from 0 to 4 do (
nm2_i=pesos2_i-TB_(3,2);
tgreta=append(tgreta,nm2_i); )
tgreta=unique tgreta;
normal2=append(tgreta,TB_(3,2));
TEt=(x_0-x_1)*(x_0-x_2);
for k from 1 to 5 do (
P2_k=normal2_(k-1);
P2E_k=delete(P2_k,normal2);
TE_k=(TEt)*(P2_k);
for i from 0 to (#P2E_k - 1) do (
TP2_(k,i)=(P2E_k)_i - P2_k;
TE_k=TE_k * TP2_(k,i); );
deA3_k=(TE_k)*(tgflag);
nuA3_k=5*x_1+TB_(3,2)+P2_k-x_0; )

-- ------------------------------------------------
-- ------ contributions of the isolated -----------
-- ------- fixed points over fixed P1s ------------
-- ------------------------------------------------

-- over (x_1^3,x_1^2, x_0^2/x_1^2).
normal={2*x_0-2*x_1,x_0-x_1,x_2-x_1,x_3-x_1};
TEt=(x_0-x_1)*(x_2-x_0);
for k from 1 to 4 do (
P2_k=normal_(k-1);
P2E_k=delete(P2_k,normal);
TE_k=-1*(TEt)*(P2_k)^2;
for i from 0 to (#P2E_k - 1) do (
TP2_(k,i)=(P2E_k)_i - P2_k;
TE_k=TE_k * TP2_(k,i); );
deA3sp_k=(TE_k)*(tgflag);
nuA3sp_k=5*x_1+P2_k-x_0+2*x_0-2*x_1; )

-- point 0_(1,0),(x_2/x_1), (x_1/x_0).
-- over (x_0^2,x_0x_1,x_0^3,s_5',t_1=t_3)
normal={x_3-x_2,x_1-x_2,2*x_1-x_0-x_2};
TEt=(2*x_1-x_0-x_2)*(x_1-x_0);

for k from 1 to 3 do (
P2_k=normal_(k-1);
P2E_k=delete(P2_k,normal);
P2E_k=append(P2E_k,x_2-x_1);
TE_k=-1*(TEt)*(P2_k)^2;
for i from 0 to (#P2E_k - 1) do (
TP2_(k,i)=(P2E_k)_i - P2_k;
TE_k=TE_k * TP2_(k,i); );
de01A3sp_k=(TE_k)*(tgflag);
nu01A3sp_k=4*x_0+x_2-x_1+x_1-x_0+P2_k; )

-- point 0_(1,0),(x_2/x_1), (x_1/x_0), (x_2/x_1).
-- over (x_0^2,x_0x_1,x_0^3,s_5',(t_1=t_3),v_0')
normal={x_1-x_0,x_1-x_2,2*x_1-2*x_2,x_1+x_3-2*x_2,3*x_1-x_0-2*x_2,
2*x_1-x_0-x_2};  TEt=(x_2-x_1);
for k from 1 to 6 do (
P2_k=normal_(k-1);
P2E_k=delete(P2_k,normal);
TE_k=(TEt)*(P2_k);
for i from 0 to (#P2E_k - 1) do (
TP2_(k,i)=(P2E_k)_i - P2_k;
TE_k=TE_k * TP2_(k,i); );
dev0A4_k=(TE_k)*(tgflag);
nuv0A4_k=4*x_0+x_2-x_1+x_1-x_0+x_2-x_1+P2_k; )

-- point 0_(1,0),(x_1/x_0) over (x_0^2,x_0x_1,x_0^3,p_l)
normal={x_3-x_1,x_2-x_1};
TEt=(-2*x_1+x_0+x_2)*(x_1-x_0);
for k from 1 to 2 do (
P2_k=normal_(k-1);
P2E_k=delete(P2_k,normal);
P2E_k=append(P2E_k,x_1-x_0);
P2E_k=append(P2E_k,x_1-x_0);
TE_k=-1*(TEt)*(P2_k)^2;
for i from 0 to (#P2E_k - 1) do (
TP2_(k,i)=(P2E_k)_i - P2_k;
TE_k=TE_k * TP2_(k,i); );
de01plA2_k=(TE_k)*(tgflag);
nu01plA2_k=4*x_0+x_1-x_0+P2_k; )

-- point 0_(1,0),p_l=(x_1/x_0),q_l2=(x_1/x_0)
-- over (x_0^2,x_0x_1,x_0^3,p_l,q_l2)
normal={x_0-x_1,x_0+x_3-2*x_1,x_0+x_2-2*x_1,x_1-x_0};
TEt=(-2*x_1+x_0+x_2)*(x_1-x_0);
for k from 1 to 4 do (
P2_k=normal_(k-1);
P2E_k=delete(P2_k,normal);
TE_k=-1*(TEt)*(P2_k)^2;
for i from 0 to (#P2E_k - 1) do (
TP2_(k,i)=(P2E_k)_i - P2_k;
TE_k=TE_k * TP2_(k,i); );
de01qlA3_k=(TE_k)*(tgflag);
nu01qlA3_k=4*x_0+2*x_1-2*x_0+P2_k; )

-- point 0_(2,0),p_l=(x_2/x_0)
-- over (x_0^2,x_0x_2,x_0^3,p_l)
normal={2*x_1-x_0-x_2,x_1-x_0,x_3-x_2,x_1-x_2,x_2-x_0};
TEt=(x_1-x_2);
for k from 1 to 5 do (
P2_k=normal_(k-1);
P2E_k=delete(P2_k,normal);
TE_k=-1*(TEt)*(P2_k)^2;
for i from 0 to (#P2E_k - 1) do (
TP2_(k,i)=(P2E_k)_i - P2_k;
TE_k=TE_k * TP2_(k,i); );
de02plnA4_k=(TE_k)*(tgflag);
nu02plnA4_k=4*x_0+x_2-x_0+P2_k; )

-- ---------------------------------
-- ------ fixed lines --------------
-- ---------------------------------

-- contributions of fixed P1
Rh=Rx[h]/h^2;
contl = (Nl,Nll,peso) -> (
denl=product(Nl+h*Nll)*tgflag;
conjl=denl-2*(denl%h);
denl=denl*conjl;
denl=sub(denl,Rx);
numl=peso^13;
numl=numl*conjl; 
numl=numl//h; 
numl=sub(numl,Rx);
)

-- fixed P^1 over (x_1^3,x_1^2)
Nl={x_0-x_1,2*x_0-2*x_1,x_2-x_0,x_0-x_1,x_2-x_1,x_3-x_1};
Nll=apply(6,j->d_(j+1));
peso=3*x_1+x_0;
contl(Nl,Nll,peso);
bp_1=numl;
cp_1=denl;

-- fixed P^1 over  (x_0^3,x_0x_1,x_0^2,s_5')
Nl={x_2-x_1,2*x_1-x_0-x_2,x_1-x_0,x_1-x_2,x_3-x_2,2*x_1-x_0-x_2};
Nll=apply(6,j->d_(j+7));
peso=3*x_0+x_2;
contl(Nl,Nll,peso);
bp_2=numl;
cp_2=denl;

-- fixed P^1 over (x_0^3,x_0x_1,x_0^2)
Nl={x_1-x_0,x_1-x_0,x_0+x_2-2*x_1,x_2-x_1,x_3-x_1,x_1-x_0};
Nll=apply(6,j->d_(j+13));
peso=3*x_0+x_1;
contl(Nl,Nll,peso);
bp_3=numl;
cp_3=denl;

-- fixed P^1 over  (x_0^3,x_0x_1,x_0^2,p_l)
Nl={x_1-x_0,x_0+x_2-2*x_1,x_1-x_0,x_0+x_2-2*x_1,x_0+x_3-2*x_1,x_0-x_1};
Nll=apply(6,j->d_(j+19));
peso=2*x_1+2*x_0;
contl(Nl,Nll,peso);
bp_4=numl;
cp_4=denl;

-- fixed P^1 over  (x_0^3,x_0x_2,x_0^2)
Nl={x_1-x_2,2*x_1-x_0-x_2,x_2-x_0,x_1-x_2,x_3-x_2,x_1-x_0};
Nll=apply(6,j->d_(j+25));
peso=3*x_0+x_2;
contl(Nl,Nll,peso);
bp_5=numl;
cp_5=denl;

-- --------------------------------------
-- --------------- sum ------------------
-- --------------------------------------
-- --------------------------------------

-- trick to vary the flag

Ry=R[y_0..y_3];

soma=0;
pesos={y_0=>0,y_1=>1,y_2=>5,y_3=>25};
pntsfixos=0;
P1s=0;

Y=ideal(y_0..y_3);

for i from 0 to 3 do (  
bandeira_i=append({},x_0=>y_i);
Ys_i=sub(Y,Ry/y_i);
Ys_i=trim(sub(Ys_i,Ry));
for j from 0 to 2 do (
bandeira_(i,j)=append(bandeira_i,x_1=>(Ys_i)_(2-j));
Ys_(i,j)=sub(Ys_i,Ry/((Ys_i)_(2-j)));
Ys_(i,j)=trim(sub(Ys_(i,j),Ry));
for k from 0 to 1 do (
bandeira_(i,j,k)=append(bandeira_(i,j),x_2=>(Ys_(i,j))_k);
Ys_(i,j,k)=sub(Ys_(i,j),Ry/(Ys_(i,j))_k);
Ys_(i,j,k)=trim(sub(Ys_(i,j,k),Ry));
bandeira_(i,j,k)=append(bandeira_(i,j,k),x_3=>(Ys_(i,j,k))_0);
); ); )

time
for a from 0 to 3 do (
for b from 0 to 2 do (
for c from 0 to 1 do (
for k from 1 to 3 do (
for i from 0 to 4 do (
nun_(k,i)=sub(nu_(k,i),bandeira_(a,b,c));
nun_(k,i)=sub(nun_(k,i),pesos);
nun_(k,i)=sub(nun_(k,i),QQ);
den_(k,i)=sub(de_(k,i),bandeira_(a,b,c));
den_(k,i)=sub(den_(k,i),pesos);
den_(k,i)=sub(den_(k,i),QQ);
ctr_(k,i)=nun_(k,i)^13 / den_(k,i);
soma=soma+ctr_(k,i); 
pntsfixos=pntsfixos + 1; ); );
soma=soma-ctr_(3,4);
pntsfixos=pntsfixos-1;
for k from 1 to 3 do (
for i from 1 to 4 do (
nu0n_(k,i)=sub(nu0_(k,i),bandeira_(a,b,c));
de0n_(k,i)=sub(de0_(k,i),bandeira_(a,b,c));
nu0n_(k,i)=sub(nu0n_(k,i),pesos);
nu0n_(k,i)=sub(nu0n_(k,i),QQ);
de0n_(k,i)=sub(de0n_(k,i),pesos);
de0n_(k,i)=sub(de0n_(k,i),QQ);
ctr0_(k,i)=nu0n_(k,i)^13 / de0n_(k,i);
soma=soma+ctr0_(k,i); 
pntsfixos=pntsfixos + 1;
); );
for k from 1 to 3 do (
nun_(3,4,k)=sub(nu_(3,4,k),bandeira_(a,b,c));  
den_(3,4,k)=sub(de_(3,4,k),bandeira_(a,b,c));
nun_(3,4,k)=sub(nun_(3,4,k),pesos);
nun_(3,4,k)=sub(nun_(3,4,k),QQ);
den_(3,4,k)=sub(den_(3,4,k),pesos);
den_(3,4,k)=sub(den_(3,4,k),QQ); 
ctr_(3,4,k)=nun_(3,4,k)^13 / den_(3,4,k);
soma=soma+ctr_(3,4,k);
pntsfixos=pntsfixos + 1;
);
for k from 1 to 3 do (
nun_(1,0,k)=sub(nu_(1,0,k),bandeira_(a,b,c));  
den_(1,0,k)=sub(de_(1,0,k),bandeira_(a,b,c));
nun_(1,0,k)=sub(nun_(1,0,k),pesos);
nun_(1,0,k)=sub(nun_(1,0,k),QQ);
den_(1,0,k)=sub(den_(1,0,k),pesos);
den_(1,0,k)=sub(den_(1,0,k),QQ);
ctr_(1,0,k)=nun_(1,0,k)^13 / den_(1,0,k);
soma=soma+ctr_(1,0,k);
pntsfixos=pntsfixos + 1;
);
for k from 1 to 3 do (
nun01A2_k=sub(nu01A2_k,bandeira_(a,b,c));  
den01A2_k=sub(de01A2_k,bandeira_(a,b,c));
nun01A2_k=sub(nun01A2_k,pesos);
nun01A2_k=sub(nun01A2_k,QQ);
den01A2_k=sub(den01A2_k,pesos);
den01A2_k=sub(den01A2_k,QQ); 
ctr01A2_k=nun01A2_k^13 / den01A2_k;
soma=soma+ctr01A2_k;  
pntsfixos=pntsfixos + 1;
);
for k from 1 to 3 do (
nun02A2_k=sub(nu02A2_k,bandeira_(a,b,c));  
den02A2_k=sub(de02A2_k,bandeira_(a,b,c));
nun02A2_k=sub(nun02A2_k,pesos);
nun02A2_k=sub(nun02A2_k,QQ);
den02A2_k=sub(den02A2_k,pesos);
den02A2_k=sub(den02A2_k,QQ); 
ctr02A2_k=nun02A2_k^13 / den02A2_k;
soma=soma+ctr02A2_k;  
pntsfixos=pntsfixos + 1;
);
for k from 1 to 5 do (
nun03A2_k=sub(nu03A2_k,bandeira_(a,b,c));  
den03A2_k=sub(de03A2_k,bandeira_(a,b,c));  
nun03A2_k=sub(nun03A2_k,pesos);
nun03A2_k=sub(nun03A2_k,QQ);
den03A2_k=sub(den03A2_k,pesos);
den03A2_k=sub(den03A2_k,QQ); 
ctr03A2_k=nun03A2_k^13 / den03A2_k;
soma=soma+ctr03A2_k;  
pntsfixos=pntsfixos + 1;
);
for k from 1 to 5 do (
nunA3_k=sub(nuA3_k,bandeira_(a,b,c));  
denA3_k=sub(deA3_k,bandeira_(a,b,c));
nunA3_k=sub(nunA3_k,pesos);
nunA3_k=sub(nunA3_k,QQ);
denA3_k=sub(denA3_k,pesos);
denA3_k=sub(denA3_k,QQ); 
ctrA3_k=nunA3_k^13 / denA3_k;
soma=soma+ctrA3_k;  
pntsfixos=pntsfixos + 1;
);
for k from 1 to 4 do (
nunA3sp_k=sub(nuA3sp_k,bandeira_(a,b,c));  
denA3sp_k=sub(deA3sp_k,bandeira_(a,b,c));
nunA3sp_k=sub(nunA3sp_k,pesos);
nunA3sp_k=sub(nunA3sp_k,QQ);
denA3sp_k=sub(denA3sp_k,pesos);
denA3sp_k=sub(denA3sp_k,QQ); 
ctrA3sp_k=nunA3sp_k^13 / denA3sp_k;
soma=soma+ctrA3sp_k;  
pntsfixos=pntsfixos + 1;
);
for k from 1 to 3 do (
nun01A3sp_k=sub(nu01A3sp_k,bandeira_(a,b,c));  
den01A3sp_k=sub(de01A3sp_k,bandeira_(a,b,c));
nun01A3sp_k=sub(nun01A3sp_k,pesos);
nun01A3sp_k=sub(nun01A3sp_k,QQ);
den01A3sp_k=sub(den01A3sp_k,pesos);
den01A3sp_k=sub(den01A3sp_k,QQ); 
ctr01A3sp_k=nun01A3sp_k^13 / den01A3sp_k;
soma=soma+ctr01A3sp_k;  
pntsfixos=pntsfixos + 1;
);
for k from 1 to 6 do (
nunv0A4_k=sub(nuv0A4_k,bandeira_(a,b,c));  
denv0A4_k=sub(dev0A4_k,bandeira_(a,b,c));
nunv0A4_k=sub(nunv0A4_k,pesos);
nunv0A4_k=sub(nunv0A4_k,QQ);
denv0A4_k=sub(denv0A4_k,pesos);
denv0A4_k=sub(denv0A4_k,QQ); 
ctrv0A4_k=nunv0A4_k^13 / denv0A4_k;
soma=soma+ctrv0A4_k;  
pntsfixos=pntsfixos + 1;
);
for k from 1 to 2 do (
nun01plA2_k=sub(nu01plA2_k,bandeira_(a,b,c));  
den01plA2_k=sub(de01plA2_k,bandeira_(a,b,c));
nun01plA2_k=sub(nun01plA2_k,pesos);
nun01plA2_k=sub(nun01plA2_k,QQ);
den01plA2_k=sub(den01plA2_k,pesos);
den01plA2_k=sub(den01plA2_k,QQ); 
ctr01plA2_k=nun01plA2_k^13 / den01plA2_k;
soma=soma+ctr01plA2_k;  
pntsfixos=pntsfixos + 1;
);
for k from 1 to 4 do (
nun01qlA3_k=sub(nu01qlA3_k,bandeira_(a,b,c));  
den01qlA3_k=sub(de01qlA3_k,bandeira_(a,b,c));
nun01qlA3_k=sub(nun01qlA3_k,pesos);
nun01qlA3_k=sub(nun01qlA3_k,QQ);
den01qlA3_k=sub(den01qlA3_k,pesos);
den01qlA3_k=sub(den01qlA3_k,QQ); 
ctr01qlA3_k=nun01qlA3_k^13 / den01qlA3_k;
soma=soma+ctr01qlA3_k;  
pntsfixos=pntsfixos + 1;
);
for k from 1 to 5 do (
nun02plnA4_k=sub(nu02plnA4_k,bandeira_(a,b,c));  
den02plnA4_k=sub(de02plnA4_k,bandeira_(a,b,c));
nun02plnA4_k=sub(nun02plnA4_k,pesos);
nun02plnA4_k=sub(nun02plnA4_k,QQ);
den02plnA4_k=sub(den02plnA4_k,pesos);
den02plnA4_k=sub(den02plnA4_k,QQ); 
ctr02plnA4_k=nun02plnA4_k^13 / den02plnA4_k;
soma=soma+ctr02plnA4_k;  
pntsfixos=pntsfixos + 1;
);
for k from 1 to 5 do ( 
cpn_k=sub(cp_k,bandeira_(a,b,c));  
bpn_k=sub(bp_k,bandeira_(a,b,c));
cpn_k=sub(cpn_k,pesos);
cpn_k=sub(cpn_k,R);
bpn_k=sub(bpn_k,pesos);
bpn_k=sub(bpn_k,R);
ctrP1_k=bpn_k / cpn_k;
soma=soma + ctrP1_k;
P1s=P1s+1;
);
); ); )
-- soma returns the degree as a linear expression on d's
-- the relations for d's are

II1=ideal(d_29,d_27+2,2*d_25+d_26+2*d_28+2*d_30+1);
II2=ideal(d_23,d_20+d_22,d_19+d_21-d_24+2,d_17,d_16,d_15-1);
II3=ideal(d_13+d_14+d_18+2,d_11,d_9+d_30+1,d_8+d_12);
II4=ideal(2*d_7-2*d_10-d_26-2*d_30+1,d_6,d_5,d_3-1,2*d_1+d_2+2*d_4+2);
II=II1+II2+II3+II4;

-- substitute these relations in soma
-- to obtain the degree of the exceptional component
soma=sub(-soma, R/II)

\end{verbatim}

The minus sign is due to the minus in Theorem \ref{integral}. The variable soma returns the result of the main Theorem \ref{GRAUE3}.

\chapter{Appendix B}

\section{Induced action}\label{acaoquadcub}

For a given action of $\bb{C}^*$ in $\bb{P}^3$,
$$ \begin{array}{cccl}
\xi: & \mathbb{C}^* \times \bb{P}^3 &  \rightarrow	& \bb{P}^3\\
& [t,(q_0:q_1:q_2:q_3)] & \mapsto & (t^{-w_0}q_0:t^{-w_1}q_1:t^{-w_2}q_2:t^{-w_3}q_3)\ ,
\end{array} $$
we have, for each $t \in \mathbb{C}^*$, a linear map $T_t: \bb{C}^4 \rightarrow \bb{C}^4$ defined on canonical basis $\{e_0,e_1,e_2,e_3\}$ by $T_t(e_i)=t^{-w_i}e_i$.

The dual space $(\bb{C}^4)^\vee$ has dual basis $\{x_0,x_1,x_2,x_3\}$, the linear forms with $$x_i(e_j) = \delta_{ij} = \left\{\begin{array}{c}
1, \textrm{ if } i=j \\
0, \textrm{ if } i\not=j  
\end{array}.\right.$$ 

Under $T_t$ we obtain another basis of $\bb{C}^4$, the basis $\beta=\{t^{-w_0}e_0,t^{-w_1}e_1,t^{-w_2}e_2,t^{-w_3}e_3\}$. 
The dual basis of $\beta$ is the basis $\alpha=\{y_0,y_1,y_2,y_3\}$ of $(\bb{C}^4)^\vee$ given by linear forms with $$y_i(t^{-w_j}e_j)=\delta_{ij} \  \therefore \ y_i(e_j)=t^{w_j}\delta_{ij}. $$
Thus, the dual basis is $\alpha=\{t^{w_0}x_0,t^{w_1}x_1,t^{w_2}x_2,t^{w_3}x_3\}$, and so the action $\xi$ induces an action $\xi^\vee$ on $\pd3$ just by
\begin{equation}\label{xidual}
\begin{array}{cccccccc}
x_0 & \mapsto & t^{w_0}x_0 & \ & \ & x_2 & \mapsto & t^{w_2}x_2 \\
x_1 & \mapsto & t^{w_1}x_1 & \ & \ & x_3 & \mapsto & t^{w_3}x_3 \\
\end{array}
\end{equation}

For a point $(1:q_1:q_2:q_3) \in \bb{P}^3$ the action $\xi$ gives
$$t\circ(1:q_1:q_2:q_3) = (t^{-w_0}:t^{-w_1}q_1:t^{-w_2}q_2:t^{-w_3}q_3) = (1:t^{w_0-w_1}q_1:t^{w_0-w_2}q_2:t^{w_0-w_3}q_3).$$

If the weights $w_i, \ 0\leq i \leq 3$, are all distinct, then for $t\not= \pm 1$ the point will be changed. Therefore, the only fixed point is the point (1:0:0:0). In a similar way, the other fixed points by the action $\xi$ are $(0:1:0:0), \ (0:0:1:0)$ and $(0:0:0:1)$. 

The same rule holds for planes: for distinct weights $w_i$ the only fixed planes are the standard ones, $x_i=0$, $0 \leq i \leq 3$, and so there is only four fixed points on $\pd3$ (Note: the planes $x_i=0$ are fixed, but not \textit{point-a-point}).

However, if we set $w_0=w_1$ then all points of the line ${(q_0:q_1:0:0)}$ will be fixed by the action, point-a-point, and so we gain a fixed component in $\mathbb{P}^3$ of positive dimension. 

\subsection{Action on lines}

Let $\bb{G} = \bb{G}(1,3)$ the grassmannian of lines in $\bb{P}^3$, and consider the Pl\"ucker embedding 
$\bb{G} \xrightarrow{\sim} Q \subset \bb{P}^5,$ where $Q = \{(p_{01}:p_{02}:p_{03}:p_{12}:p_{13}:p_{23}) \in \bb{P}^5 | \ p_{01}p_{23}-p_{02}p_{13}+p_{03}p_{12}=0 \}$ (\cite{3264}, example 2.1, pg.68).

The Pl\"ucker embedding associates a line  $$\ell = \left\{\begin{array}{c}
	a_0x_0+a_1x_1+a_2x_2+a_3x_3=0\\
	b_0x_0+b_1x_1+b_2x_2+b_3x_3=0  
\end{array}\right.$$ 
to the point of $Q$ whose coordinates are the $(2\times 2)$--minors of the coefficients matrix
 $$\left[\begin{array}{cccc}
	a_0 & a_1 & a_2 & a_3\\
	b_0 & b_1 & b_2 & b_3  
\end{array}\right].$$ 

The action (\ref{xidual}) applied on $\ell $  gives the line
$$t \circ \ell = \left\{\begin{array}{c}
	a_0t^{w_0}x_0+a_1t^{w_1}x_1+a_2t^{w_2}x_2+a_3t^{w_3}x_3=0\\
	b_0t^{w_0}x_0+b_1t^{w_1}x_1+b_2t^{w_2}x_2+b_3t^{w_3}x_3=0  
\end{array},\right.$$
and the coefficients matrix of $t \circ \ell$ is
$$\left[\begin{array}{cccc}
	t^{w_0}a_0 & t^{w_1}a_1 & t^{w_2}a_2 & t^{w_3}a_3\\
	t^{w_0}b_0 & t^{w_1}b_1 & t^{w_2}b_2 & t^{w_3}b_3  
\end{array}\right].$$ 
Computing the (2$\times$2)--minors of this matrix, we conclude that the action passes to $Q$ by
$$t \circ \underline{p} = (t^{w_0+w_1}p_{01}:t^{w_0+w_2}p_{02}:t^{w_0+w_3}p_{03}:t^{w_1+w_2}p_{12}:t^{w_1+w_3}p_{13}:t^{w_2+w_3}p_{23}).$$
Thus, for distinct weights, the Pl\"ucker relation $p_{01}p_{23}-p_{02}p_{13}+p_{03}p_{12}=0$ guarantees that the only fixed lines by the action are the standard ones $x_i=x_j=0$, $i,j \in \{0,1,2,3\}$.

\subsection{Action on quadrics and cubics}\label{acaoqc}

The set of quadric surfaces of $\mathbb{P}^3$ form a $\mathbb{P}^9$, with induced action from $\xi^\vee$ (\ref{xidual}).
For  $g=b_0x_0^2+b_1x_0x_1+\cdots+b_9x_3^2$, the induced action gives
$$t \circ g = t^{2w_0}b_0x_0^2+ t^{w_0+w_1}b_1x_0x_1+\cdots+t^{2w_3}b_9x_3^2\ .$$
For fixed quadrics to be only the monomials $ x_ix_j = 0 $, it is not enough that the weights are different. Now, it is also necessary that the sums of weights, taken 2--a--2, to be different. 

This is similar to cubics. For $f=a_0x_0^3+a_1x_0^2x_1+\cdots+a_{19}x_3^3$, the induced action is
$$t\circ f=t^{3w_0}a_0x_0^3+t^{2w_0+w_1}a_1x_0^2x_1+\cdots+t^{3w_3}a_{19}x_3^3.$$
Therefore, to have only isolated fixed cubics we need that the sums of weights, taken 3--a--3, to be different. 

As an example, we can take 
\begin{equation}\label{pesospart}
(w_0,w_1,w_2,w_3)=(0,1,5,25)
\end{equation} and the components of fixed flags, quadrics and cubics are all isolated points.

\section{Three planes in $\bb{P}^3$}\label{3planos}

The intention here is to show with a simpler example some methods used in calculations involving Bott's formula.

A basic example is to find the number of points living in the intersection of three generic planes in $\mathbb{P}^3$. The answer to this question is widely known as 1, and can be expressed in terms of the Chow ring of $\mathbb{P}^3$ 
as
$$\int\limits_{\mathbb{P}^3}c_1(\mathcal{O}_{\mathbb{P}^3}(1))^3.$$     

Now allow us to complicate the calculations a little more by computing this number on the blowup of $\mathbb{P}^3$ along a line, let's say, along the line
$$\ell: x_1=x_2=0.$$

Denote by $\pi: X \rightarrow \mathbb{P}^3$ this blowup, and $\mathcal{O}(1)$ the pullback bundle $$\mathcal{O}(1)=\pi^*(\mathcal{O}_{\mathbb{P}^3}(1)).$$

Now, we want to show how to compute $\int\limits_{X}c_1(\mathcal{O}(1))^3 = 1$  applying Bott's formula,
\begin{equation}\label{bottex}
\int\limits_{X}c_1(\mathcal{O}(1))^3 \cap [X] = \sum\limits_{F} \int\limits_{F} \dfrac{c_1^T(\mathcal{O}(1)|_{F})^3 \cap [F]_T}{c_{top}^T(\mathcal{N}_{F|X})}
\end{equation}
where the sum runs through all fixed components $F$ under an action of the torus $\mathbb{C}^*$.

To make things more fun, let's set the action
$$\begin{array}{ccc}
	\mathbb{C}^* \times \mathbb{P}^3 & \longrightarrow & \mathbb{P}^3  \\
\left(t,(s_0:s_1:s_2:s_3)\right) & \longmapsto & (t^{-1}s_0:t^{-1}s_1:t^{-2}s_2:t^{-3}s_3)
\end{array}$$

This action induces an action on the linear forms
just by change the signals because of duality, {\em i.e,} 
\begin{equation} \label{pesosex}
\begin{array}{cccccccc}
x_0 & \mapsto & tx_0 & \ & \ & x_2 & \mapsto & t^2x_2 \\
x_1 & \mapsto & tx_1 & \ & \ & x_3 & \mapsto & t^3x_3 \\
\end{array}
\end{equation}

Notice that the action on $\mathbb{P}^3$ gives two isolated fixed points and one fixed line
$$p_2=(0:0:1:0) \notin \ell;$$
$$p_3=(0:0:0:1) \in \ell;$$ 
$$d=(s_0:s_1:0:0)$$
and the fixed component $d$ (the line of equations $x_2=x_3=0$) cuts the blowup center $\ell$ in one point $p_0=(1:0:0:0)$. 
We have to compute the contributions of the components of fixed points in the sum (\ref{bottex}).

In the numerator of (\ref{bottex}) we have $c_1^T(\mathcal{O}(1)|_{F})$. Since $\mathcal{O}(1)$ is a line bundle, the \linebreak equivariant Chern class of this bundle over a fixed component $F$ is 
\begin{equation} \label{ctex} c_1^T(\mathcal{O}(1))|_{F} = \lambda + c_1(\mathcal{O}(1)|_{F}),
\end{equation}
where $\lambda$ is the weight of the fiber of the bundle on the fixed component $F$ and $c_1(\mathcal{O}(1)|_{F})$ denotes the classical Chern class of $\mathcal{O}(1)$ restricted to $F$.   

%Even more, f
For an equivariant rank $n$ bundle $E$ over $F=\p1
$, %if 
one can split $E$ as a direct sum of $n$ line bundles $E = \oplus _{\lambda} E^\lambda$, for which the weight of the fiber of each $E^\lambda$ is $\lambda$. Hence (\ref{ctex})  applies for each $E^\lambda$, 
$$c_1^T(E^\lambda)|_{F} = \lambda + c_1(E|_{F}).
$$
It follows that the $n$-th equivariant Chern class of $E$ can be obtained by the product
\begin{equation}\label{cntasprod}
c_n^T(E)|_F = \prod _\lambda c_1^T(E^\lambda)|_F.
\end{equation}

$\bullet$ Contribution of $p_2$ \\

The simplest point to compute the contribution in (\ref{bottex}) is $p_2$. Since it is out of the blowup center, its contribution can be computed already in $\mathbb{P}^3$.

The fiber of $\mathcal{O}(1)|_{p_2}$ is generated by $x_2$, of weight $\lambda=2$ (see (\ref{pesosex})). And $\dim (\{p_2\})=0$, so $c_1(\mathcal{O}(1)|_{p_2})=0$. In view of (\ref{ctex}) we have 
\begin{equation} \label{contrp2}
c_1^T(\mathcal{O}(1)|_{p_2})=2.
\end{equation}

The denominator is
$c_3^T(\mathcal{N}_{p_2|\mathbb{P}^3})$. 
%The normal bundle fits in the exact sequence 
%$$T{p_2} \hooklar T_{p_2}\mathbb{P}^3 \twoheadlar \mathcal{N}_{p_2|\mathbb{P}^3}.$$
%But the tangent space to a point is $0$, and then
But
$\mathcal{N}_{p_2|\mathbb{P}^3} \simeq
T_{p_2}\mathbb{P}^3$ as the normal space  of a point in an
ambient is equal to the tangent space  of this ambient at the point. 
To compute $c_3^T(T_{p_2}\mathbb{P}^3)$, remember that
the tangent bundle to any Grassmannian %\linebreak 
$G=G(k,V)$ ($\rk V = r$) 
is (see \cite{3264}, Theorem 2.4, pg. 73)
\begin{equation} \label{QxS}
TG = \textrm{Hom} _G(\mathcal{S},\mathcal{Q})=\mathcal{Q} \otimes  \mathcal{S}^\vee,
\end{equation}
where 
%$$\mathcal{Q}^\vee \hooklar (\mathbb{C}^r)^\vee \twoheadlar \mathcal{S}^\vee $$
$$\mathcal{S} \hooklar (V=\mathbb{C}^r) \twoheadlar 
\mathcal{Q} 
$$
is the tautological sequence over $G$; the fiber of
$\mathcal{Q}^\vee\subset V^\vee$ over a point in $G$ 
corresponding to a $k$-linear subspace of $V$
is the rank $r-k$ space of linear forms spanned by the
equations of this space.
In particular, for $G=\p3,k=1,r=4$, at the point $p_2=(0:0:1:0)$, given by the equations $x_0=x_1=x_3=0$,
the fiber of $\mathcal{Q}$ is spanned by the duals of these three equations; and the fiber of $\mathcal{S}^\vee$ is spanned by $\overline{x_2}:=x_2\mod\id{x_0,x_1,x_3}$. Therefore, using (\ref{QxS}) the tangent space to $\mathbb{P}^3$ at the point $p_2$ can be written as
$$T_{p_2}\mathbb{P}^3 = \langle x_0^\vee, x_1^\vee, x_3^\vee \rangle \otimes \langle x_2 \rangle = \langle x_0^\vee \otimes x_2, \ x_1^\vee \otimes x_2, \ x_3^\vee \otimes x_2 \rangle.$$

The action induced on this tangent space is just 
$$t \ \circ \ (x_0^\vee \otimes x_2) = t^{-w_0}x_0^\vee \otimes t^{w_2}x_2 = t^{w_2-w_0}(x_0^\vee \otimes x_2)$$
$$t \ \circ \ (x_1^\vee \otimes x_2) = t^{-w_1}x_1^\vee \otimes t^{w_2}x_2 = t^{w_2-w_1}(x_1^\vee \otimes x_2)$$
$$t \ \circ \ (x_3^\vee \otimes x_2) = t^{-w_3}x_3^\vee \otimes t^{w_2}x_2 = t^{w_2-w_3}(x_3^\vee \otimes x_2)$$
where $w_i$ is the weight of $x_i$, $i \in \{0,1,2,3\}$
(here, $w_0=w_1=1, \ w_2=2$ and $w_3=3$). 
%Hence, these three vectors are eigenvectors of
%$T_{p_2}\mathbb{P}^3$. 
We will write 
$$T_{p_2}\mathbb{P}^3 = \underbrace{\dfrac{x_2}{x_0}}_{w_2-w_0} + \underbrace{\dfrac{x_2}{x_1}}_{w_2-w_1} + \underbrace{\dfrac{x_2}{x_3}}_{w_2-w_3}.$$

In this notation, the $T_{p_2}\mathbb{P}^3$ is
split into three line subbundles, and each one is
associated to one weight. Using the expression
(\ref{cntasprod}),
we find
$$c_3^T(T_{p_2}\mathbb{P}^3) = (w_2-w_0)(w_2-w_1)(w_2-w_3)=1 \cdot 1 \cdot (-1) = -1.$$  

Therefore, the contribution of the point $p_2$ is
\begin{equation} \label{apbc1}
\dfrac{c_1^T(\mathcal{O}(1)|_{p_2})^3}{c_3^T(T_{p_2}\mathbb{P}^3)} = \dfrac{8}{-1}=-8.
\end{equation}

$\bullet$ Contribution over $p_3$ \\

The point $p_3=(0:0:0:1)$ lives inside the blowup center, the line $\ell$.

The exceptional divisor of the
 blowup $X$ of $\mathbb{P}^3$ along the line $\ell$ 
is a $\mathbb{P}^1$-bundle over the line $\ell$, the projectivized normal bundle $\mathbb{P}(\mathcal{N}_{\ell | \mathbb{P}^3})$. So, over the point $p_3$ there is a $\mathbb{P}^1-$fiber for which we expect two fixed points.

Since $p_3$ is the point of equations $x_0=x_1=x_2=0$, we have
\begin{equation} \label{Tpexe}
T_{p_3}\mathbb{P}^3 = \dfrac{x_3}{x_0} + \dfrac{x_3}{x_1} + \dfrac{x_3}{x_2}.
\end{equation}

The blowup center is $\ell=\{(x_0:0:0:x_3)\} $ (the variables here are only $x_0$ and $x_3$). Over $\ell$ the equation of $p_3$ is only $x_0=0$, and the tangent space to $\ell $ at the point $p_3$ is
\begin{equation} \label{Tlexe}
T_{p_3}\ell = \dfrac{x_3}{x_0}.
\end{equation}
Next, beware of the fact that, although the normal bundle
$\cl N_{\ell|\p3}=\cl O_{\ell}(1)\oplus \cl O_{\ell}(1)$,
these two summands are not isomorphic as $\bb C^\star-$modules.
The short exact sequence
$$T\ell \hooklar T\mathbb{P}^3_{|\ell} \twoheadlar \mathcal{N}_{\ell|\mathbb{P}^3}$$
shows us that
$$(\mathcal{N}_{\ell|\mathbb{P}^3})_{p_3} = \underbrace{\dfrac{x_3}{x_1}}_{w_3-w_1=2} + \underbrace{\dfrac{x_3}{x_2}}_{w_3-w_2=1}.$$
Thus we have two \textit{isolated} fixed points (since
the weights are distinct). We will denote these points by 
$$q_1=\left(p_3,\dfrac{x_3}{x_1}\right)\ \text{ and }\ 
q_2=\left(p_3,\dfrac{x_3}{x_2}\right).
$$
The fibers of $\mathcal{O}(1)$ at these two points are
the same, $\cl O_{\p3}(1)_{p_3}=\id{x_3}$. 
Then, \linebreak $c_1^T(\mathcal{O}(1)|_{q_i}) = 3$, \
$i=1,2$ \ in view of our choice of
weights\, (\ref{pesosex}). 
\\

Now we look at the tangent space to $X$ at these points,
recalling \ $\mathcal{N}_{q_i|X} \simeq T_{q_i}X$. \\

For an inclusion of a smooth scheme $Y$ into a smooth
variety $Z$, denote by $\tilde{Z}$ the blowup of $Z$
along $Y$. The exceptional divisor can be identified with
$\mathbb{P}(\mathcal{N}_{Y|Z})$, and for a point $(y,s)$,
\ $y \in Y$ and $s \in \mathbb{P}(\mathcal{N}_{yY|Z})$,
the tangent space to $\tilde{Z}$ can be split as
\begin{equation}\label{tblowup}
T_{(y,s)}\tilde{Z} = T_yY \oplus \mathcal{L}_s \oplus T_{[\mathcal{L}_s]}\mathbb{P}(\mathcal{N}_{yY|Z}),
\end{equation}
where $\mathcal{L}_s \subset \mathcal{N}_{yY|Z}$ is the line represented by the point $[\mathcal{L}_s] \in \mathbb{P}(T_yZ)$.\\

For the point $q_1$, (\ref{tblowup}) gives
$$T_{q_1}X = T_{p_3}\ell \oplus \mathcal{L}_{\frac{x_3}{x_1}} \oplus T_{[\mathcal{L}_{\frac{x_3}{x_1}}]}\mathbb{P}(\mathcal{N}_{p_3\ell|\mathbb{P}^3}).$$

But $\mathcal{N}_{p_3\ell|\mathbb{P}^3}=
\dfrac{x_3}{x_1}
+ \dfrac{x_3}{x_2} $ and so 
$$T_{[\mathcal{L}_{\frac{x_3}{x_1}}]}\mathbb{P}(\mathcal{N}_{p_3\ell|\mathbb{P}^3}) = \underbrace{\dfrac{x_3}{x_2}}_{\mathcal{S}^\vee} \cdot \underbrace{\dfrac{x_1}{x_3}}_{\mathcal{Q}} = \dfrac{x_1}{x_2}.$$
This yields 
$$T_{q_1}X = \dfrac{x_3}{x_0} + \dfrac{x_3}{x_1} + \dfrac{x_1}{x_2},$$
and $$c_3^T(T_{q_1}X) = (w_3-w_0)(w_3-w_1)(w_1-w_2)=2 \cdot 2 \cdot (-1) = -4.$$
With this, we have the contribution of the fixed point $q_1:$
\begin{equation} \label{apbc2}
\dfrac{c_1^T(\mathcal{O}(1)|_{q_1})^3}{c_3^T(T_{q_1}X)} = \dfrac{27}{-4} \ .
\end{equation}
In a similar way, $$T_{q_2}X = \dfrac{x_3}{x_0} + \dfrac{x_3}{x_2} + \dfrac{x_2}{x_1},$$ 
so the contribution of $q_2$ is
\begin{equation} \label{apbc3}
\dfrac{c_1^T(\mathcal{O}(1)|_{q_2})^3}{c_3^T(T_{q_2}X)} = \dfrac{27}{2} \ .
\end{equation}  \\

$\bullet$ Contribution over the fixed line $d$ \\

The fixed component $d$ meets the blowup center $\ell$
in one point $p_0=(1:0:0:0)$. Then its contribution to the sum (\ref{bottex}) cannot be computed yet.

Like (\ref{Tpexe}) and (\ref{Tlexe}),

\begin{equation} \label{tp0ex}
T_{p_0}\mathbb{P}^3 = \dfrac{x_0}{x_1} + \dfrac{x_0}{x_2} + \dfrac{x_0}{x_3},
\end{equation}

$$T_{p_0}\ell = \dfrac{x_0}{x_3}.$$

Hence, $$\mathcal{N}_{p_0\ell|\mathbb{P}^3} = \underbrace{\dfrac{x_0}{x_1}}_{w_0-w_1=0} + \underbrace{\dfrac{x_0}{x_2}}_{w_0-w_2=-1} = 1 + \dfrac{x_0}{x_2},$$
where the notation 1 represents an eigenvector associated
to the weight 0 (due to the weights
$w_0=w_1$). Therefore, over the point $p_0$ there are two
isolated fixed points, namely $r_1=(p_0,1)$ and
$r_2=\left(p_0,\dfrac{x_0}{x_2}\right)$. It is over each
of these points that we will compute the contributions to the sum (\ref{bottex}). \\

The contribution of $r_2$ follows the same ideas like the points $q_1$ and $q_2$. Write (\ref{tblowup})
$$T_{r_2}X = \underbrace{\dfrac{x_0}{x_3}}_{T_{p_0}\ell}
\ \ +\ \ 
\underbrace{\dfrac{x_0}{x_2}}_{\mathcal{L}_{\frac{x_0}{x_2}}}
\ \ \ +
\underbrace{\dfrac{x_2}{x_0}}_{T_{[\mathcal{L}_{\frac{x_0}{x_2}}]}\mathbb{P}(\mathcal{N}_{p_0\ell|\mathbb{P}^3})};$$
so
$$c_3^T(T_{r_2}X) = (w_0-w_3)(w_0-w_2)(w_2-w_0)=(-2)\cdot(-1)\cdot 1 = 2.$$ 

The fiber of $\mathcal{O}(1)|_{r_2}$ is spanned by $x_0$, so $c_1^T(\mathcal{O}(1)|_{r_2})=1$. Then,

\begin{equation} \label{apbc4}
\dfrac{c_1^T(\mathcal{O}(1)|_{r_2})^3}{c_3^T(T_{r_2}X)} = \dfrac{1}{2} \ .
\end{equation} \\

The last point is $r_1=(p_0,1)$. But this point lives in
the strict transform $\tilde{d}$ of the line $d$, which
constitutes a full component of fixed points and we will use this $\tilde{d}$ to compute the contribution in (\ref{bottex}).   

The bundle $\mathcal{O}(1)|_{\tilde{d}} = \mathcal{O}_{\tilde{d}}(1)$ has weight 1 since both $x_0$ and $x_1$ have this weight. But $\tilde{d}$ has dimension one, and $c_1(\mathcal{O}_{\tilde{d}}(1)) = h$. In view of (\ref{ctex}),
$$c_1^T(\mathcal{O}_{\tilde{d}}(1)) = 1+h.$$

The denominator in (\ref{bottex}) is now $c_2^T(\mathcal{N}_{\tilde{d}|X}).$ In terms of weights, we can write over the point $r_1$
$$T_{r_1}X = \dfrac{x_0}{x_3} + 1 +\dfrac{x_0}{x_2}.$$
From $T_{r_1}\tilde{d}=\dfrac{x_0}{x_1}=1$, the normal space can by written as
\begin{equation}\label{normalexe}
\mathcal{N}_{r_1\tilde{d}|X} = \dfrac{x_0}{x_3} + \dfrac{x_0}{x_2}. 
\end{equation}

Notice that (\ref{normalexe}) gives a decomposition of $\mathcal{N}_{r_1\tilde{d}|X}$ in two eigensubspaces, and we can compute $c_2^T(\mathcal{N}_{\tilde{d}|X})$ like (\ref{cntasprod}).

Remember that (\ref{tp0ex})
$$T_{p_0}\mathbb{P}^3 = \dfrac{x_0}{x_1} + \dfrac{x_0}{x_2} + \dfrac{x_0}{x_3}$$
and
$$T_{p_0}d = \dfrac{x_0}{x_1} = 1.$$
Then, $$\mathcal{N}_{p_0d|X} = \dfrac{x_0}{x_3} + \dfrac{x_0}{x_2}.$$
(Note the same description as (\ref{normalexe})). But
this is the fiber of the normal bundle of a $\mathbb{P}^1$ in
$\mathbb{P}^3$, which can be described as $\mathcal{N}_{d|\mathbb{P}^3} =
\mathcal{O}_{d}(1) \oplus \mathcal{O}_{d}(1)$.\\

After the blowup, the normal bundle of $\tilde{d}$ in $X$ shall be modified, 
in accordance with the Proposition
\ref{proptensor}, page \pageref{proptensor}. It
says that the subbundles related to the blowup
center are unchanged, and the other subbundles
shall be tensorized with
$\mathcal{O}_{\tilde{d}}(-1)$. From $T_{p_0}\ell =
\dfrac{x_0}{x_3}$ it follows that \\

$\bullet $ the subbundle $\mathcal{O}_d(1)$ related to the weight $\lambda_1 = \dfrac{x_0}{x_3}$ will not change after the blowup: 
$$\mathcal{N}^{\lambda_1}_{\tilde{d}|X} = \mathcal{O}_{\tilde{d}}(1).$$

$\bullet $ the subbundle $\mathcal{O}_d(1)$ related to the weight $\lambda_2 = \dfrac{x_0}{x_2}$ will change after the blowup: 
	$$\mathcal{N}^{\lambda_2}_{\tilde{d}|X} = \mathcal{O}_{\tilde{d}}(1) \otimes \mathcal{O}_{\tilde{d}}(-1) = \mathcal{O}_{\tilde{d}}.$$

This will be denoted by
$$\mathcal{N}_{r_1\tilde{d}|X} =
\underbrace{\dfrac{x_0}{x_3}}_{\mathcal{O}_{\tilde{d}}(1)}
+ \underbrace{\dfrac{x_0}{x_2}}_{\mathcal{O}_{\tilde{d}}}.$$

Now we are able to use (\ref{cntasprod})

$$c_2^T(\mathcal{N}_{r_1\tilde{d}|X}) = (w_0-w_3+h)(w_0-w_2) = 2 - h,$$
and so
$$\int\limits_{\tilde{d}}\dfrac{c_1^T(\mathcal{O}_{\tilde{d}}(1))^3}{c_2^T(\mathcal{N}_{\tilde{d}|X})} = \int\limits_{\tilde{d}}\dfrac{(1+h)^3}{2-h} = \int\limits_{\tilde{d}}\dfrac{(1+h)^3(2+h)}{4}$$ 
where in the last equality we multiplied both numerator and denominator by $2+h$, and used that $h^2=0$ (of course, dim $\tilde{d} = 1$).

Expanding the numerator and taking the coefficient of $h$ we obtain the contribution 
\begin{equation}\label{apbc5}
\int\limits_{\tilde{d}}\dfrac{c_1^T(\mathcal{O}_{\tilde{d}}(1))^3}{c_2^T(\mathcal{N}_{\tilde{d}|X})} = \dfrac{7}{4}.
\end{equation}

Putting together the contributions (\ref{apbc1}), (\ref{apbc2}), (\ref{apbc3}), (\ref{apbc4}) and (\ref{apbc5}) the number in (\ref{bottex}) is

$$-8-\dfrac{27}{4}+\dfrac{27}{2}+\dfrac{1}{2}+\dfrac{7}{4} = 1.$$

\section{Normal bundle of strict transform }\label{normalblowup}

Now we will show the following:  \medskip

\begin{prop}\label{proptensor}
Let  $\bb{P}^1 = \ell \subset \bb{P}^N$ and $Y \subset \bb{P}^N$ be a smooth $m$--dimensional subscheme with $\ell \cap Y = \{p\}$, a single point. Let $X=Bl_Y\bb{P}^N$, $\tilde \ell$ the strict transform of $\ell$, and 
$$\mathcal{N}_{\ell|\bb{P}^N} = \mathcal{O}_{\ell}(a_1) \oplus \cdots \oplus \mathcal{O}_{\ell}(a_m) \oplus \mathcal{O}_{\ell}(a_{m+1}) \oplus \cdots \oplus \mathcal{O}_{\ell}(a_{N - 1}).$$

If, for any $\bb{C}^*$--action that fix $\ell$ we have $Y$ invariant and the tangent space to $Y$ at $p$ is isomorphic as $\mathbb{C}^*$--module with the fibers
$$T_pY \simeq \mathcal{O}_{\ell}(a_1)|_p \oplus \cdots \oplus \mathcal{O}_{\ell}(a_m)|_p\ ,$$
then we have the decomposition for $\mathcal{N}_{\tilde \ell|X} $,
$$\mathcal{N}_{\tilde \ell|X} = \mathcal{O}_{\tilde \ell}(a_1) \oplus \cdots \oplus \mathcal{O}_{\tilde \ell}(a_m) \oplus \mathcal{O}_{\tilde \ell}(a_{m+1}-1) \oplus \cdots \oplus \mathcal{O}_{\tilde \ell}(a_{N - 1}-1).$$

\begin{proof}
Let $\bb{P}^N=\{(x_0:\ldots:x_N)\}$ and $\ell=\bb{P}^1=\{(x_0:x_1:0:\ldots:0)\}$, $p=(1:0:\ldots:0)$. We will compute $\int\limits_{\bb{P}^N}c_1[\cl{O}_{\bb{P}^N}(1)]^N$ via Bott in two ways: directly in $\bb{P}^N$ and in the blowup along $Y$.

Fix $v \in \mathbb{Z}$, $v \geq 2$, and put an action on linear forms by
$$\begin{array}{cccccccccccl}
	x_0 & \mapsto & tx_0 & \ & \ & x_2 & \mapsto & t^vx_2 & \ & x_i & \mapsto & t^{v+i-2}x_i \ (2\leq i \leq N)\\
	x_1 & \mapsto & tx_1 & \ & \ & x_3 & \mapsto & t^{v+1}x_3 & \ & x_N & \mapsto & t^{v+N-2}x_N \\
\end{array}$$

Thus, the weights of $x_0$ and $x_1$ are equal to 1 and the weight of $x_i$, $i \geq 2$, is $v+i-2$. 

The line $\ell$ is a fixed component, and %since $\cl{O}_{\bb{P}^N}(1)$ is well defined along $\ell$ 
we can compute the contribution of this line now.

We have $c_1^T(\cl{O}_{\bb{P}^N}(1))$ = 1 (weight of $x_0,x_1$) + $h$ ($=c_1(\cl{O}_\ell(1))$). And from	
$$\mathcal{N}_{p\ \ell|\bb{P}^N} = \underbracket{\dfrac{x_0}{x_2}}_{\mathcal{O}_{\ell}(a_1)} + \underbracket{\dfrac{x_0}{x_3}}_{\mathcal{O}_{\ell}(a_2)} + \cdots +  \underbracket{\dfrac{x_0}{x_N}}_{\mathcal{O}_{\ell}(a_{(N-1)})},$$
$$c_{N-1}^T(\cl{N}_{\ell|\bb{P}^N}) = [c_1(\cl{O}_\ell(a_1))+1-v]\cdot [c_1(\cl{O}_\ell(a_2))+1-v-1]\cdot \ldots \cdot [c_1(\cl{O}_\ell(a_{N-1}))+1-v-N+2]$$
$$c_{N-1}^T(\cl{N}_{\ell|\bb{P}^N}) = [a_1h+1-v]\cdot [a_2h-v]\cdot \ldots \cdot [a_{N-1}h-v-N+3]$$
$$c_{N-1}^T(\cl{N}_{\ell|\bb{P}^N}) = Z + \left(\dfrac{a_1Z}{1-v} + \dfrac{a_2Z}{-v} + \cdots + \dfrac{a_{N-1}Z}{-v-N+3} \right)h$$
$$c_{N-1}^T(\cl{N}_{\ell|\bb{P}^N}) = Z\left(1 + \sum_{i=1}^{N-1} \dfrac{a_i}{2-i-v}h \right), $$
where
$$Z=(1-v)(-v)(-v-1)\cdots (-v-N+3) =
(-1)^{N-1}\dfrac{(v+N-3)!}{(v-2)!}\ .
$$

The contribution of $\ell$ is
$$\int\limits_{\ell}\dfrac{c_1^T(\cl{O}_\ell(1))^N}{c_{N-1}^T(\cl{N}_{\ell|\bb{P}^N})} = \int\limits_{\ell}\dfrac{1+Nh}{Z\cdot (1+\sum \frac{a_i}{2-i-v}h)}= \frac{1}{Z}\int\limits_{\ell}(1+Nh)(1-\sum_{i=1}^{N-1} \frac{a_i}{2-i-v}h)$$
\begin{equation}\label{SA}
\int\limits_{\ell}\dfrac{c_1^T(\cl{O}_\ell(1))^N}{c_{N-1}^T(\cl{N}_{\ell|\bb{P}^N})} = \dfrac{1}{Z}\left(N - \sum_{i=1}^{N-1}\dfrac{a_i}{2-i-v}  \right)\ .
\end{equation}

This contribution is the same if computed on the fixed points that arise over the point $p$ when we do the blowup along $Y$.

By the hypotheses of tangent space to $Y$ at $p$, this tangent space is ($1\leq m \leq N-2$)
$$T_pY=\dfrac{x_0}{x_2} + \dfrac{x_0}{x_3} + \cdots + \dfrac{x_0}{x_{m+1}}\ .$$
So, $\cl{N}_{p\ Y|\bb{P}^N} = 1 + \dfrac{x_0}{x_{m+2}} + \dfrac{x_0}{x_{m+3}} + \cdots + \dfrac{x_0}{x_N}\ .$ There are $N-m-1$ isolated fixed points, corresponding to $\dfrac{x_0}{x_i}$, and a point $\tilde p \in \tilde \ell$.

To get the contribution of $q_i:=(p,\frac{x_0}{x_i})$,\ $m+2\leq i \leq N$, note that
$$T_{q_i}X=\underbrace{\dfrac{x_0}{x_2}+\cdots + \dfrac{x_0}{x_{m+1}}}_{T_pY} + \underbrace{\dfrac{x_0}{x_i}}_{\cl{L}_{x_0/x_i}}+ \underbrace{\dfrac{x_i}{x_0} + \dfrac{x_i}{x_{m+2}}+\cdots + \overbracket{\hat{\dfrac{x_i}{x_i}}}^{\textrm{omit}} + \cdots + \dfrac{x_i}{x_N}}_{T_{\cl{L}_{x_i/x_0}}\bb{P}(\cl{N}_{p \ Y|\bb{P}^N})}$$
$$c_N^T(T_{q_i}X)=(1-v)(-v)(-v-1)\cdots(-v-m+2) \cdot (-1)(v+i-3)^2 \cdot (i-m-2)(i-m-3)\cdots(i-N)$$
$$c_N^T(T_{q_i}X)=(-1)^{N-1}\dfrac{(v+m-2)!}{(v-2)!}\cdot (v+i-3)^2 \cdot \prod_{\substack{j=m+2 \\ j \not= i}}^{N}(j-i)$$
Then, the contribution of the point $q_i$ is
\begin{equation}\label{ggi}
\dfrac{(-1)^{N-1}(v-2)!}{(v+m-2)!\cdot (v+i-3)^2 \cdot \prod\limits_{\substack{j=m+2 \\ j \not= i}}^{N}(j-i)}
\end{equation}

Next, we compute the contribution of $\tilde \ell$. Write 
$$\mathcal{N}_{p\ \tilde \ell|X} = \underbracket{\dfrac{x_0}{x_2}}_{\mathcal{O}_{\ell}(b_1)} + \underbracket{\dfrac{x_0}{x_3}}_{\mathcal{O}_{\ell}(b_2)} + \cdots +  \underbracket{\dfrac{x_0}{x_N}}_{\mathcal{O}_{\ell}(b_{(N-1)})},$$
for unknowns $b_i$. The computational aspects are the same as for $\ell$, just changing $\underline{a}$ by $\underline{b}$, so the contribution of $\tilde \ell$ is
\begin{equation}\label{SB}
\dfrac{1}{Z}\left(N - \sum_{i=1}^{N-1}\dfrac{b_i}{2-i-v}  \right)\ .
\end{equation}

The sum of contributions in (\ref{ggi}) and (\ref{SB}) is equal to (\ref{SA}), and so we get the equation
$$\dfrac{1}{Z}\left(N - \sum_{i=1}^{N-1}\dfrac{a_i}{2-i-v}  \right) = \dfrac{1}{Z}\left(N - \sum_{i=1}^{N-1}\dfrac{b_i}{2-i-v}  \right) + \sum\limits_{i=m+2}^{N}	\dfrac{(-1)^{N-1}(v-2)!}{(v+m-2)! (v+i-3)^2  \prod\limits_{\substack{j=m+2 \\ j \not= i}}^{N}(j-i)}$$
$$\Downarrow \ \textrm{ multiply everything by }Z$$
\begin{equation}\label{eqab}
\sum_{i=1}^{N-1}\dfrac{a_i}{2-i-v} = \sum_{i=1}^{N-1}\dfrac{b_i}{2-i-v} - \dfrac{(v+N-3)!}{(v+m-2)!}\sum\limits_{i=m+2}^{N}	\dfrac{1}{(v+i-3)^2 \cdot \prod\limits_{\substack{j=m+2 \\ j \not= i}}^{N}(j-i)} \ .
\end{equation}
%$$\Downarrow$$
%\begin{equation}\label{eqab}
%\sum_{i=1}^{N-1}\dfrac{a_i}{2-i-v} = \sum_{i=1}^{N-1}\dfrac{b_i}{2-i-v} + %\sum_{i=m+1}^{N-1}\dfrac{1}{2-i-v}
%\end{equation}

Now, we can verify that the solutions 
\begin{equation}\label{solab}
b_1=a_1, \ b_2 = a_2, \ldots , b_m=a_m, \ b_{m+1}=a_{m+1}-1, \ldots, b_{N-1}=a_{N-1}-1.
\end{equation}
holds for any $v \in \bb{Z}, \ v \geq 2$.

Setting the $N-1$ values $v=2, \ v=3, \ldots , v=N$ on equation (\ref{eqab}) we obtain the linear system $CB-CA=D$, where
$$A=[a_1 \cdots a_{N-1}]^t, \ B=[b_1 \cdots b_{N-1}]^t, \ C=
\left[\begin{array}{cccc}
	1           & \frac{1}{2} & \cdots & \frac{1}{N-1}\\
	\frac{1}{2} & \frac{1}{3} & \cdots & \frac{1}{N}\\
	\vdots      &   \vdots    & \ddots & \vdots     \\
	\frac{1}{N-1} & \frac{1}{N} & \cdots & \frac{1}{2N-3}
\end{array}\right]_{(N-1)\times(N-1)}$$
Since det($C)\not=0$, the solution (\ref{solab}) is the unique solution for $b's$ with respect to $a's$. 
\end{proof}
\end{prop}

\begin{exem}\normalfont
For $N=4$ and $m=1$ ({\em i.e,} blowup of $\bb{P}^4$ along a smooth curve), equation (\ref{eqab}) reduces to
$$\dfrac{a_1}{1-v}+\dfrac{a_2}{-v}+\dfrac{a_3}{-1-v}=\dfrac{b_1}{1-v}+\dfrac{b_2}{-v}+\dfrac{b_3}{-1-v}-\dfrac{2v+1}{v(v+1)}\ .$$
 
 It is a straightforth calculation the verification of the solution, for all $v \in \mathbb{Z}, \ v \geq 2$,
 \begin{equation}\label{abn4m1}
 b_1=a_1, \ \ b_2=a_2-1, \ \ b_3=a_3-1.
 \end{equation}
   
   Moreover, setting values 2, 3 and 4 to $v$ we obtain the linear system
   $$
   \begin{array}{ccc} \label{eqn4m1}
   	v=2 & \Rightarrow & 6a_1+3a_2+2a_3=6b_1+3b_2+2b_3+5 \\
   	v=3 & \Rightarrow & 6a_1+4a_2+3a_3=6b_1+4b_2+3b_3+7 \\
   	v=4 & \Rightarrow & 20a_1+15a_2+12a_3=20b_1+15b_2+12b_3+27 
   \end{array}$$
   whose solutions are unique (\ref{abn4m1}).
  \qed
\end{exem}

Next, we present scripts for \textit{Macaulay2} to perform the calculations in the above proof.

First, set values for $N,m$ as you desire, provided that $2 \leq m \leq N-2$.
\begin{verbatim}
N=4; M=1; -- 1<M<N-1
-- change these values if you want to.

R=QQ[x_0..x_N,v,a_1..b_N];
Rh=R[h]/h^2;

-- setting the weights
p01={x_0=>1,x_1=>1};
poth=apply(N-1,i->x_(i+2)=>v+i);
pesos=flatten(append(p01,poth));

-- tangent and normal spaces to l and Y(=d here)
TPN=apply(N,i->x_0-x_(i+1));
Tl=x_0-x_1; Nl=delete(Tl,TPN);
Td=apply(M,i->x_0-x_(2+i)); Nd=TPN-set(Td);

-- computing the contribution contl of line l
cn=product(apply(N-1,i->a_(i+1)*h+sub(Nl_i,pesos)));
conj=cn-2*(cn%h); den=cn*conj;
num=((1+h)^N)*conj;
num=num//h;
den=sub(den,R);
num=sub(num,R);
contl=num/den;

-- after the blowup

-- contributions of points x_i/x_0
for i from 0 to N-M-1 do (
eq=Nd_i;
if eq!=x_0-x_1 then (
Teq=delete(eq,Nd);
Teq=apply(N-M-1,j->Teq_j-eq);
Teq=append(Teq,eq);
Teq=flatten(append(Teq,Td));
Teq=apply(N,j->sub(Teq_j,pesos));
Te_i=product Teq; ) );

-- the contribution contL of strict transform ~l
cnL=product(apply(N-1,i->b_(1+i)*h+sub(Nl_i,pesos)));
conjL=cnL-2*(cnL%h); denL=cnL*conjL;
numL=((1+h)^N)*conjL;
numL=numL//h;
denL=sub(denL,R);
numL=sub(numL,R);
contL=numL/denL;

-- sum of contributions
eq=0;
for i from 0 to N-M-1 do (
if Nd_i!=x_0-x_1 then (
eq=eq+(1/Te_i); ) );

-- matching the contributions before and after the blowup
-- eq represents the obtained equation
eq=eq+contL-contl;
eq=numerator eq;

-- solutions of linear system and the equations like in exemple
I=ideal();
for i from 0 to N-2 do (
eqt_i=sub(eq,v=>i+2);
I=I+eqt_i;  
print(i+2, eqt_i);   )
trim I
\end{verbatim}

\begin{remk}\normalfont
For $m=0$, {\em i.e,} blowup of $X$ along a point, a similar result of Proposition \ref{proptensor} can be found in \cite{Ionescu}, Lemma 1.6, page 5.
\end{remk}

%\nocite{israel1}
%\nocite{3264}
%\nocite{atiyah1994introduction}
%\nocite{meirelles}
%\nocite{Edidin96equivariantintersection}
%\nocite{Fulton}
%\nocite{GomezMeNeto}
%\nocite{Harris}
%\nocite{IsraelTI}
%\nocite{Jouanolou}
%\nocite{Wanderson}
%\nocite{NetoeEndonvein}
%\nocite{edidin1998localization}
%\nocite{eisenbud2013commutative}
%\nocite{israelevivi}
%\nocite{hartshorne2013algebraic}
%\nocite{viviane}
%\nocite{CukAceaMassri}
%\nocite{AluffiSegre}
%\nocite{israel1}
%\nocite{Helmer}

\bibliographystyle{abntex2-num}
\bibliography{referencias}
\end{document}